%% file: piece6_optimale_JB_arXivX6.tex
\newcommand{\cras}{0}

\newcommand{\definitE}{%
l'ensemble $\mathcal{E}$ des courbes à courbure algébrique positive, dont les extrémités et les tangentes en leurs extrémités sont données%
}
\newcommand{\definitEe}{%
the set $\mathcal{E}$ of curves with positive algebraic curvature, whose extremities and tangents in their extremities are given%
}

\newcommand{\resumefrancais}{%
On considère \definitE. \`A chacune des courbes de $\mathcal{E}$, on associe le minimum du rayon de courbure algébrique. 
Il existe une unique courbe de $\mathcal{E}$ qui maximise ce minimum et cette courbe
est égale  à l'unique courbe de $\mathcal{E}$, formée d'un arc de cercle et d'un segment de droite, éventuellement réduit à un point. 
Cette courbe correspond aussi à un cas particulier des courbes de Dubins et sera 
utilisée pour améliorer la conception d'une pièce intervenant  dans un brevet.
}
\newcommand{\resumeanglais}{%
We consider \definitEe. 
For each of the curves of $\mathcal{E}$, we define the minimum of the radius of curvature. 
There exists a unique curve of  $\mathcal{E}$ which maximizes this minimum and this curve
is equal  to the unique curve of $\mathcal{E}$ composed  of an  arc of circle 
and a line segment, 
where appropriate reduced to a point. This curve corresponds also to a particular case of Dubins's curve and will be used to improve the conception of a piece of a patent.
}

\newcommand{\adresse}{%
Laboratoire Inter-universitaire de Biologie de la Motricité\\
   POLYTECH\\
   Université Claude Bernard - Lyon 1\\
   15 Boulevard André LATARJET\\
   69622 Villeurbanne Cedex\\
France
}

\newcommand{\nom}{Jérôme Bastien}

\ifcase \cras
\newcommand{\ftitre}{\MakeUppercase{Existence et unicité d'une courbe à courbure positive maximisant le minimum du rayon de courbure}}
\or
\newcommand{\ftitre}{Existence et unicité d'une courbe à courbure positive maximisant le minimum du rayon de courbure}
\fi
\newcommand{\atitre}{Existence and uniqueness of a curve with positive curvature maximizing the minimum radius of curvature}
\newcommand{\ftitrefacul}{{\MakeUppercase{Existence et unicité d'une courbe maximisant le minimum du rayon de courbure}}}
                                          
\newcommand{\aeL}{\text{p.p. sur  } ]0,L[}

\newcommand{\numerosectioncasdiscretgeneral}{5.1}

\ifcase \cras

\documentclass[a4paper,11pt,reqno,final]{amsart}
\newcommand{\pathfilesty}{.}
\usepackage{\pathfilesty/articlefrancais}
\graphicspath{{./dessins/}}

\bibliography{%
piece6_optimale_JB_arXiv[X6]_new,%
../../../divers/circuit_rail/piece6_optimale/description_piece6_optimale,%
../../../divers/circuit_rail/piece6_optimale/description_piece6_optimale_maths,%
../../../divers/rapport_activite/mabibliographiebis,%
../../../divers/rapport_activite/monenseignementbib%
}

\begin{document}
\selectlanguage{french}



\title[\ftitrefacul]{\ftitre}

\author{\nom}

\date{\today}

\address{\adresse}

\email{\href{mailto:jerome.bastien@univ-lyon1.fr}{\nolinkurl{jerome.bastien@univ-lyon1.fr}}}

\begin{abstract}
\resumefrancais
\end{abstract}

\begin{abstract}
\selectlanguage{french}
\resumefrancais
\vspace{0.5 cm}

\noindent
\textsc{Abstract}.
\selectlanguage{english}
{\textit{ (\atitre ) }}
\resumeanglais
\selectlanguage{french}
\end{abstract}


\maketitle

\newcommand{\sout}[1]{} 
\newcommand{\uwave}[1]{#1} 
\newcommand{\rar}{0}
\input{ajouts}


\or


%

\documentclass{elsart3-1}


\usepackage{epsfig}

\usepackage{amssymb}

\usepackage[french,english]{babel}

\newtheorem{theorem}{Theorem}[section]
\newtheorem{lemma}[theorem]{Lemma}
\newtheorem{e-proposition}[theorem]{Proposition}
\newtheorem{corollary}[theorem]{Corollary}
\newtheorem{e-definition}[theorem]{Definition\rm}
\newtheorem{remark}{\it Remark\/}
\newtheorem{example}{\it Example\/}
\newtheorem{theoreme}{Th\'eor\`eme}[section]
\newtheorem{lemme}[theoreme]{Lemme}
\newtheorem{proposition}[theoreme]{Proposition}
\newtheorem{corollaire}[theoreme]{Corollaire}
\newtheorem{definition}[theoreme]{D\'efinition\rm}
\newtheorem{remarque}{\it Remarque}
\newtheorem{exemple}{\it Exemple\/}
\renewcommand{\theequation}{\arabic{equation}}
\setcounter{equation}{0}

\usepackage[T1]{fontenc}
\usepackage[latin1]{inputenc}
\usepackage{url}
\usepackage[centerlast,small]{subfigure}
\usepackage{amsmath}
\usepackage{amsfonts}
\usepackage{psfrag}
\newenvironment{preuve}[1][]{\textsc{Preuve#1}}{\hfill \null \hfill  \qed}
\newenvironment{preuveidee}[1][]{\textsc{Idées de la preuve#1}}{\hfill \null \hfill \qed}
\newcommand{\Er}{\mathbb{R}}
\newcommand{\Erp}{\mathbb{R_{\, +}}}
\newcommand{\Erm}{\mathbb{R_{\, -}}}
\newcommand{\Erps}{\mathbb{R_{\, +} ^{\, *}}}
\newcommand{\Ce}{\mathbb{C}}
\newcommand{\En}{\mathbb{N}}
\newcommand{\norm}[2][{}]{\lVert#2\rVert_{{#1}}}
\newcommand{\vnorm}[2][{}]{\left\Vert#2\right\Vert_{#1}}
\newcommand{\prodsca}[2]{\left\langle #1,#2\right\rangle}  

\def\og{\leavevmode\raise.3ex\hbox{$\scriptscriptstyle\langle\!\langle$~}}
\def\fg{\leavevmode\raise.3ex\hbox{~$\!\scriptscriptstyle\,\rangle\!\rangle$}}

\journal{the Acad\'emie des sciences}
\begin{document}
\centerline{Géométrie/ Geometry}
\centerline{Problèmes mathématiques de la mécanique/Mathematical Problems in Mechanics}
\begin{frontmatter}

\selectlanguage{french}
\title{\ftitre}

\author{\nom}
\address{\adresse}

\ead{jerome.bastien@univ-lyon1.fr}

\medskip
\selectlanguage{french}

\begin{abstract}
\selectlanguage{french}
\resumefrancais
\vskip 0.5\baselineskip

\selectlanguage{english}
\noindent{\bf Abstract}
\vskip 0.5\baselineskip
\noindent
{\bf  \atitre}
\resumeanglais
\end{abstract}
\end{frontmatter}

\graphicspath{{./dessins/}}


\fi


\selectlanguage{english}
\section*{Abridged English version}

In the framework of the patent \cite{brevetJB}, we had to define a curve 
whose  extremities $A$  et $B$  and the tangents at its points $A$ and $B$ are given,
both these  tangents  not being parallel.
The chosen curve is a parabola, or equivalently, a Bézier curve
(see also  
\cite{enumeration_circuit_JB_2016,%
bastienbrevetforum,%
bastienbrevetmmi,%
bastienbrevetapmep}). 
The disadvantage of this curve is that it has too small radii of curvature and we attempted  to find a  less incurved curve.
For this, we define \definitEe.
For each of curves of $\mathcal{E}$, we define the minimum of the radius of curvature.

This problem is very close to the problem of Dubins's curves \cite{MR0089457,MR0132460}, but it is not equivalent.
Dubins  also looks for curves whose extremities and tangents at the extremities are given. Under the assumption that the radius of curvature is everywhere
on the curve greater than a given radius of curvature $R>0$, he determines the curve which minimizes the length.
He proves that these curves, entirely defined by $R$, called geodesic and denoted
\begin{equation}
\label{rapeldubinseq01}
\mathcal{G}(R),
\end{equation}
are necessarily the union of three arcs of circle of radius $R$,
or the union of two arcs of circle of radius $R$ and of line segment or  a subset of these curves.
In our case, 
we impose the 
\ajoutB
algebraic curvature and 
we do not consider \textit{a priori} this radius 
of curvature $R$.
We 
prove that there exists a unique curve of $\mathcal{E}$ composed of an arc of circle of radius $R_a$ and of a line segment, denoted $\mathcal{J}$.
This case corresponds to the limit case of figure \ref{coducl1}.
This radius $R_a$ depends only on the points $A$ and $B$ and on the tangents to these points and 
 is the greatest value of $R$, for which the Dubins's curve 
$\mathcal{G}(R)$ belongs to  $\mathcal{E}$.
Next, we prove  that the curve of $\mathcal{E}$, maximizing the minimum of the radius of curvature is unique and is precisely equal to  $\mathcal{J}$.

\selectlanguage{french}


\ifcase \cras

\section*{Version abrégée en français}

Dans le cadre du brevet  \cite{brevetJB}, il a été nécessaire de construire une courbe
passant par deux points $A$ et $B$ du plan, dont les directions des tangentes en $A$  et $B$ 
sont données en étant non parallèles. La courbe choisie est une parabole, ou de façon équivalente, une courbe de Bézier 
(voir aussi 
\cite{enumeration_circuit_JB_2016,%
bastienbrevetforum,%
bastienbrevetmmi,%
bastienbrevetapmep}). 
Cette courbe présentant l'inconvénient d'avoir des rayons de courbures trop petits, on a essayé de trouver une courbe moins incurvée. Pour cela, 
on définit \definitE\ et à chacune des courbes de $\mathcal{E}$, on associe le minimum du rayon de courbure.

Ce problème ressemble fortement aux courbes de Dubins \cite{MR0089457,MR0132460}
sans lui être équivalent.
Dubins cherche des courbes  passant aussi par deux points $A$ et $B$ du plan, dont les directions des tangentes en $A$ en $B$ 
sont données. Sous l'hypothèse qu'en tout point, le rayon de courbure est supérieur à $R$ où $R>0$ est donné à l'avance, il cherche la courbe qui minimise la longueur.
Il montre que de telles courbes, entièrement définie par $R$ appelées  géodésiques et notées \eqref{rapeldubinseq01},
ne peuvent être que la réunion de trois arcs de cercles de rayon $R$
ou la réunion de deux  arcs de cercle de rayon $R$ et d'un segment, ou une sous partie de ces deux courbes.
Dans notre cas, 
nous imposons la courbure algébrique positive,
nous ne donnons pas ce rayon de courbure $R$ \textit{a priori}.
Nous montrons 
tout d'abord 
qu'il existe une unique courbe de $\mathcal{E}$ formée d'un arc de cercle de rayon $R_a$ et d'un segment de droite,
notée $\mathcal{J}$.
Ce cas   correspond au cas limite de la figure \ref{coducl1}.
Ce rayon $R_a$  dépend uniquement des points $A$ et $B$ et des tangentes données en ces points et  est la plus grande valeur possible de $R$, pour laquelle 
la courbe de Dubins
$\mathcal{G}(R)$ est dans $\mathcal{E}$.
Ensuite, nous montrons que $\mathcal{J}$ est l'unique courbe de  $\mathcal{E}$, maximisant le minimum du rayon de courbure.

\fi

\section{Introduction}
\label{introduction}

Dans le cadre du brevet  \cite{brevetJB}, il a été nécessaire de construire six courbes 
de classe ${\mathcal{C}}^1$  chacune d'elles passant par un point $A$ et un point $B$ en étant tangente 
respectivement en $A$ et $B$ aux droites $(AO)$ et $(BO)$, où, pour chacune d'elle, $A$, $B$ et $O$ sont trois points donnés du plan.
Plus de détails pourront être trouvés dans 
\cite{enumeration_circuit_JB_2016,%
bastienbrevetforum,%
bastienbrevetmmi,%
bastienbrevetapmep}. Chacune de ces courbes doit relier un des sommets ou un des milieux de côté d'un carré de centre
$O$ et de côté $1$, le point $O$ est fixé, centre du carré et conventionnellement choisi comme repère,  et les points $A$ et $B$ sont l'un des sommets ou un des milieux de côté du carré.
Compte tenu des isométries laissant invariant le carré, seules six courbes ont dû être définies : 
deux segments de droites, deux arcs de cercles et deux portions de paraboles, comme représenté sur 
\cite[Figure 1]{enumeration_circuit_JB_2016}.
On définit 
les deux vecteurs $\vec \alpha$ et $\vec \beta$  et l'angle $\Omega$ de la façon suivante : 
\begin{subequations}
\label{eq10tottot}
\begin{align}
\label{eq10tot}
&\vec \alpha=\frac{1}{OA}\overrightarrow{AO},\quad 
\vec \beta=\frac{1}{OB}\overrightarrow{OB},
\intertext{et}
\label{eq10totnew}
&\Omega=
\left( \widehat{\vec \alpha,\vec \beta} \right).
\end{align}
\end{subequations}
De façon plus générale, on 
se donne trois \ajoutC du plan, $O$,  $A$ et $B$, deux à deux distincts, deux vecteurs unitaires $\vec \alpha$ et $\vec \beta$
\ajoutK
où $\Omega$ n'est pas un multiple de $\pi$. 
Quitte 
à changer de sens de parcours de la courbe, donc à intervertir $A$ et $B$ et à  multiplier $\vec \alpha$ et $\vec \beta$ par $-1$, 
 on peut supposer, sans perte de généralité, que $\Omega$ vérifie 
\begin{equation}
\label{eq20}
\Omega \in ]0,\pi[.
\end{equation}
On s'intéresse à une courbe au moins de classe ${\mathcal{C}}^1$ passant par $A$, tangente à $(OA)$ en $A$, tangente à 
$(OB)$ en $B$. 
On peut choisir une parabole, ce qui a été fait par exemple dans le cas du brevet sur 
les \cite[figures 1e) et 1f)]{enumeration_circuit_JB_2016}.

Ces courbes ont servi à définir des rails \textit{Easyloop} \textregistered\ aptes à faire rouler un train miniature.
Lors de la fabrication des pièces,
la dernière forme, donnée dans 
\cite[La  figure 1f)]{enumeration_circuit_JB_2016}
n'a pas été retenue, puisque trop incurvée, c'est-à-dire, 
que le minimum du rayon de courbure  est trop petit.
Il est en effet nécessaire que le rayon de courbure ne soit pas trop petit pour deux raisons. 
Tout d'abord, les courbes ainsi définies correspondent aux lignes médianes des rails construits : 
les passages des roues et les bords du rails sont définis comme des courbes à distance constante de ces courbes et si le rayon de courbure est trop petit, des points stationnaires avec
des discontinuités de la tangente peuvent apparaître. En outre, les roues des véhicules qui empruntent ces rails doivent pouvoir tourner par rapport au châssis du véhicule et si le rayon de courbure
est trop petit, l'angle de braquage, c'est-à-dire, l'angle entre les essieux qui supportent les paires de roues et l'axe longitudinal du véhicule, est trop important.
Nous proposons donc de déterminer une courbe pour définir la pièce correspondant à 
\cite[figure 1f)]{enumeration_circuit_JB_2016}
qui soit optimale, au sens où le minimum du rayon de courbure est choisi le plus grand possible. 
Nous imposerons aussi à la courbe recherchée d'être à courbure positive. Sans cette hypothèse, le problème est mal posé, puisque
l'on peut construire une courbe formée de trois arcs de cercles, chacun de rayon $R$, avec $R$ arbitrairement grand. 
\ajoutP

Ce problème ressemble fortement à celui des courbes de Dubins \cite{MR0089457,MR0132460}
sans lui être équivalent.
Dubins cherche des courbes  passant aussi par deux points $A$ et $B$ du plan, dont les directions des tangentes en $A$ en $B$ 
sont données. Sous l'hypothèse qu'en tout point, le rayon de courbure est supérieur à $R$ où $R>0$ est donné à l'avance, il cherche la courbe qui minimise la longueur.
Il montre que de telles courbes, entièrement définie par $R$ appelées  géodésiques et notées 
\eqref{rapeldubinseq01}
ne peuvent être que la réunion de trois arcs de cercles de rayon $R$
ou la réunion de deux  arcs de cercle de rayon $R$ et d'un segment, ou une sous partie de ces deux courbes.
\newcommand{\lengthdiftypcodu}{6.5}
\begin{figure}
\psfrag{A}{$A$}
\psfrag{B}{$B$}
\psfrag{O}{$O$}
\psfrag{C}{}
\psfrag{D}{}
\psfrag{E}{}
\psfrag{F}{}
\psfrag{G}{}
\psfrag{a}{$\vec \alpha$}
\psfrag{b}{$\vec \beta$}
\centering
\subfigure[\label{codu0}Le cas $0$ : $0<R<R_a$\ifcase \cras \or /\textit{{The case $0$ : $0<R<R_a$}}\fi]
{\epsfig{file=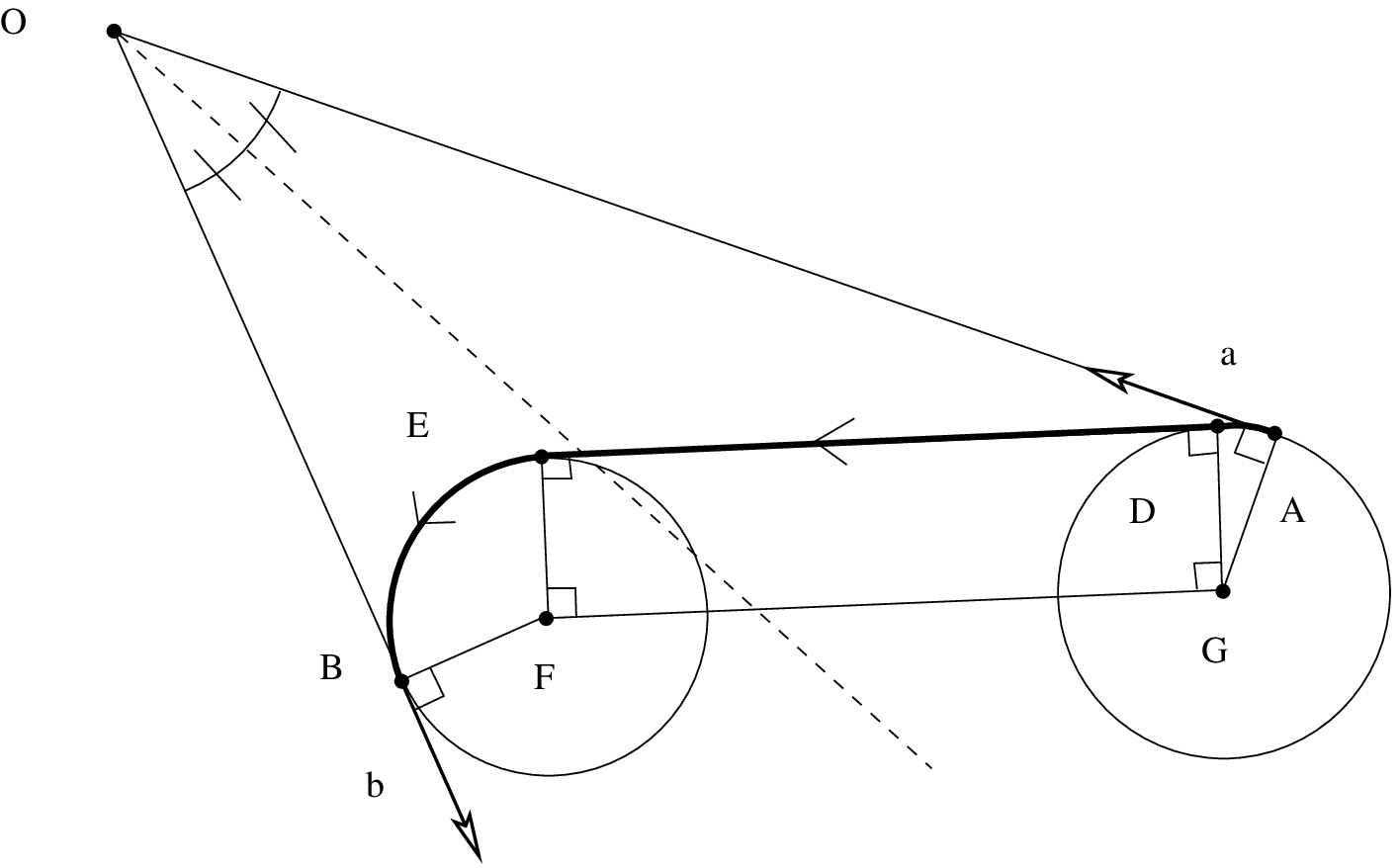, width=\lengthdiftypcodu cm}}
\quad
\subfigure[\label{coducl1}Le cas limite $0-1$ : $R=R_a$\ifcase \cras \or //\textit{{The limit case $0-1$ : $R=R_a$}}\fi]
{\epsfig{file=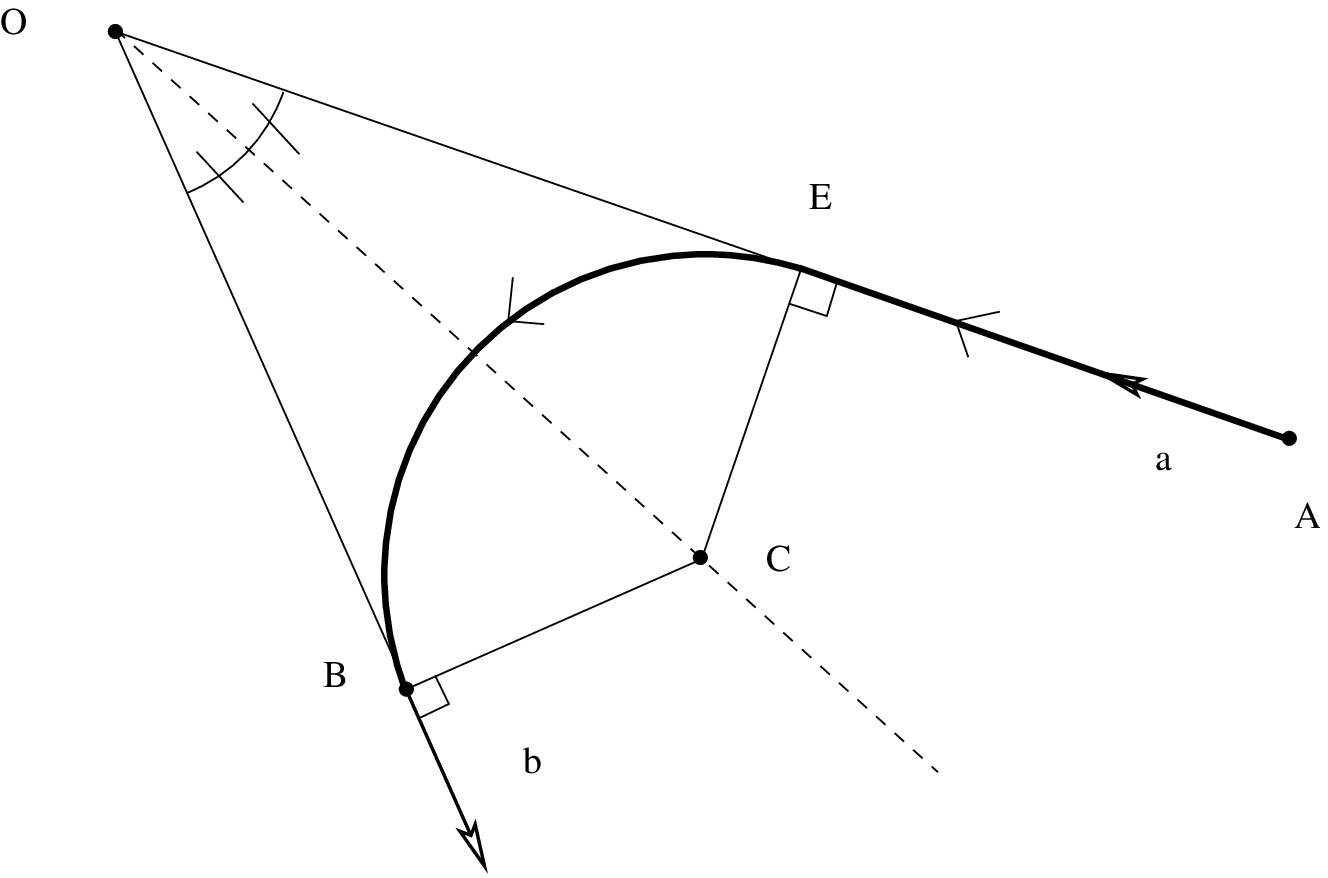, width=\lengthdiftypcodu  cm}}
\caption{\label{diftypcodu}Les différents types de courbes de Dubins définies par $R$ correspondant aux points $A$, $B$ et les vecteurs $\alpha$ et $\beta$ dans les cas $0<R\leq R_a$.\ifcase \cras \or/\textit{The different kinds of Dubins's curves 
defined by $R$ corresponding to points $A$ and $B$ and to vectors  $\alpha$ and $\beta$ in the cases  $0<R\leq R_a$.}\fi}
\end{figure}
Les courbes de Dubins correspondant à $0<R\leq R_a$ (où le rayon $R_a$ ne dépend que des points $O$, $A$ et $B$) sont représentées sur la figure \ref{diftypcodu}.
\ajoutY
considère \definitE. 
Dans notre cas, à la différence des travaux de Dubins, nous ne donnons pas ce rayon de courbure $R$ \textit{a priori} et nous imposons une courbure algébrique positive. 
Nous montrerons qu'il existe une unique courbe de $\mathcal{E}$ formée d'un arc de cercle de rayon $R_a$ et d'un segment de droite.
Ce cas   correspond au cas limite de la figure \ref{coducl1}.
Notre problème, qui ne me semble pas évoqué dans la littérature\footnote{Cette question 
a néanmoins été partiellement soulevée, mais visiblement non résolue sur 
\url{https://math.stackexchange.com/questions/1391778/connect-two-points-given-their-angles-with-a-maximum-radius}},
est donc distinct \textit{a priori} de celui de Dubins.
Ces travaux de Dubins ont été retrouvés plus tard autrement en utilisant le principe de maximum de Pontryagin par exemple dans 
\cite{Boissonnatetal1991,Boissonnatetal1994,sussmann_tang_1991}. 
Très utilisées en robotique et en automatique, ces courbes de Dubins ont fait l'objet de nombreuses recherches.
Voir par exemple les deux thèses suivantes
\cite{lazardtel00442770,Jalel2016}.
De nombreuses variantes existent sur ces courbes de Dubins : 
si des points de rebroussement sont possibles (ce que l'on n'a pas ici, puisque le paramétrage est normal) dans le cas 
où le robot peut inverser sa vitesse \cite{MR1069892} ; 
les recherches prenant en compte les obstacles ont été initiées par Laumond dans 
\cite{Laumond1987};  des problèmes analogues avec des courbes constituées d'arcs de cercle et de segments de longueurs minimales imposées
sont présentés dans 
\cite{MR3158582}.
Notons \sout{enfin} que dans \cite{Moser2009}, un problème proche de  notre problème est évoqué : il s'agit de trouver la
courbe, de longueur donnée (ou inférieure à une longueur donnée) 
qui maximise le minimum du rayon de courbure divisé par le rayon de courbure d'une courbe de référence, donnée à l'avance, 
comme frontière d'un convexe donné.
\ajoutD

Dans la section \ref{enonce}, nous énonçons le problème. 
En section \ref{uniquecercledroite}, nous construirons  l'unique courbe de $\mathcal{E}$ formée d'un segment de droite et d'un arc de cercle 
et nous montrerons ensuite  que cette courbe correspond à une courbe 
 optimale parmi les courbes de Dubins. 
Le \cite[théorème 3.2]{piece6_optimale_JB_2019_X4}  affirme  
qu'il  existe une courbe $X$ de $\mathcal{E}$ qui minimise le maximum de la courbure.
De plus, des simulations numériques présentées dans \cite[section 5]{piece6_optimale_JB_2019_X4} ont corroboré le fait que 
l'unique  courbe de $\mathcal{E}$ formée d'un segment de droite et d'un arc de cercle est bien l'une de celles
qui maximisent le minimum du rayon de courbure.
En section \ref{vraiprobleme}, nous confirmerons cette observation et nous montrerons de façon plus précise  
que l'unique  courbe de $\mathcal{E}$ formée d'un segment de droite et d'un arc de cercle
est l'unique courbe de $\mathcal{E}$,  qui maximise le minimum du rayon de courbure.

\ifcase \cras
On pourra consulter \cite{piece6_optimale_JB_CRAS_2019_ICJ}
où sont présentés de façon condensée les résultats de cet article.
\fi

\section{\'Enoncé du problème}
\label{enonce}

Reprenons le formalisme du papier historique de Dubins  \cite{MR0089457}.
On se donne $O$, $A$ et $B$ trois points deux à deux distincts du plan et $\vec \alpha$ et $\vec \beta$ 
deux vecteurs, vérifiant 
 \eqref{eq10tottot} et  \eqref{eq20}.
Nous cherchons une courbe, paramétrée de façon normale par son abscisse curviligne : 
$s\in [0,L=L(X)]\mapsto X(s)\in \Er^2$, de classe ${\mathcal{C}}^1$.
Pour toute la suite, pour toute telle fonction $X$ de classe  ${\mathcal{C}}^1$, on notera 
\begin{equation*}
L=L(X).
\end{equation*}
On note $\vnorm{.}$ la norme Euclidienne de $\Er^2$. 
On suppose que l'on a 
\begin{subequations}
\label{eq100tot}
\begin{align}
\label{eq100a}
&\forall s\in [0,L],\quad \vnorm{X'(s)}=1,\\
\label{eq100b}
&X(0)=A,\\
\label{eq100c}
&X(L)=B,\\
\label{eq100d}
&X'(0)=\vec \alpha,\\
\label{eq100e}
&X'(L)=\vec \beta.
\end{align}
\end{subequations}
On supposera que 
\begin{equation}
\label{eq110}
X\in W^{2,\infty}(0,L;\Er^2),
\end{equation}
ce qui permet de définir  la valeur absolue de la courbure  $c$  par 
\begin{equation}
\label{eq115old}
|c(s)|= \vnorm{X''(s)}, \quad \aeL.
\end{equation}
La fonction $X$ est dans $ {\mathcal{C}}^1([0,L];\Er^2)$, on a $X'(s)\not =0$ 
et 
on considère  
une détermination continue  de 
l'angle $\phi$ défini par 
\begin{equation}
\label{eq120}
\forall s\in [0,L],\quad 
 \phi(s)=
\left( \widehat{\vec \alpha,X'(s)} \right).
\end{equation}
La fonction $\phi$ est donc dans 
$W^{1,\infty}(0,L)\subset  {\mathcal{C}}^0([0,L])$
et on a 
\begin{equation}
\label{eq130}
\frac{d \phi}{ds}=c,
 \quad \aeL.
\end{equation}
où ici $c$ désigne la courbure algébrique. 
On impose alors 
\begin{equation}
\label{eq140}
\text{$c$ est presque partout positive.}
\end{equation}
Ainsi, \eqref{eq140}
 est équivalent à 
\begin{equation}
\label{eq150}
\text{$\phi$ est croissant.}
\end{equation}
Dans ce cas, on peut réécrire \eqref{eq115old} sous la forme 
\begin{equation}
\label{eq115}
c(s)= \vnorm{X''(s)}, \quad \aeL.
\end{equation}
\ajoutA
Notons que 
\ajoutQ
la condition suivante :
\begin{equation}
\label{eq160}
\forall s\in [0,L],\quad 
 \phi(s)\in [0,\Omega].
\end{equation}
\ajoutG

Dire 
que le minimum du rayon de courbure est le plus grand possible revient à dire
que le maximum de la valeur absolue de la courbure est  le plus petit possible, soit encore, selon 
\eqref{eq140},  que le maximum de la courbure, défini comme $\vnorm[L^\infty(0,L)]{c}$  est  le  plus petit possible.

\begin{definition}
\label{defiE}
\ajoutHb
définit l'ensemble $\mathcal{E}$ des courbes $X$ du plan 
vérifiant 
\eqref{eq100tot},
\eqref{eq110},
\eqref{eq140} (ou \eqref{eq150}).
\end{definition}

D'après \eqref{eq115}, le problème consistera finalement à déterminer une courbe $X$ de $\mathcal{E}$ minimisant 
$\vnorm[L^{\infty}(0,L)]{\vnorm{X''}}$, c'est-à-dire le sup essentiel de la fonction de $[0,L]$ dans $\Er$ : 
$s\mapsto \vnorm{X''}$ : 
\begin{equation}
\label{eq200}
\vnorm[L^{\infty}(0,L)]{\vnorm{X''}}
=\min_{Z\in \mathcal{E}} \vnorm[L^{\infty}(0,L)]{\vnorm{Z''}}.
\end{equation}

Notons enfin que $X$ est de classe ${\mathcal{C}}^1$ et non nécessairement ${\mathcal{C}}^2$.
D'autres travaux  utilisent  
les courbes à courbures continues, utilisant par exemples les clothoïdes comme courbes de raccordement, 
qui permettent une croissance continue de la courbure \cite{MR1717119,scheuertel00001746}. 
En effet, 
dans le cadre de la robotique ou du transport,
il n'est pas possible d'avoir une discontinuité de la courbure, qui implique une discontinuité des accélérations normales, et 
donc des chocs, ce qui use le matériel prématurément ou gêne le voyageur ; un robot ou un véhicule ne peut pas non plus changer instantanément d'angle de braquage.
Au contraire ici,  la discontinuité de la courbure  ne nous gêne pas pour plusieurs raisons. 
Dans le domaine du jouet, 
les masses et les vitesses des véhicules sont très faibles, donc les chocs dus aux discontinuité de l'accélération normale sont négligeables. 
De plus, la notion de confort du voyageur n'a pas de sens.
Les roues des véhicules peuvent subir une discontinuité de l'angle de braquage 
parce qu'elles présentent un léger jeu par rapport au châssis. 
Enfin, la courbe construite dans le cas du brevet 
 \cite{brevetJB,enumeration_circuit_JB_2016} est de classe ${\mathcal{C}}^1$, ${\mathcal{C}}^2$ par morceaux,
mais non ${\mathcal{C}}^2$, puisque formée de portions de segments, de cercles et de paraboles. 
Il n'est donc pas nécessaire de restreindre notre étude aux courbes ${\mathcal{C}}^2$.

\ifcase \cras
\begin{remark}
\or
\begin{remarque}
\fi
Quitte à parcourir, le cas échéant, la courbe dans l'autre sens, on peut 
remplacer respectivement
\eqref{eq20} par "$\Omega \in ]-\pi,\pi[\setminus \{0\}$", 
\eqref{eq140} par "$c$ est presque partout de signe constant", 
\eqref{eq150} par "$\phi$ est monotone" et 
\eqref{eq160} par "$
\forall s\in [0,L],\quad 
 |\phi(s)|\in [0,|\Omega|]$". On cherche toujours la courbe qui minimise le maximum de la courbure géométrique.
\ifcase \cras
\end{remark}
\or
\end{remarque}
\fi

\ifcase \cras
\begin{remark}
\or
\begin{remarque}
\fi
\label{defiErem}
Notons que les résultats de Dubins sont valables pour tout couple de points $(A,B)$ et pour 
tout couple de vecteurs unitaires
$\left(\vec \alpha, \vec \beta\right)$. Ici, on impose les conditions supplémentaires \eqref{eq10tottot} et \eqref{eq20}.
On pourrait croire que le point $O$ peut être construit à partir des points distincts $A$ et $B$ et des vecteurs
$\vec \alpha$ et $\vec \beta$ 
deux vecteurs unitaires donnés vérifiant 
 \eqref{eq20} de la façon suivante :
\begin{equation}
\label{defiEremeq01}
\text{$O$ est l'unique intersection des droites passant respectivement par $A$ et $B$ et dirigées par 
$\vec \alpha$ et $\vec \beta$.}
\end{equation}
Si on le  définit ainsi, le point $O$ n'est pas nécessairement distinct de $A$ et de $B$ 
et \eqref{eq10tot} est alors remplacé \textit{a priori} par 
\begin{equation}
\label{eq10totbis}
\vec \alpha=\pm \frac{1}{OA}\overrightarrow{AO},\quad 
\vec \beta=\pm \frac{1}{OB}\overrightarrow{OB},\quad 
\end{equation}
Cette généralisation est inutile, comme le montre le lemme \ref{lemmeA}, qui servira à plusieurs reprises. 
\ifcase \cras
\end{remark}
\or
\end{remarque}
\fi
\ifcase \cras
\begin{lemma}
\or
\begin{lemme}
\fi
\label{lemmeA}
\ajoutH\
  S'il existe une courbe $X$ du plan 
vérifiant 
\ajoutHc\
\eqref{eq100tot},
\eqref{eq110},
\eqref{eq140} (ou \eqref{eq150}), 
alors $O$ est distinct de $A$ et de $B$ et si on considère les réels $u_0$ et $v_0$ tels que
$\overrightarrow{AO}=u_0 \vec \alpha$ et $\overrightarrow{OB}=v_0 \vec  \beta$, alors $u_0 $ et $v_0$  sont strictement positifs. 
\ifcase \cras
\end{lemma}
\or
\end{lemme}
\fi
\ifcase \cras
Voir la preuve du lemme \ref{lemmeA} en annexe \ref{preuve}, page \pageref{lemmeApreuve}.
\or
Voir la preuve du lemme \ref{lemmeA} en annexe \ref{preuve}.
\fi

\ajoutI

\section{Construction de l'unique courbe de $\mathcal{E}$ formée d'un segment de droite et d'un arc de cercle.}
\label{uniquecercledroite}

\begin{definition}
\label{uniquecercledroitedef01}
On se donne  $O$ et $A$, $B$, deux à deux distincts  et $\vec \alpha$ et $\vec \beta$ vérifiant \eqref{eq10tottot} et \eqref{eq20}.
Nous dirons que nous sommes dans le cas symétrique si $OA =OB$ et dans le cas non symétrique si $OA\not =OB$.
\end{definition}

\ifcase \cras
\begin{lemma}
\or
\begin{lemme}
\fi
\label{uniquecercledroitelem01}
Il existe une unique courbe $X$ de $\mathcal{E}$
formée d'un arc de cercle de rayon $R_a>0$ et de longueur appartenant à  $]0,R_a\pi[$ 
dans le cas symétrique
et 
formée d'un arc de cercle de rayon $R_a>0$ et de longueur appartenant à  $]0,R_a\pi[$ 
et d'un segment de droite
de longueur non nulle dans le cas non symétrique. 
Le rayon $R_a$ du cercle est unique. Il ne dépend que de $O$, $A$ et $B$ et on a 
\begin{equation}
\label{uniquecercledroitelem01eq01}
\vnorm[L^{\infty}(0,L)]{\vnorm{X''}}=\frac{1}{R_a}.
\end{equation}
Pour toute la suite, cette courbe est notée sous la forme $X=\mathcal{J}(O,A,B)$ et le réel $R_a$ sous la forme $R_a(O,A,B)$.
\ifcase \cras
\end{lemma}
\or
\end{lemme}
\fi

La démonstration se fait de façon purement géométrique et est donnée 
\ifcase \cras
en annexe \ref{preuve}, page \pageref{uniquecercledroitelem01preuve}.
\or
en annexe \ref{preuve}.
\fi

\ifcase \cras
\begin{example}
\or
\begin{exemple}
\fi
\label{exemple02}
Traitons le cas particulier donné par 
$
A=(1/2,-1/2),\quad
O=(0,0),\quad
B=(0,-1/2),\quad
\linebreak[1]
\Omega=3\pi/4$
\ajoutOb.
\begin{figure}[h] 
\psfrag{O}{$O$}
\psfrag{A}{$A$}
\psfrag{B}{$B$}
\psfrag{C}{$C$}
\psfrag{D}{$D$}
\psfrag{a}{$\alpha$}
\psfrag{t}{$\theta$}
\begin{center} 
\epsfig{file=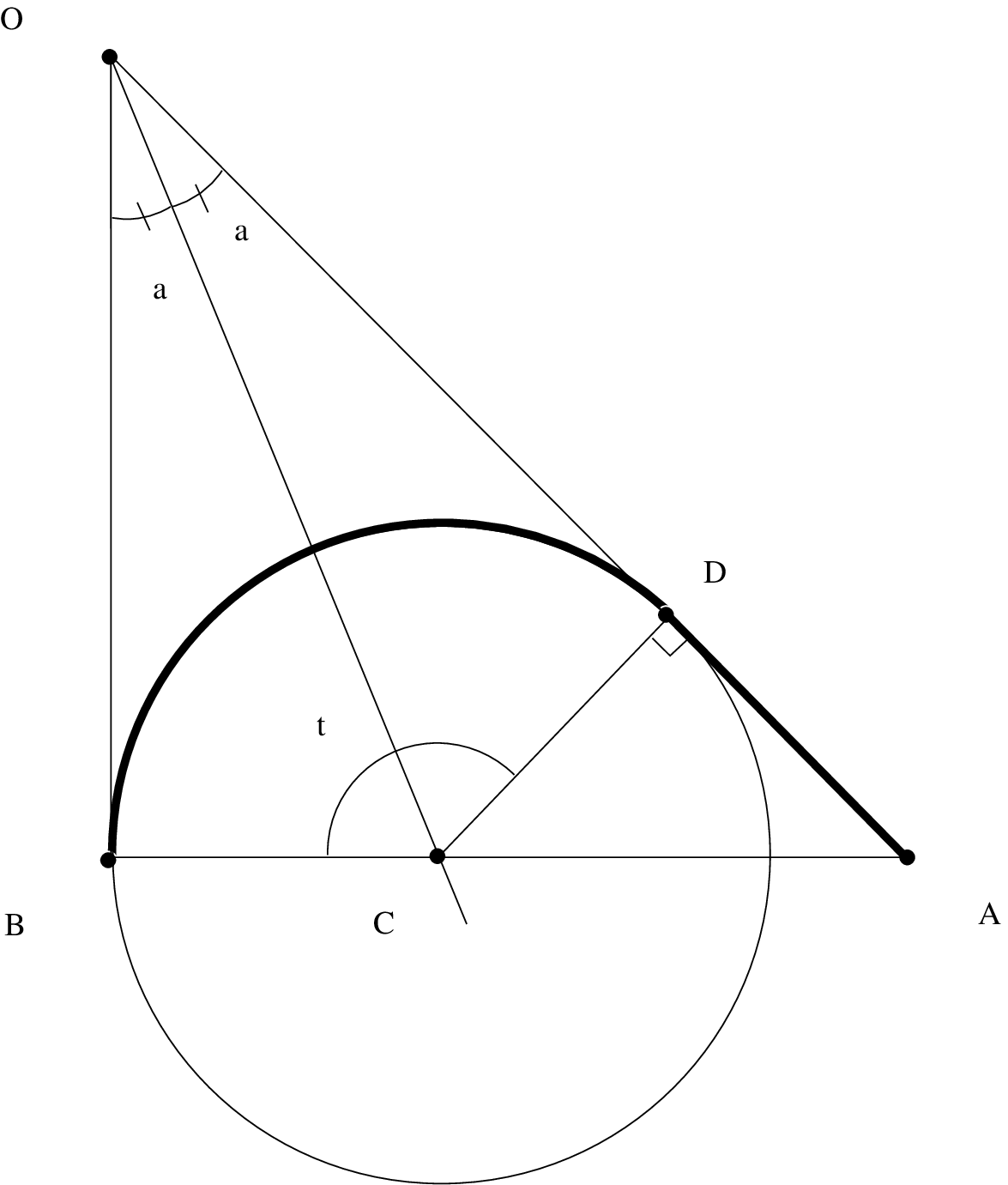, width=4  cm} 
\end{center} 
\caption{\label{forme5poptdet}La description géométrique de la courbe formée d'un arc de cercle et d'un segment de droite.\ifcase \cras \or /\textit{The geometric description of the curve composed of an arc of circle and a line segment.}\fi}
\end{figure}
La construction reprend la méthode donnée dans la preuve du lemme \ref{uniquecercledroitelem01}. 

Comme indiqué sur la figure \ref{forme5poptdet}, la courbe constituée par  la réunion  d'un arc de cercle et d'un segment de droite est définie de la façon suivante (le triangle $OAB$ étant  isocèle rectangle en $B$ avec $OB=BA=1/2$) :  
$(OC)$ est la bissectrice de l'angle $\widehat O$ avec $\alpha=\pi/8$ ;
$\theta =3\pi/4$ ; 
l'arc de cercle a pour centre  $C$ et pour rayon $R$ donné par 
$
R=\frac{\sqrt{2}-1}{2}\approx  0.20710678,
$
 et est limité par les points $B$  et $D$ ;
le segment de droite est le segment $[DA]$ avec $DA=(\sqrt{2}-1)/2$.

La parabole définie dans le cadre du brevet \cite{brevetJB}  \ajoutO\ a été remplacée par cette courbe et cela nous a permis  de faire passer le minimum du rayon de courbure
de $1/25\,\sqrt {5} \approx 0.089$ à $(\sqrt{2}-1)/2\approx 0.207$.
\ajoutL
\ifcase \cras
\end{example}
\or
\end{exemple}
\fi

En reprenant la notation \eqref{rapeldubinseq01}, construisons maintenant autrement la courbe de $\mathcal{E}$, définie dans le lemme \ref{uniquecercledroitelem01}, en utilisant
les courbes de Dubins.

\ifcase \cras
\begin{lemma}
\or
\begin{lemme}
\fi
\label{optdubinslem01}
Considérons 
$
\mathcal{F}=
\left\{R\in \Erps, \quad   \mathcal{G}(R)\in     \mathcal{E}\right\}
$.
Si $R_a(O,A,B)$ est le nombre défini dans le lemme \ref{uniquecercledroitelem01}, alors 
$
\mathcal{F}=]0,R_a(O,A,B)]
$
et la courbe de Dubins   
$\mathcal{G}(R_a(O,A,B))$  
est l'unique courbe définie dans le lemme \ref{uniquecercledroitelem01}.
\ifcase \cras
\end{lemma}
\or
\end{lemme}
\fi

Autrement dit, $\mathcal{G}(R_a(O,A,B))$ est optimale : elle correspond à la plus grande valeur possible de $R$, pour laquelle 
la courbe de Dubins
$\mathcal{G}(R)$ est dans $\mathcal{E}$.
Voir la figure \ref{coducl1} qui correspond au cas optimal.
La démonstration se fait de façon purement géométrique et est donnée en 
\ifcase \cras
annexe \ref{preuve}, page \pageref{optdubinslem01preuve}.
\or
annexe \ref{preuve}.
\fi

\section{Résultat principal : Existence, unicité et caractérisation de la  courbe minimisant le maximum de la courbure.}
\label{vraiprobleme}

On suppose désormais que sont fixés $O$ et $A$, $B$, deux à deux distincts  et $\vec \alpha$ et $\vec \beta$ vérifiant \eqref{eq10tottot} et \eqref{eq20}.

Donnons maintenant le  résultat  essentiel de 
\ifcase \cras
cet article : 
\or
cette Note :  
\fi

\ifcase \cras
\begin{theorem}
\or
\begin{theoreme}
\fi
\label{existenceunicitetheoprop01}
La courbe $\mathcal{J}(O,A,B)$,  donnée dans le lemme \ref{uniquecercledroitelem01}, est l'unique  courbe $X$ de $\mathcal{E}$ qui minimise le maximum de la courbure, c'est-à-dire vérifiant 
\eqref{eq200}. 
\ifcase \cras
\end{theorem}
\or
\end{theoreme}
\fi

Donnons tout d'abord la proposition suivante : 

\begin{proposition}
\label{lemme100}
Avec les notations habituelles, nous supposerons, sans perte de généralité que $OA\leq OB$, de sorte 
que la courbe $X=\mathcal{J}(O,A,B)\in \mathcal{E}$, de longueur $L$ décrite dans le lemme \ref{uniquecercledroitelem01}, commence d'abord par un arc de cercle de rayon $R_a=R_a(O,A,B)$. Notons 
\ajoutZ\
l'angle associé à cette courbe.
Soit maintenant une autre courbe $Z\in \mathcal{E}$, associé à l'angle 
\ajoutZb\
et de longueur $M$. On note
$e=\vnorm[L^{\infty}(0,L)]{\vnorm{Z''}}$ et on suppose que 
\begin{equation}
\label{lemme100eq10}
e\leq \frac{1}{R_a}.
\end{equation}
On se place dans le 
\ajoutT\
 et  on note respectivement $(x(s),  y(s))$ et $\left(\widehat x(s),  \widehat y(s)\right)$  les coordonnées de $X(s)$ et de $Z(s)$  dans ce repère.
On considère la longueur $l$ de la partie circulaire de la courbe $X$ définie par 
\begin{equation}
\label{lemme100eq20}
l=R_a\Omega.
\end{equation}
Sous l'hypothèse
\begin{equation}
\label{eq20bis}
\Omega \in ]0,\pi/2[,
\end{equation}
alors on a 
\begin{subequations}
\label{lemme100eq30tot}
\begin{align}
\label{lemme100eq30a}
&M\geq l,
\intertext{et on peut donc poser :}
\label{lemme100eq30b}
&\zeta_0=-\left(\widehat x\left(l\right) -x\left(l\right)\right)\sin\left(\phi\left(l\right)\right)
+
\left(\widehat y\left(l\right)-  y\left(l\right) \right)\cos\left(\phi\left(l\right)\right)\, ;
\intertext{on a:}
\label{lemme100eq30c}
& \zeta_0 \leq 0,\\
\label{lemme100eq30d}
& e<\frac{1}{R_a}\Longrightarrow \zeta_0 < 0,\\
\label{lemme100eq30e}
&  \zeta_0=0 \Longrightarrow  \left(e=\frac{1}{R_a},\quad M=L,\quad Z=X\right).
\end{align}
\end{subequations}
\end{proposition}

Seules, les idées principales sont données. La 
preuve complète est fournie en annexe \ref{preuvecomplete}, page \pageref{tugudu00}.

\ifcase \cras
\begin{proof}[Idées de la démonstration]
\or
\begin{preuveidee}
\fi

On a les relations habituelles  
\begin{subequations}
\label{eqderphinbbbb}
\begin{align}
\label{eqderphinab}
&\frac{d  x}{ds}=\cos   \phi,\\
\label{eqderphinbb}
&\frac{d y}{ds}=\sin   \phi,\\
\label{eqderphinabnew}
&\frac{d  \widehat x}{ds}=\cos   \theta,\\
\label{eqderphinbbnew}
&\frac{d \widehat y}{ds}=\sin   \theta.
\end{align}
\end{subequations}
On a  $M\geq l$ et il est légitime de poser :
\begin{equation*}
\forall s\in [0,l],\quad
\zeta(s)=
-\left(\widehat x(s) -x(s)\right)\sin\left(\phi(s)\right)
+
\left(\widehat y(s)-  y(s) \right)\cos\left(\phi(s)\right).
\end{equation*}
On note $\prodsca{.}{.}$, le produit scalaire Euclidien de $\Er^2$ (qui induit la norme Euclidienne  $\vnorm{.}$  de $\Er^2$).  Ainsi, 
géométriquement, $\zeta$ correspond à 
 composante du vecteur $\overrightarrow{X(s)Z(s)}$ sur $\vec N(s)$, qui désigne la normale extérieure à la courbe $X$ au point d'abscisse curviligne $S$.
On montre tout d'abord que 
\begin{equation*}
\forall s\in[0,l],\quad
\theta(s)-\phi(s)
\leq \left(e-\frac{1}{R_a}\right)s.
\end{equation*}
Ensuite, grâce à \eqref{eqderphinbbbb}, on montre que 
\begin{align*}
&\forall s\in[0,l],\quad
\frac{d}{ds}\left(\widehat x(s) -x(s)\right)
=
-2\sin\left(\frac{\theta(s)+\phi(s)}{2}\right)\sin\left(\frac{\theta(s)-\phi(s)}{2}\right).\\
&\forall s\in[0,l],\quad
\frac{d}{ds}\left(\widehat y(s) -y(s)\right)
=
2\cos\left(\frac{\theta(s)+\phi(s)}{2}\right)\sin\left(\frac{\theta(s)-\phi(s)}{2}\right),
\end{align*}
puis que
\begin{align*}
 \forall s\in]0,l],\quad
&\zeta(s)\leq 0,
\intertext{et}
e<\frac{1}{R_a}
\Longrightarrow 
 \forall s\in]0,l],\quad
&\zeta(s)<0.
\end{align*}
En particulier, en $s=l$, on obtient 
\eqref{lemme100eq30c}
et 
\eqref{lemme100eq30d}.
Enfin, 
grâce aux hypothèses  \eqref{eq150} et \eqref{eq20bis}, on obtient la monotonie des fonctions
$\widehat x(s) -x(s)$ et $\widehat y(s) -y(s)$, ainsi 
que leur signe constant.

La nullité de $\zeta_0$ impose 
$
\widehat x\left(l\right) -x\left(l\right)=
\widehat y\left(l\right) -y\left(l\right)=0$ 
et la monotonie de $\widehat x-x$ et $\widehat y-y$ impose donc leur nullité sur tout l'intervalle $[0,l]$. 
On en déduit 
 que $\phi$ et $\theta$ coïncident sur$[0,l]$ puis  sur  $[0,M]$ et il en est de même pour $x$ et $\widehat  x$ et $y$ et $\widehat  y$, 
en particulier en $s=M$ où $x$ et $y$ valent $x_B$ et $y_B$. Les deux courbes $X $ et $Z$ finissent donc au même point $B$ et donc $L=M$
et $X=Z$. On en déduit alors $e=1/R_a$. On a donc montré \eqref{lemme100eq30e}.
\ifcase \cras
\end{proof}
\or
\end{preuveidee}
\fi

On est en mesure maintenant d'utiliser le lemme \ref{lemmeconvexite}, conséquence 
de la croissance des angles $\theta$ et $\phi$ pour donner la preuve du théorème  \ref{existenceunicitetheoprop01}.

Seules, les idées principales sont données. La 
preuve complète est fournie en annexe \ref{preuvecomplete}, page \pageref{tugudu10}.

\ifcase \cras
\begin{proof}[Idées de la démonstration du théorème \ref{existenceunicitetheoprop01}]
\or
\begin{preuveidee}[ du théorème \ref{existenceunicitetheoprop01}]
\fi

Notons tout d'abord que 
\ajoutAA\
existe, puisque $\mathcal{E}$ est non vide et que 
\ajoutAAb\
est toujours positif.

\begin{enumerate}

\item
Premier cas :  
\ajoutRb\

On peut alors utiliser directement la proposition 
\ref{lemme100}.
Pour toute la suite, on considère la fonction $g$ de $\mathcal{E}$ dans $\Erp$ définie par 
\begin{equation}
\label{pdueq01}
\forall X\in \mathcal{E},\quad
g(X)=\vnorm[L^{\infty}(0,L(X))]{\vnorm{X''}}.
\end{equation}
On a exhibé dans le lemme \ref{uniquecercledroitelem01} une fonction $X=\mathcal{J}(O,A,B)$ de $\mathcal{E}$ vérifiant $g(X)=1/R_a$ où $R_a=R_a(O,A,X)$. On a donc
$\inf_{Z\in \mathcal{E}} g(Z)\leq {1/R_a}$.

Pour montrer que le minimum de $g$ est atteint est vaut $1/R_a$, l'idée simple est de démontrer que  si la courbure maximale d'une courbe
est trop faible alors $\zeta_0$ est négatif, ce qui contredit son aspect convexe (lemme \ref{lemmeconvexite}). 
S'il  existe  une courbe de $\mathcal{E}$ notée $Z$ telle que 
$ g(Z)< 1/{R_a}$, d'après la proposition  
\ref{lemme100}, on a donc $\zeta_0 < 0$.
Il existe donc un point de la courbe $Z\left(l\right)$ de $Z$ dans le demi-plan ouvert $\Pi$ défini par la tangente à la courbe $X$, notée $\mathcal{D}$
au point d'abscisse curviligne $l$, du côté opposé à la normale. 
Cela contredit l'appartenance de  la courbe $Z$  à  $\mathcal{E}$.
En effet, d'après le lemme \ref{lemmeconvexite},
la courbe $Z$ est incluse dans le demi-plan délimitée par la droite tangente à la courbe au point $B$, 
du côté la normale extérieure à la courbe en $B$, donc de l'autre côté du demi-plan $\Pi$  délimité par $\mathcal{D}$. 
On a donc 
\begin{equation*}
\inf_{Z\in \mathcal{E}}g(Z)=g(X)=\min_{Z\in \mathcal{E}}g(Z)=\frac{1}{R_a}.
\end{equation*}

Montrons enfin l'unicité de la courbe $X$ vérifiant 
\eqref{eq200} en montrant que cette courbe est la courbe $\mathcal{J}(O,A,B)$  donnée dans le lemme \ref{uniquecercledroitedef01}.
On utilise là encore la proposition  \ref{lemme100}.
Supposons qu'il existe une autre courbe $Z$ de $\mathcal{E}$ telle que $
g(Z)=1/R_a$,
où $R_a=R_a(O,A,B)$. 
D'après le lemme \ref{lemmeconvexite} appliqué à cette courbe $Z$ et à la droite $\mathcal{D}$, la tangente à la courbe $Z$ au point $B$, on a donc
 avec les notations de la proposition 
\ref{lemme100}, $\zeta_0\geq 0$.
D'après cette  même proposition, on a $\zeta_0\leq 0$ et donc $\zeta_0=0$ et de nouveau d'après  la proposition
\ref{lemme100}, on a $Z=X$.

Le théorème \ref{existenceunicitetheoprop01} est donc vrai sous l'hypothèse \eqref{eq20bis}.

\item
Second cas : 
\ajoutR\

Nous allons décomposer le problème en deux sous-problèmes et à chacun d'eux, nous appliquerons 
le théorème \ref{existenceunicitetheoprop01}, sous l'hypothèse \eqref{eq20bis}.

Commençons par montrer un résultat similaire à la proposition 
\ref{lemme100}. 
Supposons qu'il existe une courbe  $\mathcal{E}$ notée $Z$, de longueur $M$,  telle que, en notant $e=g(Z)$, $
e \leq 1/R_a.$
Comme dans la proposition \ref{lemme100},
notons $\theta$ l'angle associé à cette courbe.
Considérons $s_1\in ]0,M[$  tel que  $
\theta(s_1)={\Omega}/{2}$.
Considérons $Z_1$, la restriction de $Z$ à $[0,s_1]$ et $Z_2$ la restriction de $Z$ à $[s_1,M]$ et posons 
$
C=Z(s_1)$.
La courbe $Z_1$ relie donc les deux points $A$ et $C$ tandis que la courbe
$Z_2$ relie $C$ et $B$. D'après le lemme \ref{lemmeA}, les deux droites respectives passant par $A$, portée par $\vec \alpha$
et passant par $C$ portée par $Z'(s_1)$ se coupent en point $E$, distinct de $C$ et de $A$.
De même, les deux droites respectives passant par $B$, portée par $\vec \beta$
et passant par $C$ portée par $Z'(s_1)$ se coupent en point $F$, distinct de $C$ et de $B$.
On peut donc appliquer le théorème \ref{existenceunicitetheoprop01}
 sous l'hypothèse \eqref{eq20bis}
aux trois points $A$, $E$ et $C$ puisque,
$
\left(\widehat{\vec \alpha,{\overrightarrow{EC}}/{{EC}}}\right)=
{\Omega}/{2}\in]0,\pi/2[.$
Il existe donc une unique courbe à courbure positive
$\widetilde Z_1=\mathcal{J}(E,A,C)$ passant par $A$ et $C$, dont les tangentes sont portées en ces points par $\vec\alpha$ et $Z'(s_1)$,
qui minimise 
\ajoutAB\
D'après le théorème \ref{existenceunicitetheoprop01}
\ajoutAC , 
cette courbe est la réunion d'un arc de cercle de rayon $R_1$, 
d'angle $\Omega/2$ et d'un segment de droite de longueur $l_1\geq 0$. 
Puisque $Z_1$ est dans ${\mathcal{E}}_1$, 
on a donc 
$
 \vnorm[L^{\infty}(0,s_1)]{\vnorm{Z_1''}}\geq   \vnorm[L^{\infty}(0,L(\widetilde Z_1))]{\vnorm{\widetilde Z_1}}  =1/R_1,
$
ce qui implique
$$
R_a\leq R_1\text{ et } \left(e<{1}/{R_a} \Longrightarrow R_a<R_1\right).$$
De même, 
il existe une unique courbe à courbure positive
$\widetilde Z_2=\mathcal{J}(F,C,B)$ passant par $C$ et $B$, dont les tangentes sont portées en ces points par  et $Z'(s_1)$ et $\vec\beta$. On a aussi 
$$
R_a\leq R_2 \text{ et } \left(e<{1}/{R_a} \Longrightarrow R_a<R_2\right).$$
$\widetilde Z$ est la réunion de la courbe $\widetilde Z_1$ et de $\widetilde Z_2$ et est donc la réunion, au plus, de deux arcs de cercles et de deux segments de droites. 
On pose  de façon analogue à 
\eqref{lemme100eq50}
\begin{equation}
\forall s\in [0,L(\widetilde Z)],\quad
\zeta(s)=
-\left(\widetilde  x(s) -x_B\right)\sin \Omega
+
\left(\widetilde   y(s)-  y_B \right)\cos\Omega,
\end{equation}
où $\widetilde Z(s)=(\widetilde  x(s),\widetilde   y(s))$.
On pose enfin
 $\zeta_0= \zeta(\widehat L)$,
où $\widehat L$ est la somme des longueurs des deux premiers segments de droite et des deux arcs de cercle.
Comme dans la preuve de la proposition 
\ref{lemme100}, on obtient 
\begin{align*}
&\zeta_0 \leq 0,\\
& \left(R_1 >R_a \text{ et } R_2>R_a\right)\Longrightarrow \zeta_0 < 0,\\
&  \zeta_0=0 \Longrightarrow  \left(R_1=R_a \text{ et } R_2=R_a,\quad Z=\mathcal{J}(O,A,B)\right).
\end{align*}
En effet, il existe $a$, $b$, $c$ et $f$, des réels strictement négatifs tels que $
\zeta_0=a(R_1-R_a)+b(R_2-R_a)+cd_1+fd_2$.
Le fait que $\zeta_0\leq 0$ signifie que la courbe $\widetilde Z$ a des points de l'autre côté (par rapport à la normale) de la tangente en $B$  la courbe.

On conclut comme dans le premier cas.
\end{enumerate}
\ifcase \cras
\end{proof} 
\or
\end{preuveidee}
\fi

\ifcase \cras

\section{Construction effective de la pièce du circuit et exemple d'un circuit}
\label{construceffct}

Si on choisit les dimensions de la section standard Brio, choisis pour les rails Easyloop, on obtient donc finalement, 
en utilisant la courbe définie dans l'exemple \ref{exemple02},  la 
pièce 6 représentée sur la figure \ref{fin}.

\begin{figure}
\centering
\subfigure[\label{fin} La forme optimale.\ifcase \cras \or /\textit{The optimal shape.}\fi] 
{\epsfig{file=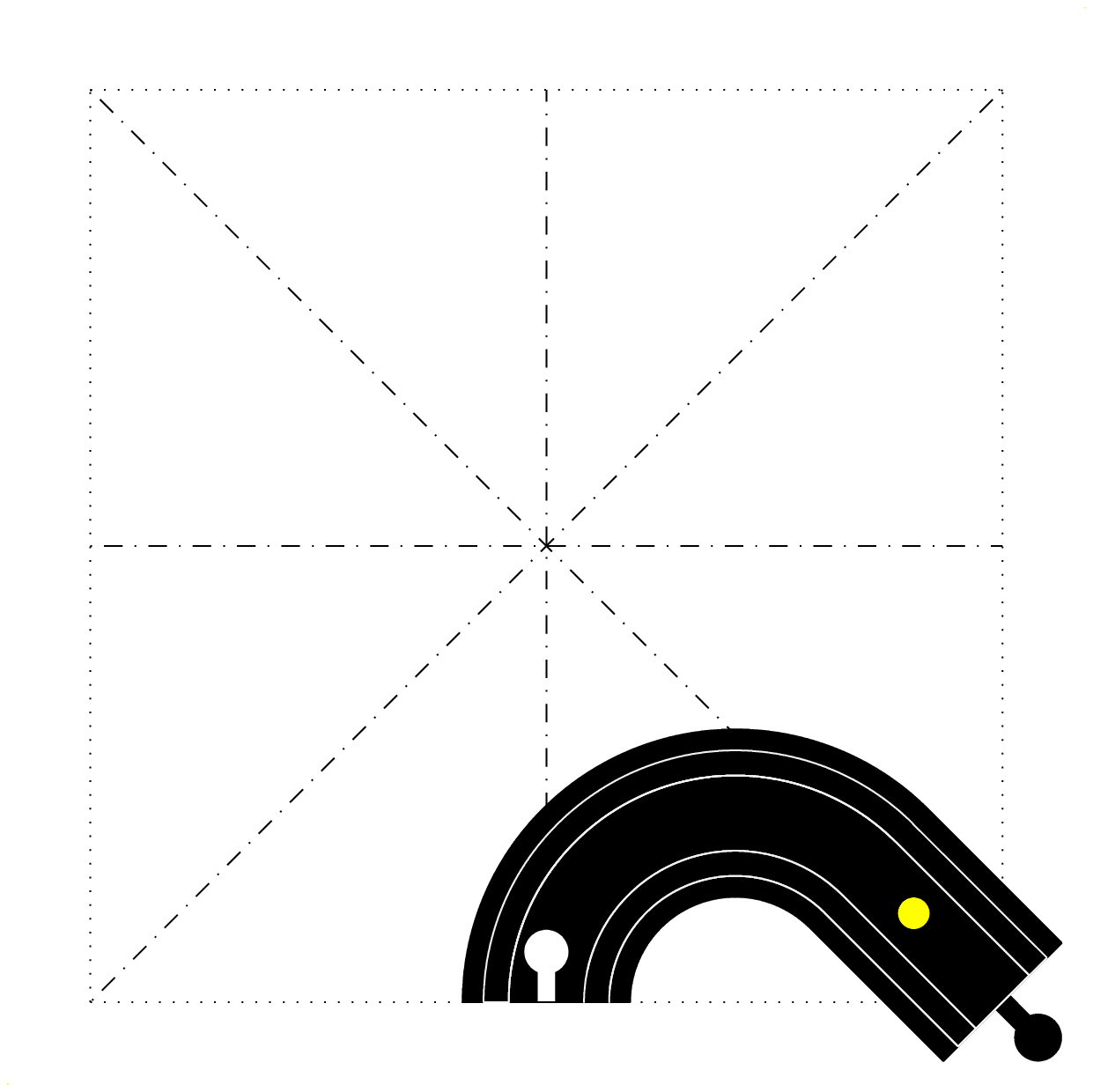, width=6 cm}}
\qquad
\subfigure[\label{exemplemcarder01b}Un exemple de circuit avec la pièce 6 optimale.\ifcase \cras \or /\textit{An exemple of track with optimal shape 6.}\fi]
{\epsfig{file=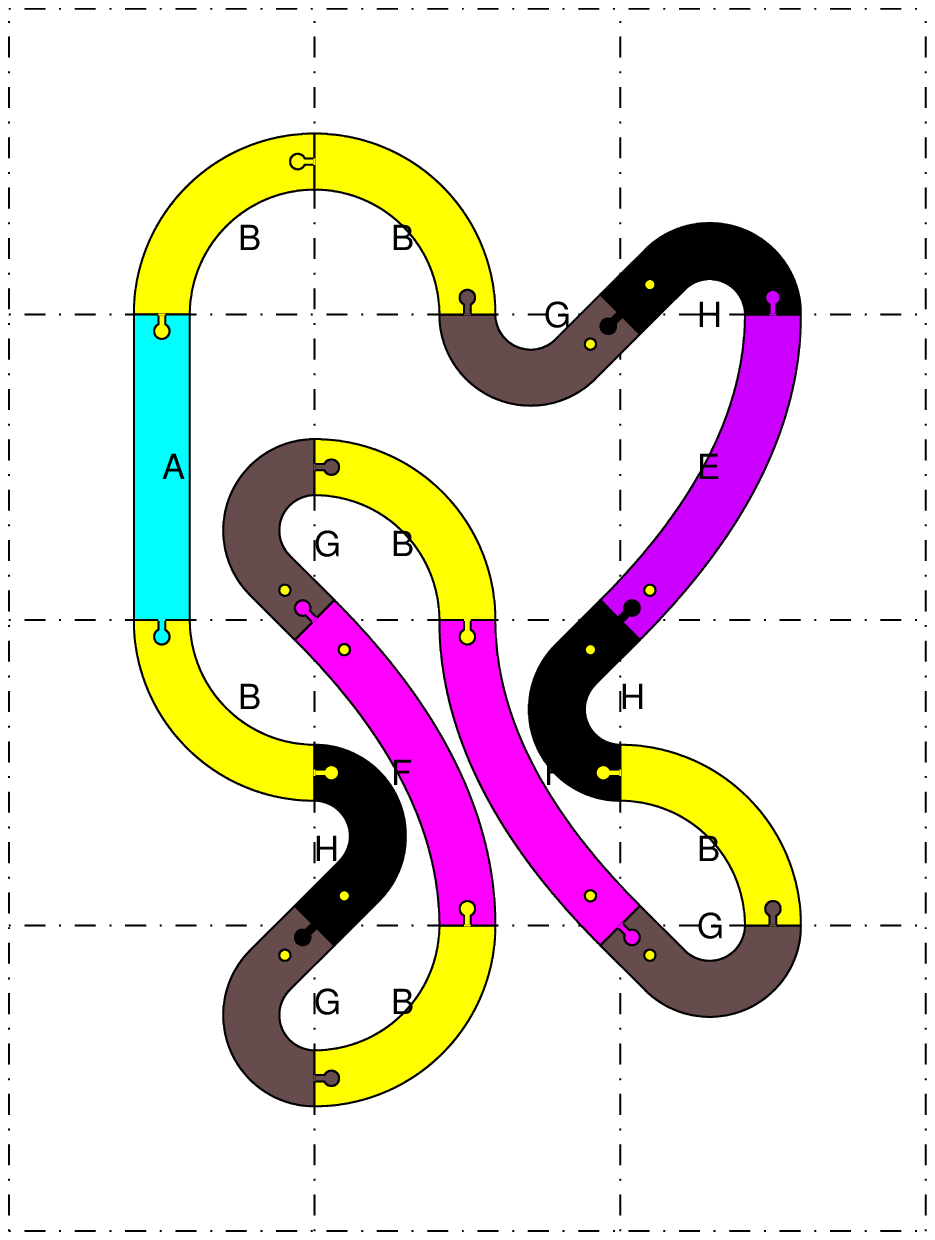, width=5 cm} }
\caption{\label{finexemplemcarder01b} Utilisations de la courbe optimale. \ifcase \cras \or /\textit{Usings the optimal curve}\fi}
\end{figure}

On pourra aussi consulter la figure \ref{exemplemcarder01b} qui montre un exemple d'un ciruit contenant cette  pièce optimale.

\fi

\section{Conclusion}
\label{conclusion}

On a montré qu'il existe une unique courbe de  \definitE, 
minimisant le maximum du rayon de courbure qui est l'unique courbe de $\mathcal{E}$  formée d'un arc de cercle et d'un segment de droite.

\ajoutM questions restent en suspens : 
\begin{itemize}
\item
L'unique  courbe minimisante 
$\mathcal{J}$ trouvée constitue un cas particulier des courbes de Dubins.
Cependant, autant dans la formulation du problème que sa résolution, ce résultat semble différent des célèbres travaux de Dubins. Fondamentalement, pourquoi, pour $R=R_a(O, A, B)$ donné, chercher une courbe,
de rayon de courbure   supérieur à $R_a$
\ajoutS\
 donne le même résultat que chercher une courbe, à courbure positive, à maximum du rayon de courbure minimal ?
\item
Comme le précise la 
\cite[remarque A.2]{piece6_optimale_JB_2019_X4}, le
\cite[théorème 3.2]{piece6_optimale_JB_2019_X4}
est encore valable si l'on remplace l'hypothèse \eqref{eq150}
par l'hypothèse plus générale \eqref{eq160}. En revanche,
\ajoutSb\
si on ne fait plus l'hypothèse \eqref{eq150} qui assure la convexité des éléments de $\mathcal{E}$.
Une question pour l'instant ouverte 
serait de reprendre le résultat 
essentiel de 
\ifcase \cras
cet article 
\or
cette Note, 
\fi
le théorème \ref{existenceunicitetheoprop01},  en remplaçant l'hypothèse 
\eqref{eq150}
par l'hypothèse plus générale \eqref{eq160}, en cherchant toujours à maximiser le minimum  du rayon de courbure 
 de la courbe. 
\end{itemize}

\ifcase \cras
\or

\section*{Remerciements}

Je remercie vivement les rapporteurs anonymes pour la qualité et la précision de leurs remarques qui m'ont aidé à améliorer cette Note et les travaux préliminaires \cite{piece6_optimale_JB_2019_X4}.

\fi

\ifcase \cras

\or

\bibliographystyle{unsrt}
\bibliography{piece6_optimale_JB_arXiv[X6]}

\fi

\appendix

\newcommand{\refpreuve}{A}
\newcommand{\refpreuvecomplete}{B}

\numberwithin{figure}{section}
\numberwithin{equation}{section}

\clearpage

\ifcase \cras
\markboth{J\'ER\^OME BASTIEN}%
{\ftitrefacul}
\fi

\ifcase \cras
\or
\textit{Les annexes \ref{preuve} 
 et \ref{preuvecomplete} sont destinées aux rapporteurs et ne font pas partie du projet de Notes.
Il est en effet  écrit sur le site des CRAS : 
"Il n'est parfois pas possible, en raison de la concision exigée, de donner les démonstrations ou les preuves complètes du résultat énoncé, spécialement dans la série I ; il est alors recommandé de joindre à l'appui de la Note un texte, si possible dactylographié, explicitant les compléments nécessaires à une bonne compréhension qui facilite l'examen de la Note par le présentateur et qui sera conservé cinq ans dans les Archives de l'Académie pour pouvoir être communiqué à tout lecteur des Comptes rendus qui en fait la demande. Il est alors fait mention au bas de la Note de l'existence de ce document par la formule suivante ; "Résumé d'un texte qui sera conservé cinq ans dans les Archives de l'Académie et dont copie peut être obtenue."  "
Les preuves de l'annexe \ref{preuve} viennent d'être déposées\ dans \url{https://arxiv.org}. Voir \cite{piece6_optimale_JB_2019_X4}.
L'annexe \ref{preuvecomplete} contient toutes les preuves des résultats essentiels.
Le cas échéant, ces preuves feront l'objet d'une publication ultérieure.}
\fi

\section{Preuves annexes}
\label{preuve}


\ifcase \cras
\or
\fi

Rappelons tout d'abord que, comme les courbes représentatives  d'applications convexes, 
on a le lemme suivant :


\input{textelemmeconvexite}

\ifcase \cras

\ifcase \cras
\begin{proof}\
\or
\begin{preuve}\
\fi

\begin{figure}[h]    
\psfrag{A}{$A$}
\psfrag{B}{$B$}
\psfrag{a}{$\vec \alpha$}
\psfrag{b}{$\vec \beta$}
\begin{center} 
\epsfig{file=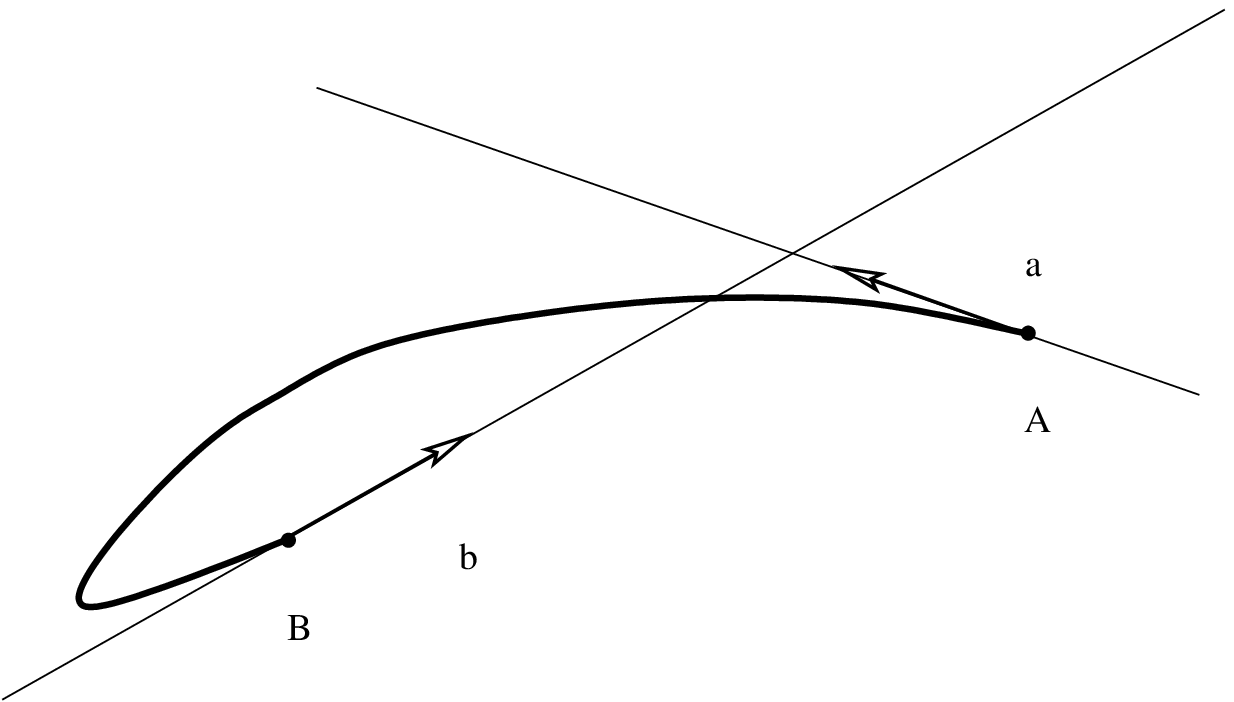, width=6 cm}  
\end{center} 
\caption{\label{contre_exemple_convexite}Contre-exemple au lemme  \ref{lemmeconvexite} sans l'hypothèse \eqref{eq20}.\ifcase \cras \or /\textit{Counter-example to Lemma  \ref{lemmeconvexite} without assumption 
\eqref{eq20}.}\fi}
\end{figure}

Notons que, sans l'hypothèse \eqref{eq20}, cette propriété devient fausse comme le montre le contre-exemple de la figure 
\ref{contre_exemple_convexite}.

\begin{figure}[h]    
\psfrag{a}{$X'(s)$}
\psfrag{b}[][l]{$\vec N(s)$} 
\psfrag{P}{$X(s)$}
\psfrag{Q}{$X(t)$}
\begin{center} 
\epsfig{file=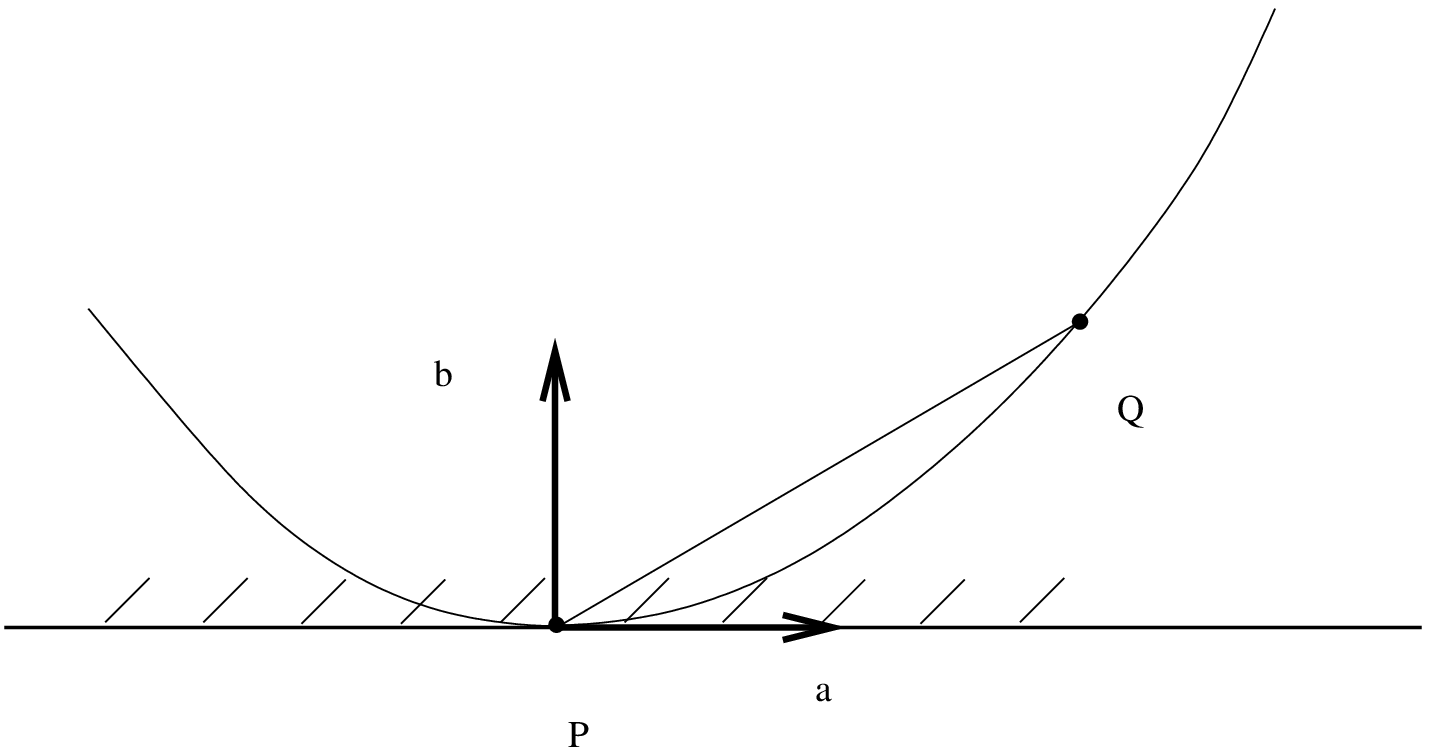, width=9 cm}  
\end{center} 
\caption{\label{courbe_meme_cote}La courbe est toujours du côté de la normale extérieure.\ifcase \cras \or /\textit{The curve is ever in the side of the outer-pointing normal.}\fi}
\end{figure}

Notons que 
\begin{equation*}
N(s)=\sigma(X'(s)),
\end{equation*}
où 
\ifcase \cras
$\sigma$    est la rotation vectorielle d'angle $\pi/2$.
\or
$\sigma$    est la rotation vectorielle d'angle $\pi/2$, déjà utilisée dans la preuve de la proposition \ref{lemme100}.
\fi
Voir sur la figure \ref{courbe_meme_cote}, la situation représentée. 
\ifcase \cras
On note $\prodsca{.}{.}$, le produit scalaire Euclidien de $\Er^2$ (qui induit la norme Euclidienne  $\vnorm{.}$  de $\Er^2$). 
\or
$\prodsca{.}{.}$ est le produit scalaire Euclidien de $\Er^2$.
\fi
Pour $s\in [0,L]$ fixé et pour tout $t\in [0,L]$, on pose 
\begin{equation*}
\gamma(t)=\prodsca{X(t)-X(s)}{\sigma(X'(s))},
\end{equation*}
dont la dérivée vaut d'après \eqref{eqderphinbbbb}
\begin{align*}
\uwave{\gamma'(t)}
&=\prodsca{X'(t)}{\sigma(X'(s))},\\
&=
\begin{pmatrix}
\cos   \phi(t) \\ \sin   \phi(t)
\end{pmatrix}
.
\begin{pmatrix}
-\sin\phi(s)\\ \cos\phi(s)
\end{pmatrix},\\
&=
\sin   \phi(t) \cos\phi(s)-\cos   \phi(t)\sin\phi(s),
\end{align*}
et donc 
\begin{equation}
\label{ffvsjhgfss}
\forall t\in [0,L],\quad
\gamma'(t)=
\sin\left(\phi(t)-\phi(s)\right).
\end{equation}
Or,  d'après 
\eqref{eq20},
\eqref{eq150} et 
\ajoutJJ
\label{ajoutJJ20}
on a, pour tout $t$, 
\begin{equation*}
-\pi<-\Omega\leq 
\phi(t)-\phi(s) \leq \Omega<\pi,
\end{equation*}
et, si $t\geq s$, 
\begin{equation*}
\phi(t)-\phi(s)\geq 0,
\end{equation*}
et donc, d'après \eqref{ffvsjhgfss}, $\gamma'(t)\geq 0$.
De même, si $t\leq s$, $\gamma'(t)\leq 0$. Or, on a $\gamma(s)=0$ et donc pour tout $t$, $\gamma(t)\geq 0$. 
$\gamma(t)$ représente la composante du vecteur $X(t)-X(s)$ sur $\sigma(X'(s))$, ce qui nous permet de conclure. 
\ifcase \cras
\end{proof}
\or
\end{preuve}
\fi

\or
\sout{La preuve, fondée sur  \eqref{eq20} et \eqref{eq150}, est laissée au lecteur.}
\fi


\ifcase \cras
\begin{proof}[Démonstration du lemme \ref{lemmeA}]\
\or
\begin{preuve}[ du lemme   \ref{lemmeA}]\
\fi

\label{tugudu80}

\ifcase \cras
\label{lemmeApreuve}
\fi

Soit $X$ 
vérifiant 
\ajoutV\
\eqref{eq20},
\eqref{eq100tot},
\eqref{eq110},
\eqref{eq140} (ou \eqref{eq150}).

\begin{figure}[h]    
\psfrag{jt}{$\vec {k}$}
\psfrag{it}{$\vec {\alpha}$}
\psfrag{O}{$O$}
\psfrag{A}{$A$}
\psfrag{B}{$B$}
\psfrag{X}{$X$}
\psfrag{Q}{$Q(s)$}
\psfrag{Xp}{$X'(s)$}
\psfrag{u}{$u$}
\psfrag{v}{$v$}
\begin{center} 
\epsfig{file=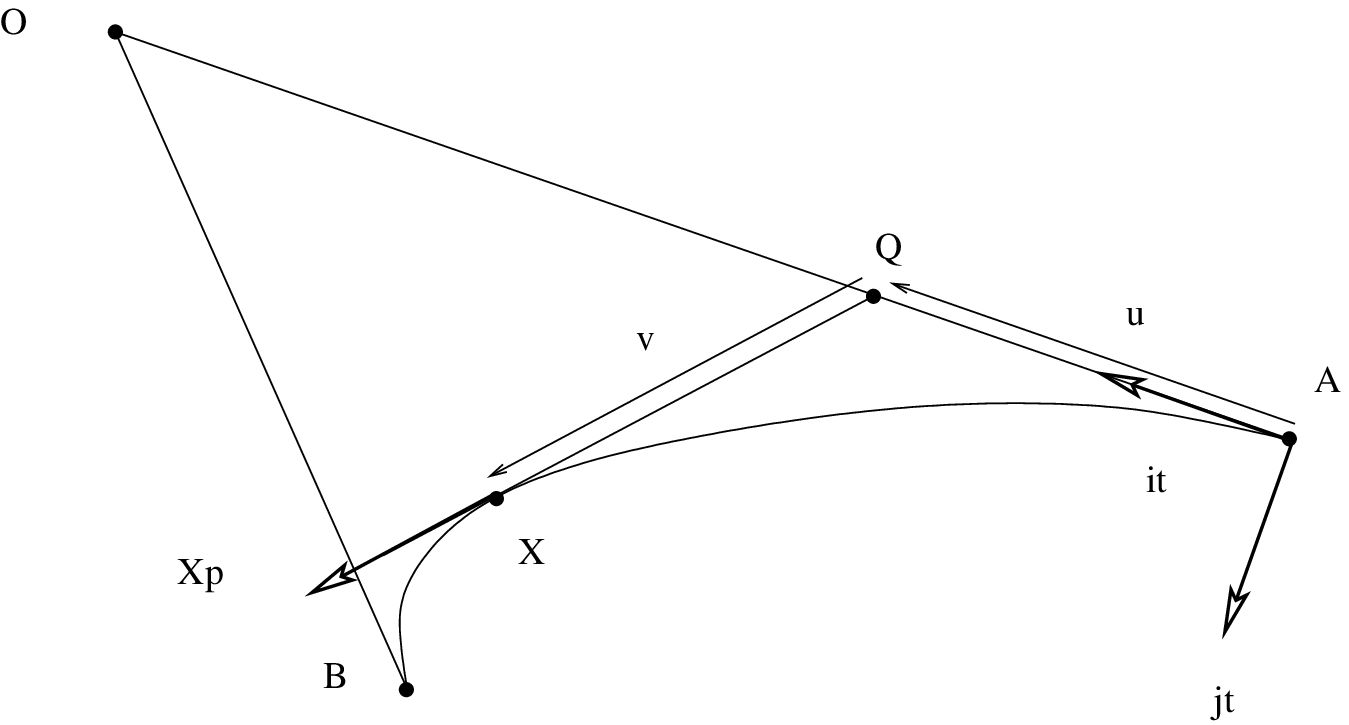, width=9 cm} 
\end{center} 
\caption{\label{longueur_majoreeb}Le repère $\left(A,\vec {\alpha}, \vec {k}\right)$.\ifcase \cras \or /\textit{{The orthonormal frame $\left(A,\vec {i}, \vec {j}\right)$.}}\fi} 
\end{figure}

On considère 
de nouveau le repère orthonormé  \ajoutU\ $\left(A, \vec {\alpha}, \vec {k}\right)$.
On note
de nouveau 
$(  x(s),  y(s))$ les coordonnées de $X$ dans ce repère,
comme indiqué sur la figure \ref{longueur_majoreeb}.
\ifcase \cras
\or
On a les relations habituelles  
\begin{subequations}
\label{eqderphinbbbbnewN}
\begin{align}
\label{eqderphinabnewN}
&\frac{d  x}{ds}=\cos   \phi,\\
\label{eqderphinbbnewN}
&\frac{d y}{ds}=\sin   \phi.
\end{align}
\end{subequations}
\fi
D'après les hypothèses 
\eqref{eq20},
\eqref{eq100d},
\eqref{eq100e},
\eqref{eq150} et 
\ajoutJJ
\label{ajoutJJ10}
il existe $s_0\in [0,L[$ tel
que 
\begin{equation}
\label{preuveq00}
\text{$\phi(s)=0$ sur \sout{$[0,s]$}\uwave{$[0,s_0]$} et  $0<\phi(s)< \pi$ sur $]s_0,L]$.}
\end{equation}
 Fixons $s\in ]s_0,L]$.
Pour tout point $P=(x_p,y_p)$ de la droite passant par $X(s)$ et porté par $X'(s)$, il existe $\lambda\in \Er$ tel que 
$\overrightarrow{PX(s)}=\lambda X'(s)$, soit d'après 
\ifcase \cras
\eqref{eqderphinbbbb}
\or
\eqref{eqderphinbbbbnewN}
\fi
\begin{equation}
\label{preuveq01}
x(s)-x_p=\lambda  \cos \phi(s),\quad y(s)-y_p=\lambda \sin\phi(s).
\end{equation}
Les deux droites respectives passant par $A$ et portée par $\vec \alpha$
et passant par $X(s)$ et portée par $X'(s)$
se coupent donc un unique point $Q(s)$ d'ordonnée $0$ et d'abscisse donnée par $x_p$
dans \eqref{preuveq01} correspondant à $y_p=0$. On a donc 
\begin{equation}
\label{preuveq10}
\lambda =\frac{y(s)}{\sin\phi(s)},\quad x_p=x(s)-y(s)\frac{ \cos \phi(s)}{\sin\phi(s)}
\end{equation}
Considérons $u(s)$, défini comme  l'abscisse $x_p$ de $Q(s)$ et $v(s)=\lambda$.
Le point $Q(s)$ vérifie donc 
\begin{equation}
\label{preuveq20}
\forall s \in ]s_0,L],\quad
\overrightarrow{AQ(s)}=u(s) \vec \alpha,\quad 
\overrightarrow{Q(s)X(s)}=v(s) \vec X'(s).
\end{equation}
où
\begin{equation}
\label{preuveq30}
u(s)= x(s)-y(s)\frac{ \cos \phi(s)}{\sin\phi(s)}  ,\quad 
v(s)=\frac{y(s)}{\sin\phi(s)} 
\end{equation}
On peut dériver $u$ presque partout et on a, compte tenu  de 
\ifcase \cras
\eqref{eqderphinbbbb}
\or
\eqref{eqderphinbbbbnewN}
\fi
\begin{align*}
\text{p.p. sur $]s_0,L[$, }
u'(s)&=x'(s)
-y'(s)\frac{ \cos  \phi(s)}{\sin\phi(s)} 
-y(s)\frac{ -\sin^2\phi(s)\phi'(s)-\cos^2\phi(s)\phi'(s)}{\sin^2\phi(s)} ,\\
&=\cos \phi(s)-\sin\phi(s)\frac{ \cos \phi(s)}{\sin\phi(s)}
+y(s)\frac{ \phi'(s)}{\sin^2\phi(s)} ,
\end{align*}
et donc 
\begin{equation}
\label{preuveq41}
u'(s)
=y(s)\frac{ \phi'(s)}{\sin^2\phi(s)}.
\end{equation}
Par ailleurs, d'après 
\ifcase \cras
\eqref{eqderphinbbbb},
\or
\eqref{eqderphinbbnewN}
\fi
on a  
\begin{equation*}
\text{$\displaystyle{y(s)=y(0)+\int_0^s y'(u)du=}$} \text{\sout{$y(0)+$}} \text{$\displaystyle{\int_0^s  \sin\phi(u)du.}$}
\end{equation*}
Ainsi, 
d'après 
\eqref{eq20},
\ajoutJJ
\label{ajoutJJ01}
 et 
\eqref{eq160}, 
$\sin(\phi)\geq 0$ et donc $y(s)\geq 0$, ce qui implique d'après \eqref{eq150} et \eqref{preuveq41} :
\begin{equation}
\label{preuveq50}
\text{$u$ est croissant sur $[s_0,L]$.}
\end{equation}
Ainsi, il existe 
\begin{equation}
\label{preuveq51} 
u(s_0+)=\lim_{\substack{s\to s_0\\ s>s_0}} u(s)\in \{-\infty\}\cap \Er.
\end{equation}
On a 
\begin{equation}
\label{preuveq60}
u(s_0+)\geq 0.
\end{equation}
\begin{figure}[h]    
\psfrag{f1}{$\phi_1$}
\psfrag{jt}{$\vec {k}$}
\psfrag{it}{$\vec {\alpha}$}
\psfrag{O}{$O$}
\psfrag{A}{$A$}
\psfrag{B}{$B$}
\psfrag{X}{$X(s_1)$}
\psfrag{Q}{$Q(s)$}
\psfrag{Xp}{$X'(s_1)$}
\psfrag{u}{$u_1$}
\psfrag{v}{$v_1$}
\begin{center} 
\epsfig{file=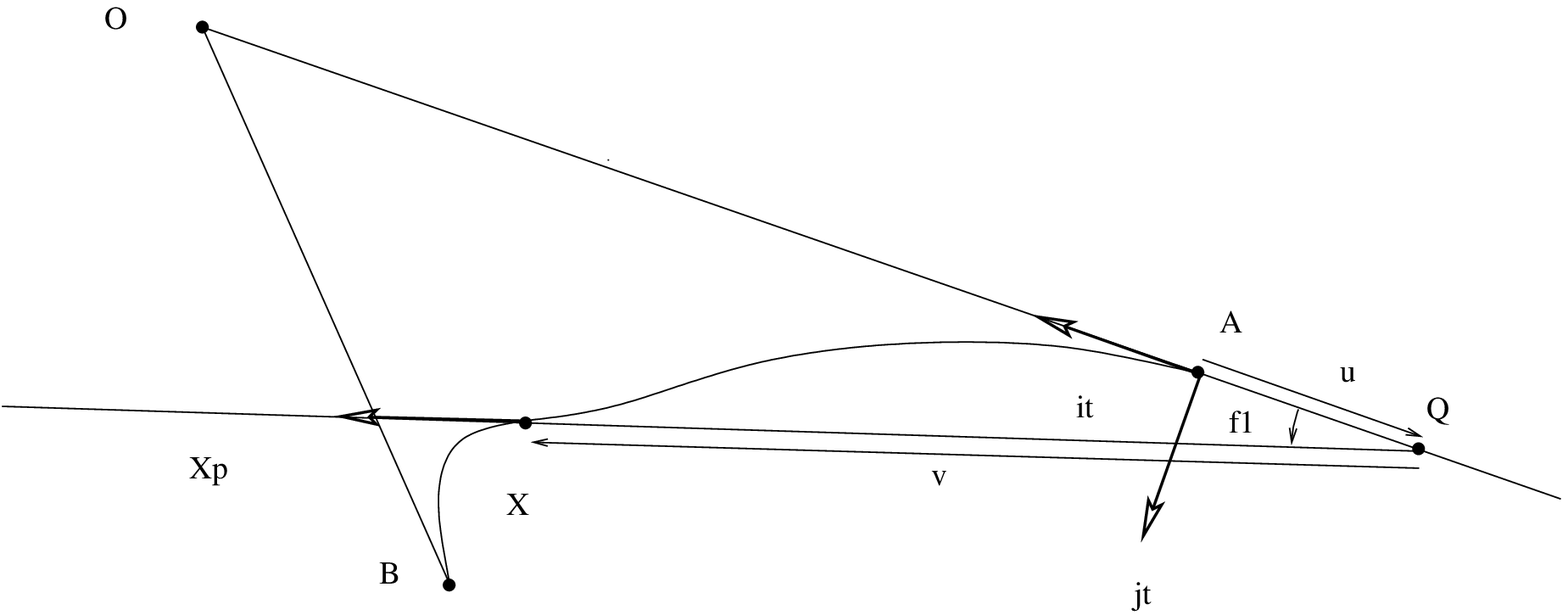, width=12 cm} 
\end{center} 
\caption{\label{longueur_majoree_bis}Le repère $\left(A,\vec {\alpha}, \vec {k}\right)$ avec un $u_1<0$.\ifcase \cras \or /\textit{{The orthonormal frame $\left(A,\vec {i}, \vec {j}\right)$} with $u_1<0$.}\fi} 
\end{figure}
Si ce n'était pas le cas, il existerait $s_1$ tel que 
\begin{equation}
\label{preuveq61}
\phi_1=\phi(s_1)\in ]0,\pi/2[,\quad 
u_1=u(s_1)<0.
\end{equation}
Cela implique que le point $A$ se situe dans le demi-plan ouvert limité par la droite passant par $X(s_1)$, portée par $X'(s_1)$, du côté opposé à la normale 
à la courbe en $X(s_1)$ (voir figure \ref{longueur_majoree_bis}).
Or, d'après le lemme \ref{lemmeconvexite}, la courbe doit aussi se trouver  dans le demi-plan  limité par la droite passant par $X(s_1)$, portée par $X'(s_1)$, du même côté que  la normale 
à la courbe en $X(s_1)$. Ainsi, \eqref{preuveq60} est vrai.
Rappelons que, d'après \eqref{eq100e}, il existe un sous-intervalle $J$ de $[s_0,L]$ tel que 
\begin{equation}
\label{preuveq50b}
\text{$u$ est strictement croissante sur $J$.}
\end{equation}
Ainsi, d'après \eqref{preuveq50}, \eqref{preuveq60} et \eqref{preuveq50b}
\begin{equation}
\label{preuveq70}
u(L)>0.
\end{equation}
Enfin, d'après \eqref{preuveq30}, $v(L)=y(L)/\sin\Omega>0$ et donc 
\begin{equation}
\label{preuveq80}
v(L)>0.
\end{equation}
On conclut en posant $u_0=u(L)$ et $v_0=v(L)$ et en utilisant \eqref{preuveq70} et \eqref{preuveq80} qui impliquent que $O$ est distinct de $A$ et de $B$.
\ifcase \cras
\end{proof}
\or
\end{preuve}
\fi


\ifcase \cras
\begin{proof}[Démonstration  du lemme \ref{uniquecercledroitelem01}]\
\or
\begin{preuve}[  du lemme \ref{uniquecercledroitelem01}]\
\fi

\label{tugudu90}

\ifcase \cras
\label{uniquecercledroitelem01preuve}
\fi

\begin{figure}[h]    
\psfrag{A}{$A$}
\psfrag{B}{$B$}
\psfrag{O}{$O$}
\psfrag{a}{$\alpha$}
\psfrag{b}{$\beta$}
\psfrag{d1}{$d_1$}
\psfrag{d2}{$d_2$}
\psfrag{C1}{$C_1$}
\psfrag{C2}{$C_2$}
\psfrag{c1}{$\mathcal{C}_1$}
\psfrag{c2}{$\mathcal{C}_2$}
\begin{center} 
\epsfig{file=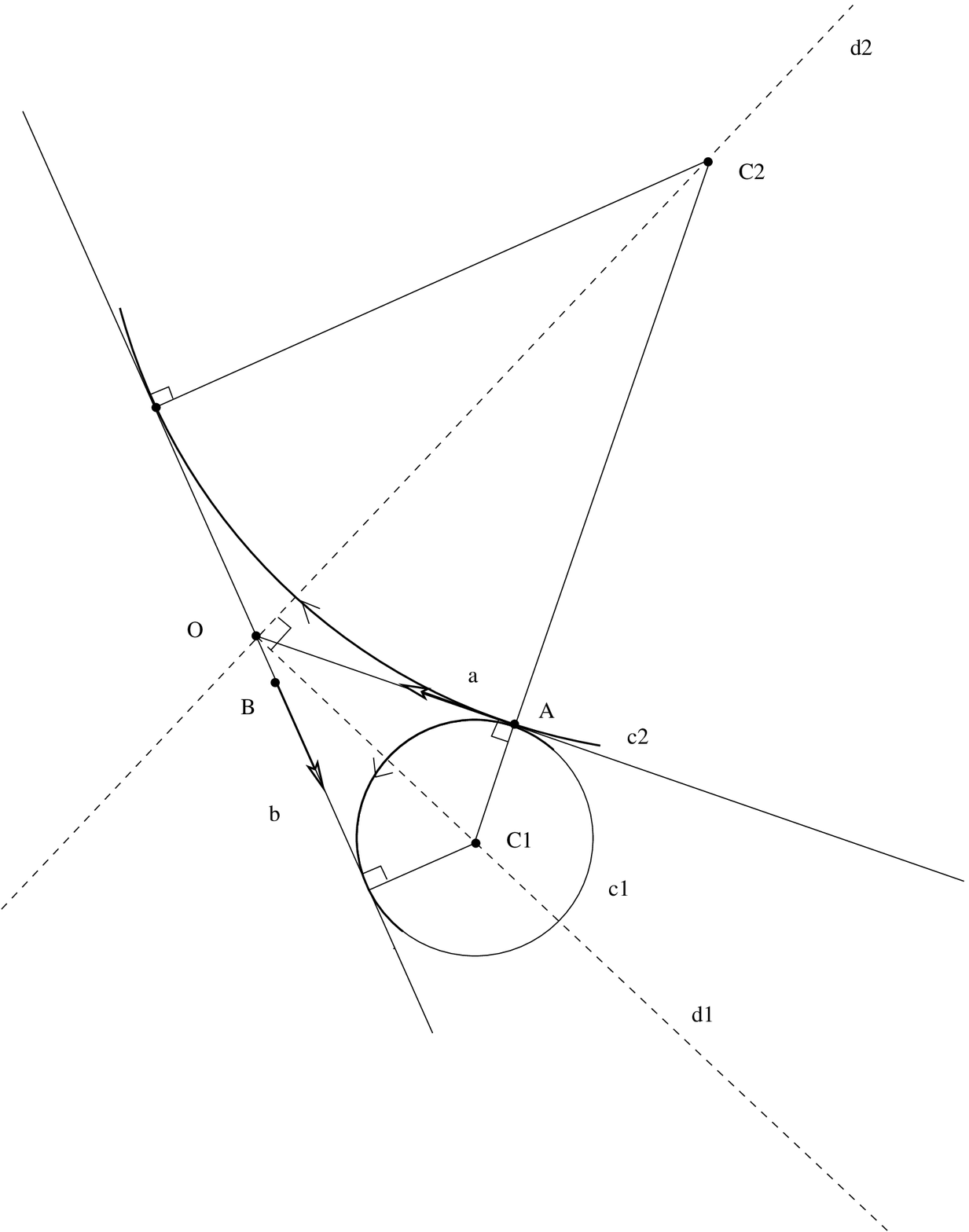, width= 10 cm} 
\end{center} 
\caption{\label{preuveuniquecercledroitelem01}$\mathcal{C}_1$ et $\mathcal{C}_2$ sont  tangents à la fois à $(OA)$ et à $(OB)$.\ifcase \cras \or /\textit{$\mathcal{C}_1$ and $\mathcal{C}_2$ are tangent to $(OA)$ and to  $(OB)$.}\fi}
\end{figure}
Supposons tout d'abord qu'il existe une courbe $Y$ de  $\mathcal{E}$
formée d'un arc de cercle $\mathcal{C}$ de longueur non nulle et d'un segment de droite $S$ éventuellement réduit à un point  et montrons qu'elle est nécessairement unique.
D'après \eqref{eq100tot}, l'arc de cercle est tangent en $A$ à $(OA)$ ou tangent en $B$ à $(OB)$. 

\begin{enumerate}

\item
Supposons  que 
$\mathcal{C}$  est tangent en $A$ à $(OA)$. Dans ce cas, le segment de droite  est  tangent en $B$ à $(OB)$ et il est donc inclus  dans la droite $(OB)$.
Puisque la courbe est de classe ${\mathcal{C}}^1$, $\mathcal{C}$ et $S$ sont tangents donc $\mathcal{C}$ est tangent à $(OB)$. 
Il n'existe que 
\ajoutW\
arcs de cercles possibles $\mathcal{C}_1$ et $\mathcal{C}_2$, dont les centres sont respectivement 
sur la bissectrice $d_1$ de $\widehat{AOB}$ ou la droite $d_2$ perpendiculaire à cette bissectrice, passant par $O$.
Voir figure \ref{preuveuniquecercledroitelem01}. Le centre $C_1$  (resp. $C_2$) $\mathcal{C}_1$ (resp. de  $\mathcal{C}_2$)   se trouve nécessairement sur la droite passant par $A$ et 
$d_1$ (resp. $d_2$), qui se coupent nécessairement, compte tenu de \eqref{eq20}.
Ces deux arcs de cercles, de part et d'autre de $(OA)$ doivent être parcourus dans le sens trigonométrique, pour respecter \eqref{eq140},
\ajoutX\
et, sur la figure  \ref{preuveuniquecercledroitelem01}, seul l'arc de cercle  $\mathcal{C}_1$ respecte ce sens.
Ainsi, si l'arc de cercle est tangent en $A$ à $(OA)$, il est unique et correspond à l'arc de cercle  $\mathcal{C}_1$ de la figure \ref{preuveuniquecercledroitelem01}. 
\begin{figure}[h]    
\psfrag{A}{$A$}
\psfrag{B}{$B$}
\psfrag{O}{$O$}
\psfrag{E}{$E$}
\psfrag{F}{$F$}
\psfrag{a}{$\vec \alpha$}
\psfrag{b}{$\vec \beta$}
\psfrag{Ca}{$\mathcal{C}_1$}
\psfrag{Cb}{$\mathcal{C}_3$}
\begin{center} 
\epsfig{file=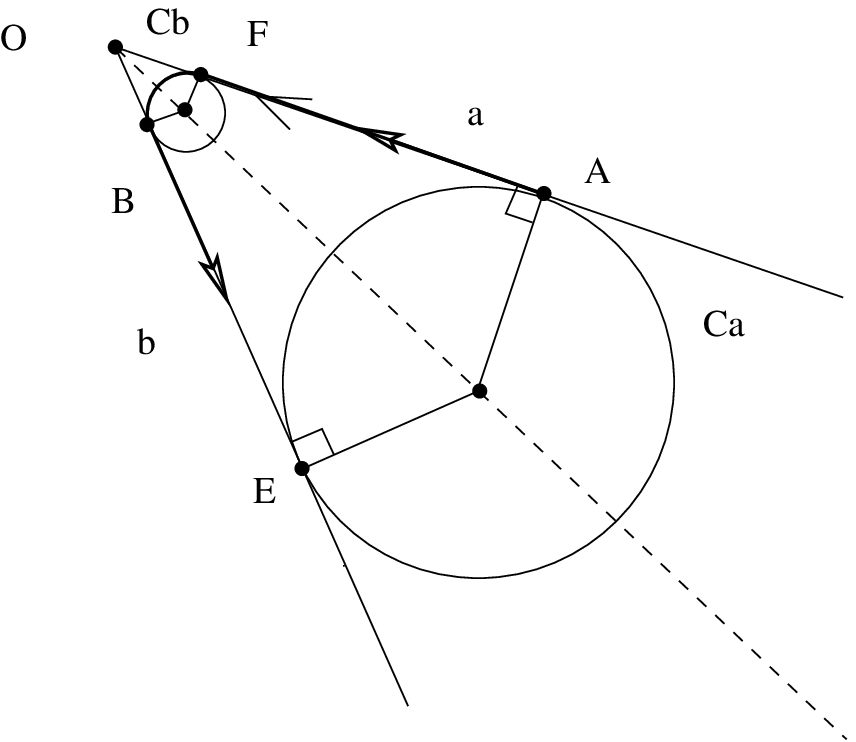, width= 6 cm} 
\end{center} 
\caption{\label{preuveuniquecercledroitelem02}Construction de $\mathcal{C}=\mathcal{C}_3$, cas 1.\ifcase \cras \or /\textit{Construction of $\mathcal{C}=\mathcal{C}_3$, case  1.}\fi}
\end{figure}
On considère ensuite le point de contact $E$ entre  $\mathcal{C}_1$, nécessairement sur la demi-droite $]OB)$. Voir figure \ref{preuveuniquecercledroitelem02}.
Si $E$ est strictement plus loin de $B$ que $O$ (ce qui implique $OB<OA$), alors la courbe $Y$ de  $\mathcal{E}$ ne peut revenir à $B$ par un segment de droite. Nécessairement,
$\mathcal{C}$  est tangent en $B$ à $(OB)$ et on construit l'unique arc de cercle $\mathcal{C}_3$, comme précédemment (voir figure \ref{preuveuniquecercledroitelem02}), tangent 
à $(OB)$ en $B$ et tangent à $(OA)$.
 Ainsi, dans ce cas, la courbe $Y$ est nécessairement formée de $\mathcal{C}=\mathcal{C}_3$ et du segment $[AF]$, non réduit à un point.
\begin{figure}[h]    
\psfrag{A}{$A$}
\psfrag{B}{$B$}
\psfrag{O}{$O$}
\psfrag{E}{$E$}
\psfrag{F}{$F$}
\psfrag{a}{$\vec \alpha$}
\psfrag{b}{$\vec \beta$}
\psfrag{Ca}{$\mathcal{C}_1$}
\begin{center} 
\epsfig{file=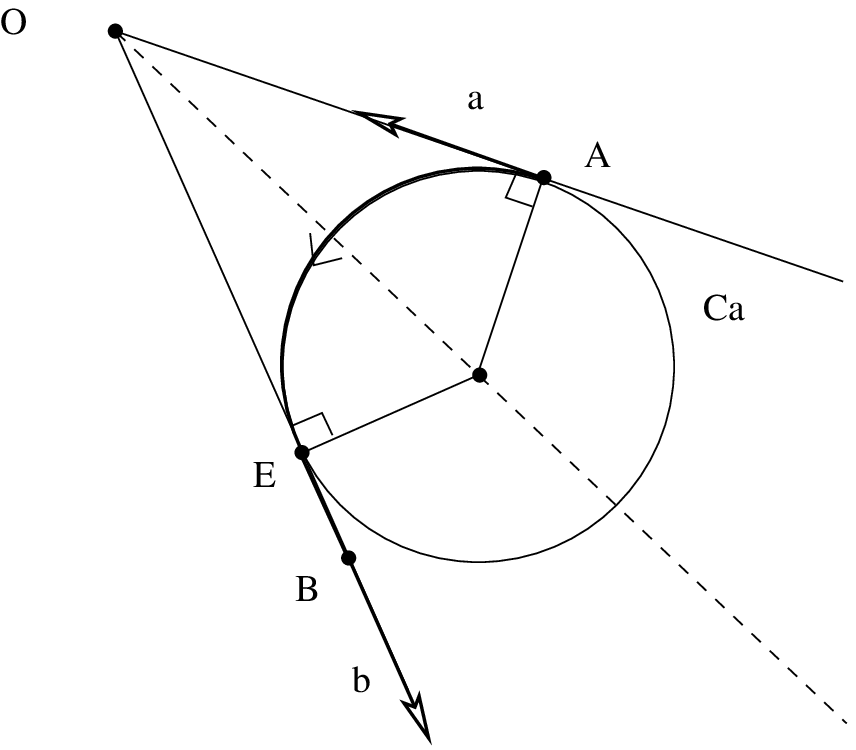, width= 6 cm} 
\end{center} 
\caption{\label{preuveuniquecercledroitelem03}Construction de $\mathcal{C}=\mathcal{C}_1$, cas 2.\ifcase \cras \or /\textit{Construction of $\mathcal{C}=\mathcal{C}_1$, case  2.}\fi}
\end{figure}
Si, au contraire $E$ est strictement moins loin de $B$ que $O$ (ce qui implique $OB>OA$), alors, de même,  la courbe $Y$ est nécessairement formée de $\mathcal{C}=\mathcal{C}_1$ et du segment $[EB]$, non réduit à un point (voir figure \ref{preuveuniquecercledroitelem03}).
\begin{figure}[h]    
\psfrag{A}{$A$}
\psfrag{B}{$B$}
\psfrag{O}{$O$}
\psfrag{a}{$\vec \alpha$}
\psfrag{b}{$\vec \beta$}
\psfrag{Ca}{$\mathcal{C}_1$}
\begin{center} 
\epsfig{file=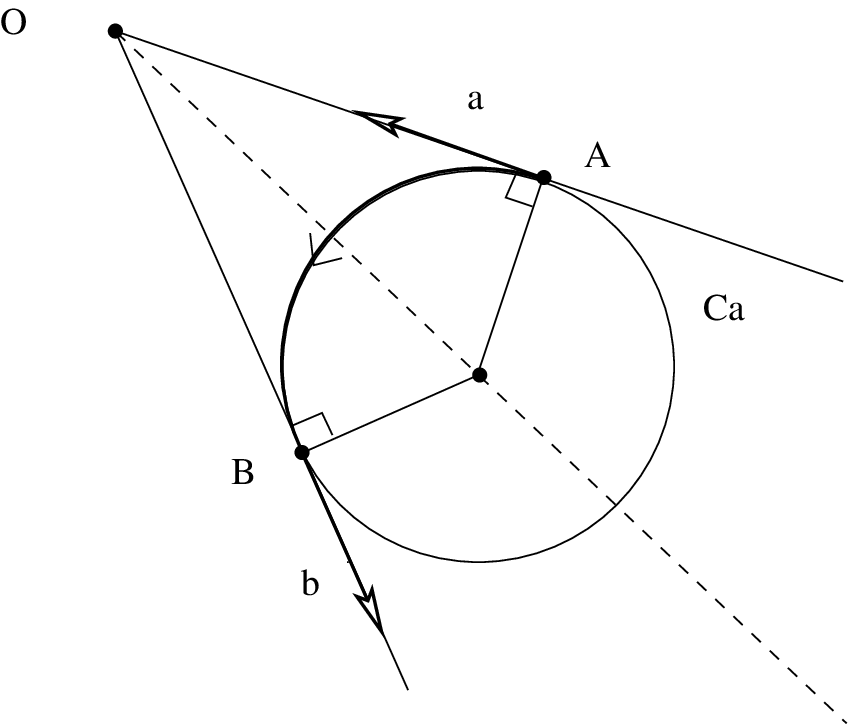, width= 6 cm} 
\end{center} 
\caption{\label{preuveuniquecercledroitelem04}Construction de $\mathcal{C}=\mathcal{C}_1$, cas 3 (cas symétrique).\ifcase \cras \or /\textit{Construction of $\mathcal{C}=\mathcal{C}_1$, case  1, symmetric  case.}\fi}
\end{figure}
Enfin, si $E$ est exactement aussi loin de $B$ que $O$ (ce qui implique $OB=OA$), alors, de même,  la courbe $Y$ est nécessairement uniquement formée de $\mathcal{C}=\mathcal{C}_1$  (voir figure \ref{preuveuniquecercledroitelem04}).
S'il existe une courbe $Y$ de  $\mathcal{E}$
formée d'un arc de cercle $\mathcal{C}$ de longueur non nulle, alors on tombe sur ce dernier cas.

\item
Si 
$\mathcal{C}$  est tangent en $B$ à $(OB)$, on arrive à la même construction. 
\end{enumerate}

Montrons maintenant l'existence de la courbe.
Dans les trois cas évoqués dans l'unicité (selon que $OB>OA$, $OB<OA$ ou $OA=OB$), on peut construire, une courbe 
$Y$ de $\mathcal{E}$,
formée d'un arc de cercle et de
périmètre dans $]0,R_a(O,A,B) \pi [$  et d'un segment de droite. 
Dans cette construction, on vérifie que 
\eqref{eq100b},
\eqref{eq100c},
\eqref{eq100d}
et 
\eqref{eq100e}
on lieu  grâce à \eqref{eq10tottot}.

L'égalité \eqref{uniquecercledroitelem01eq01} est triviale, puisque 
\sout{$\vnorm{X''}=0$ ou }
\uwave{$\vnorm{X''}=$ est égal à $0$ ou à}
$1/R_a$.

\ifcase \cras
\end{proof}
\or
\end{preuve}
\fi


\ifcase \cras
\begin{proof}[Démonstration  du lemme \ref{optdubinslem01}]\
\or
\begin{preuve}[ du lemme \ref{optdubinslem01}]\
\fi

\label{tugudu100}

\ifcase \cras
\label{optdubinslem01preuve}
\fi

Soit $R_a=R_a(O,A,B)$, le nombre défini dans le lemme \ref{uniquecercledroitelem01}.
Supposons par exemple que $OB<OA$. On est donc dans le premier cas de la démonstration du lemme \ref{uniquecercledroitelem01}.
La courbe de Dubins  $\mathcal{G}(R)$ pour $R<R_a(O,A,B)$ est unique et est formée alors de deux arcs de cercles,
reliés par un segment de droite $[ED]$ avec $E\not =D$, comme le montre la figure 
\refpreuve.\ref{codu0bis}.
Le premier arc de cercle est nécessairement tangent à la droite $(OB)$ en $B$ et son centre $F$, est tel que $(BF)$ est 
$(OB)$ soient perpendiculaires. Puisque $R<R_a(O,A,B)$, $F$ est dans le secteur de plan défini par les demi-droites
$[OB)$ et $d$, la bissectrice de l'angle $\widehat{AOB}$. Si $G$ est le centre du second cercle de la courbe
de Dubins, tangent en $A$ à $(OA)$, alors $FGDE$ est un rectangle et le segment $[ED]$ est inclus 
dans le secteur de plan défini par les demi-droites
$[OB)$ et $[OA)$. Le second cercle de la courbe de Dubins est de longueur non nulle. On vérifie que la courbe $\mathcal{G}(R)$ appartient bien à $\mathcal{E}$. Cela est vrai tant que 
$R$ appartient à $]0,R_a(O,A,B)[$. Si $R$ croît, $ED$ augmente, la longueur du premier arc de cercle augmente, celle du second diminue.
Pour le cas limite, $R=R_a(O,A,B)$, correspondant  à la figure  
\refpreuve.\ref{coducl1bis},
le premier cercle devient tangent à $(OA)$ en $E$, $D$ se confond avec $A$ et la longueur du second cercle
devient nulle. On retrouve donc l'unique courbe du lemme \ref{uniquecercledroitelem01}.
\newcommand{\lengthdiftypcodusc}{7}
\begin{figure}
\psfrag{A}{$A$}
\psfrag{B}{$B$}
\psfrag{C}{$C$}
\psfrag{D}{$D$}
\psfrag{O}{$O$}
\psfrag{D}{$D$}
\psfrag{E}{$E$}
\psfrag{F}{$F$}
\psfrag{G}{$G$}
\psfrag{a}{$\vec \alpha$}
\psfrag{b}{$\vec \beta$}
\centering
\subfigure[\label{codu0bis}Le cas 0 : $0<R<R_a(O,A,B)$\ifcase \cras \or /\textit{{The case 0 : $0<R<R_a(O,A,B)$.}}\fi]
{\epsfig{file=codu0, width=\lengthdiftypcodusc  cm}}
\quad
\subfigure[\label{coducl1bis}Le cas limite 0-1 : $R=R_a(O,A,B)$\ifcase \cras \or /\textit{{The limit case 0-1 : $R=R_a(O,A,B)$.}}\fi]
{\epsfig{file=coducl1, width=\lengthdiftypcodusc   cm}}
\quad
\subfigure[\label{codu2}Le cas  1 : $R_a(O,A,B)<R<R_b(O,A,B)$\ifcase \cras \or /\textit{{The case 1 : $R_a(O,A,B)<R<R_b(O,A,B)$.}}\fi]
{\epsfig{file=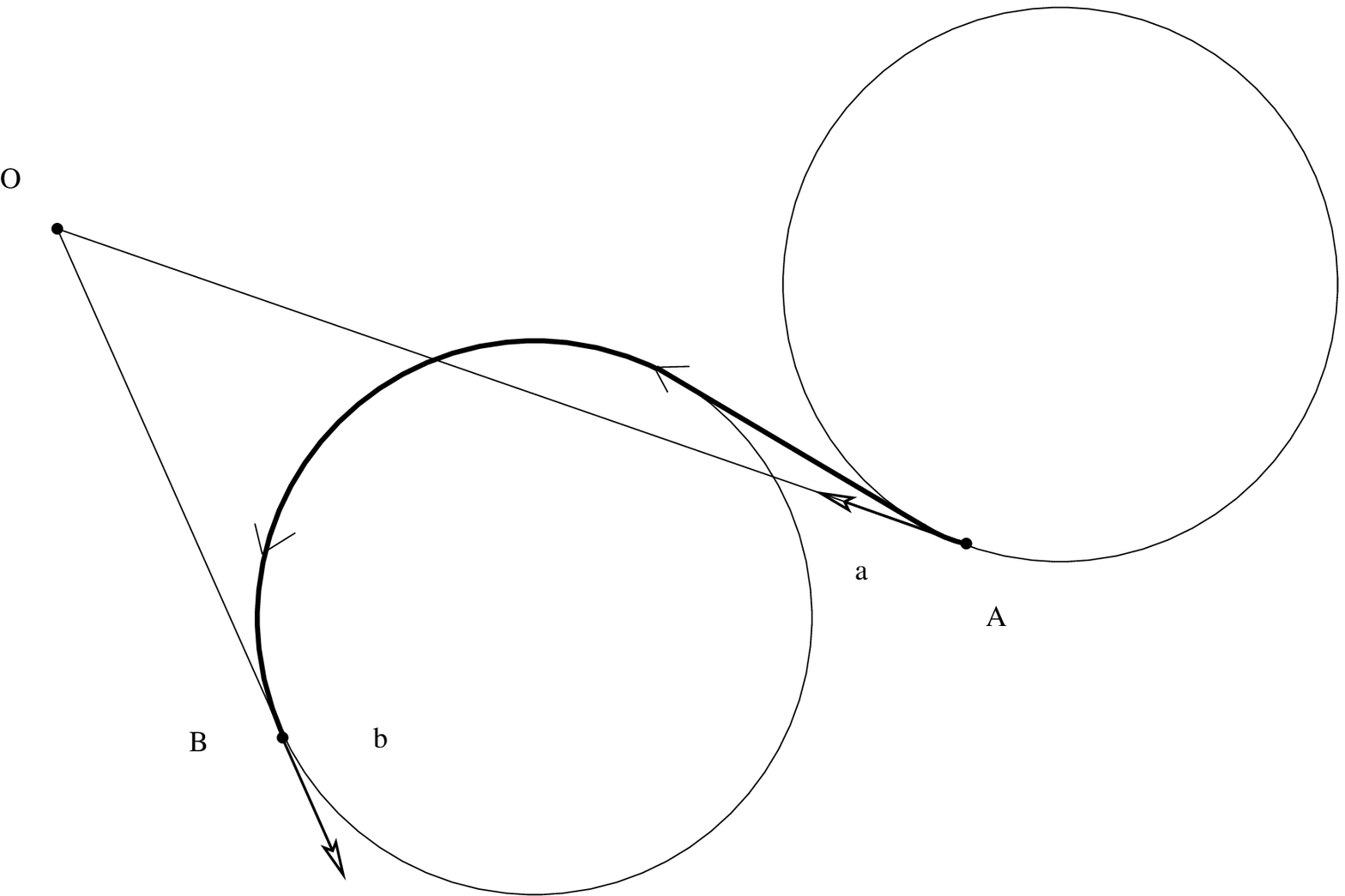, width=\lengthdiftypcodusc   cm}}
\quad
\subfigure[\label{coducl2}Le cas limite 1-2 : $R=R_b(O,A,B)$\ifcase \cras \or /\textit{{The limit case 1-2 : $R=R_b(O,A,B)$.}}\fi]
{\epsfig{file=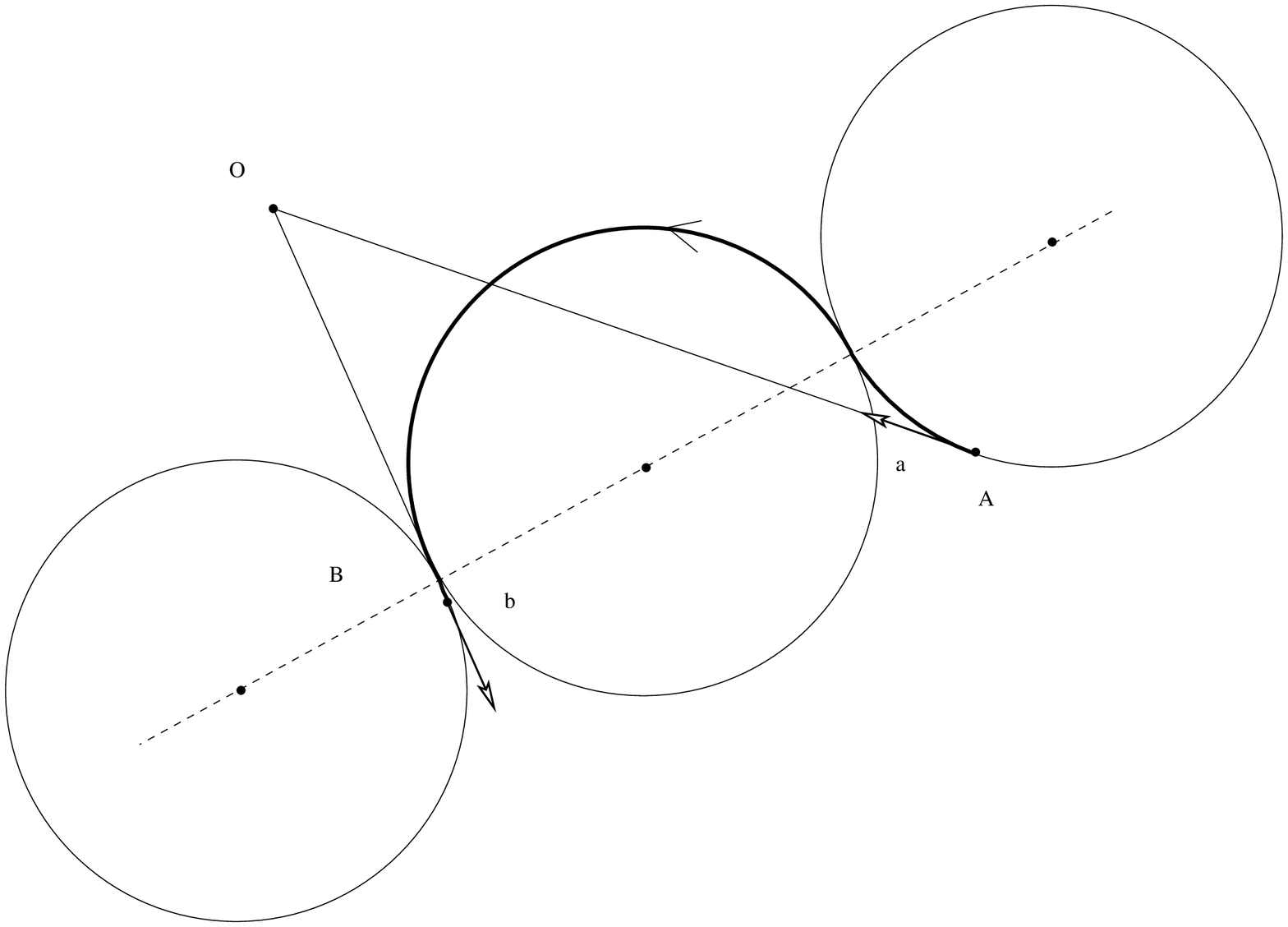, width=\lengthdiftypcodusc   cm}}
\quad
\subfigure[\label{codu3}Le cas 2 : $R_b(O,A,B)<R<R_c(O,A,B)$ \ifcase \cras \or /\textit{{The case 2 : $R_b(O,A,B)<R<R_c(O,A,B)$.}}\fi]
{\epsfig{file=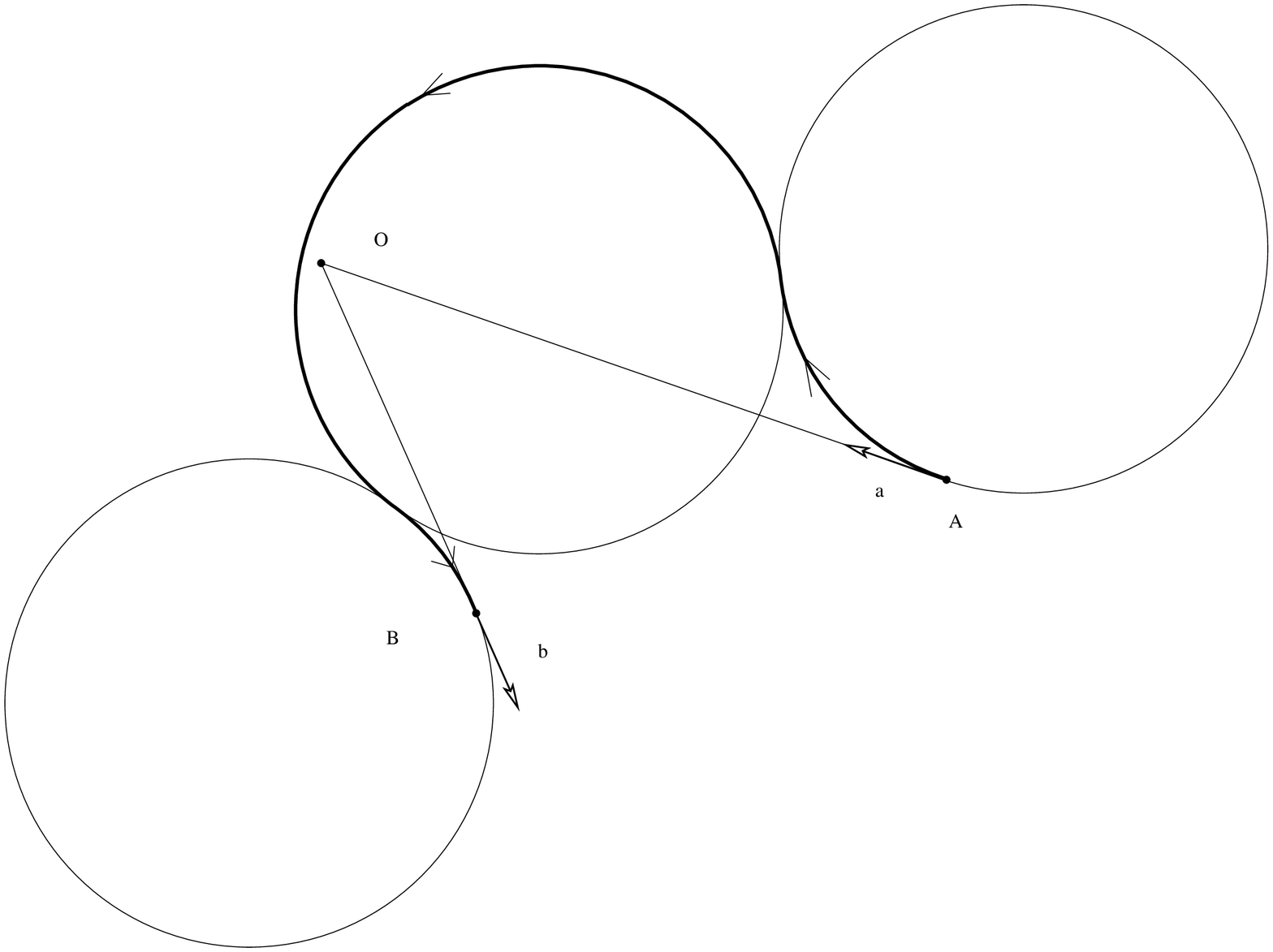, width=\lengthdiftypcodusc   cm}}
\subfigure[\label{codu4}Le cas 3 : $R>R_c(O,A,B)$ \ifcase \cras \or /\textit{{The case 2 : $R>R_c(O,A,B)$.}}\fi]
{\epsfig{file=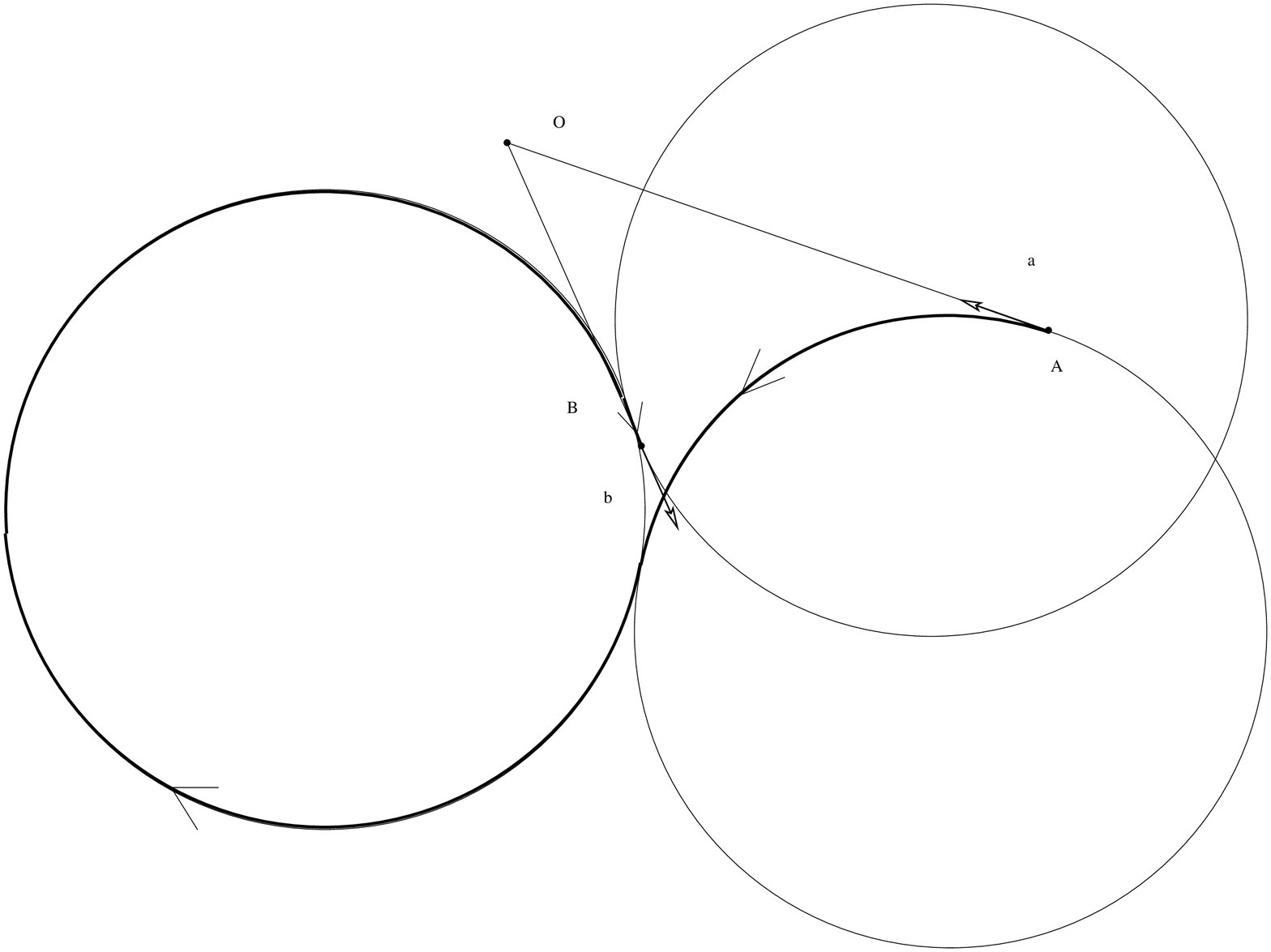, width=\lengthdiftypcodusc   cm}}
\caption{\label{diftypcodubis}Les différents types de courbes de Dubins définies par $R$  dans le cas non symétrique. 
$R_a(O,A,B)$, 
$R_b(O,A,B)$
et 
$R_c(O,A,B)$ ne dépendent que de $O$, $A$ et $B$. \ifcase \cras \or /\textit{{The different kinds of Dubins's Curve defined by $R$ in the non symmetric case. 
$R_a(O,A,B)$, 
$R_b(O,A,B)$
and
$R_c(O,A,B)$  depend only on $O$, $A$ and  $B$.}}\fi}
\end{figure}
\begin{figure}
\psfrag{A}{$A$}
\psfrag{B}{$B$}
\psfrag{C}{$C$}
\psfrag{D}{$D$}
\psfrag{O}{$O$}
\psfrag{a}{$\vec \alpha$}
\psfrag{b}{$\vec \beta$}
\centering
\subfigure[\label{codu0sc}Le cas 0 : $0<R<R_a(O,A,B)$.\ifcase \cras \or /\textit{{The case 0 : $0<R<R_a(O,A,B)$.}}\fi]
{\epsfig{file=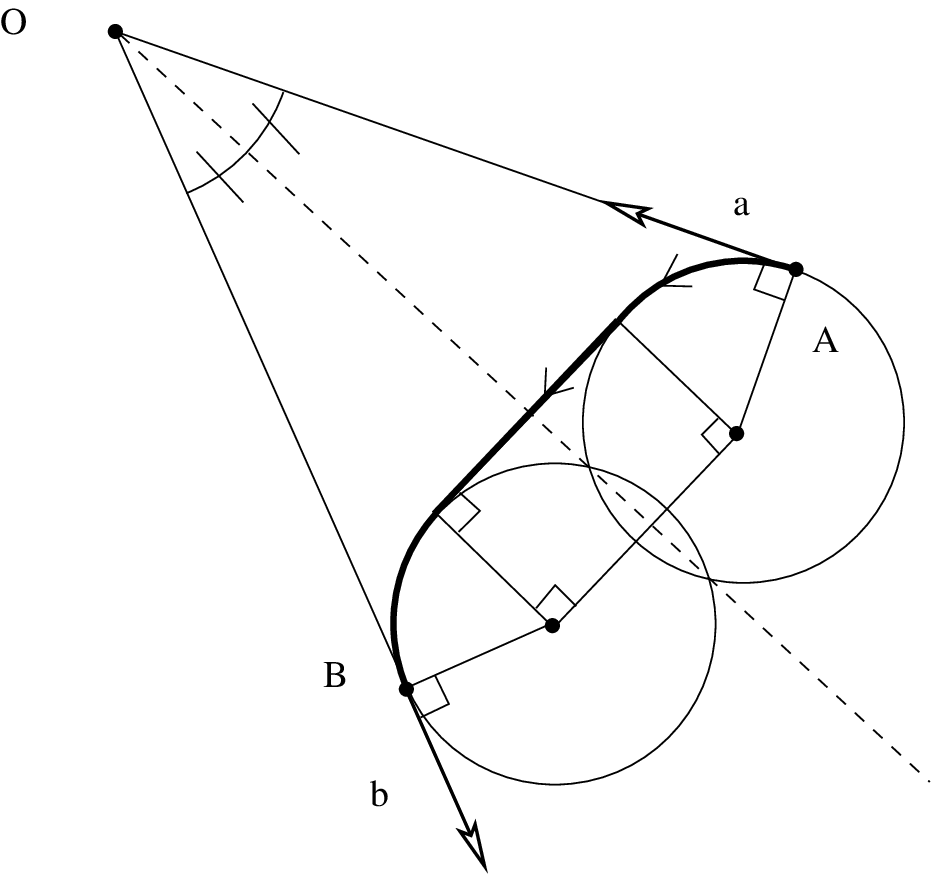, width=\lengthdiftypcodusc cm}}
\quad
\subfigure[\label{coducl1sc}Le cas limite 0-1 : $R=R_a(O,A,B)$.\ifcase \cras \or /\textit{{The limit case 0-1 : $R=R_a(O,A,B)$.}}\fi]
{\epsfig{file=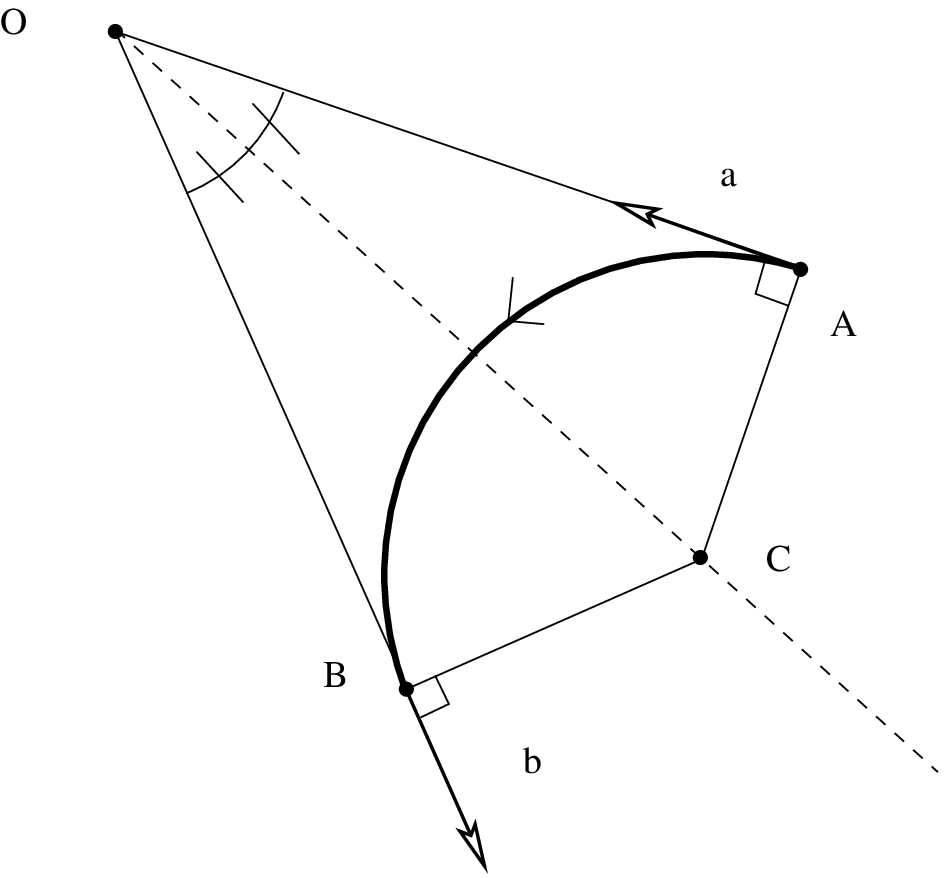, width=\lengthdiftypcodusc  cm}}
\quad
\subfigure[\label{codu3sc}Le cas  1 : $R_a(O,A,B)<R<R_b(O,A,B)$. \ifcase \cras \or /\textit{{The case 1 : $R_a(O,A,B)<R<R_b(O,A,B)$.}}\fi]
{\epsfig{file=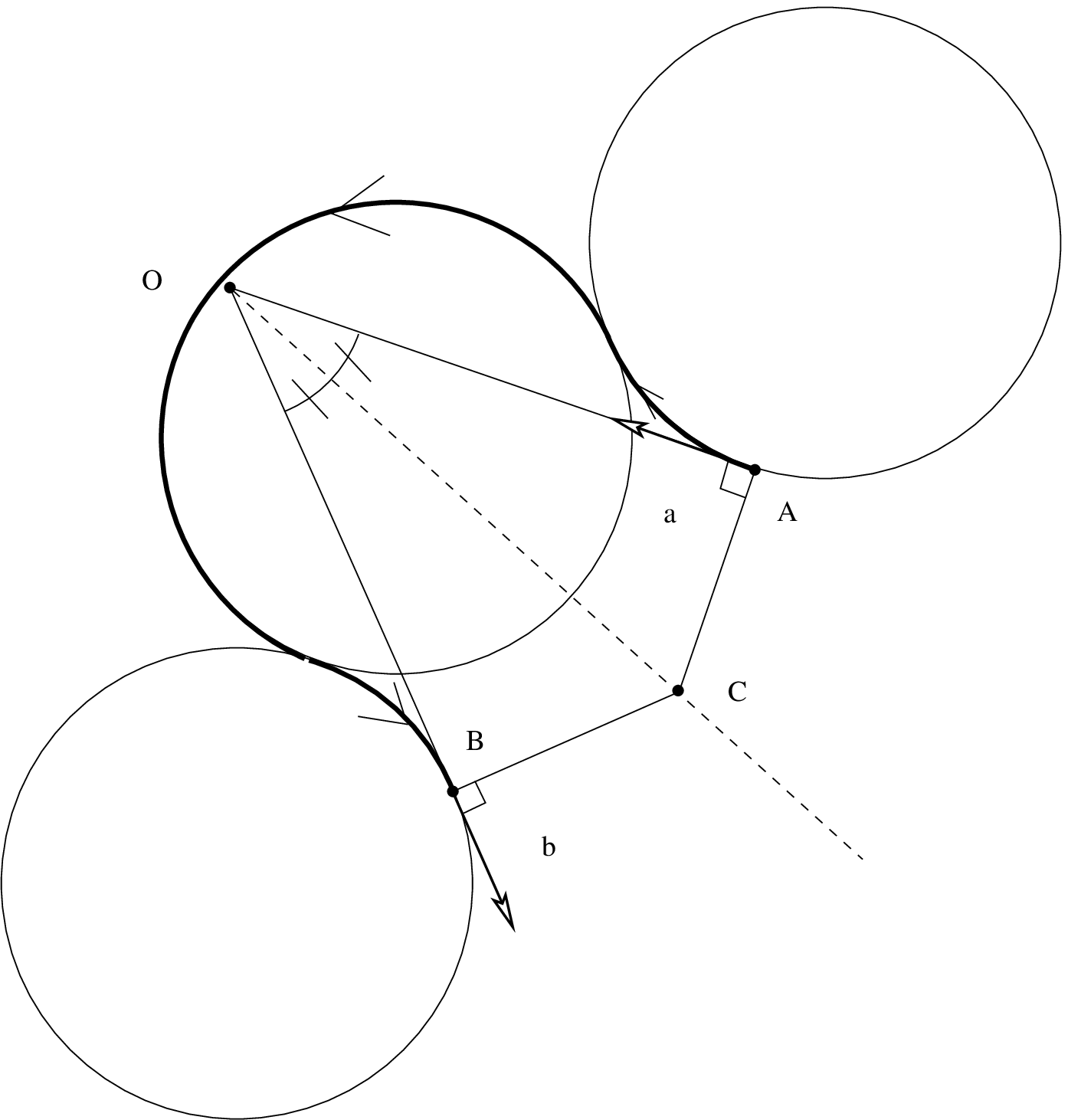, width=\lengthdiftypcodusc  cm}}
\subfigure[\label{codu2sc}Le cas  2 : $R_b(O,A,B)>R$. \ifcase \cras \or /\textit{{The case 1 : $R>R_b(O,A,B)$.}}\fi]
{\epsfig{file=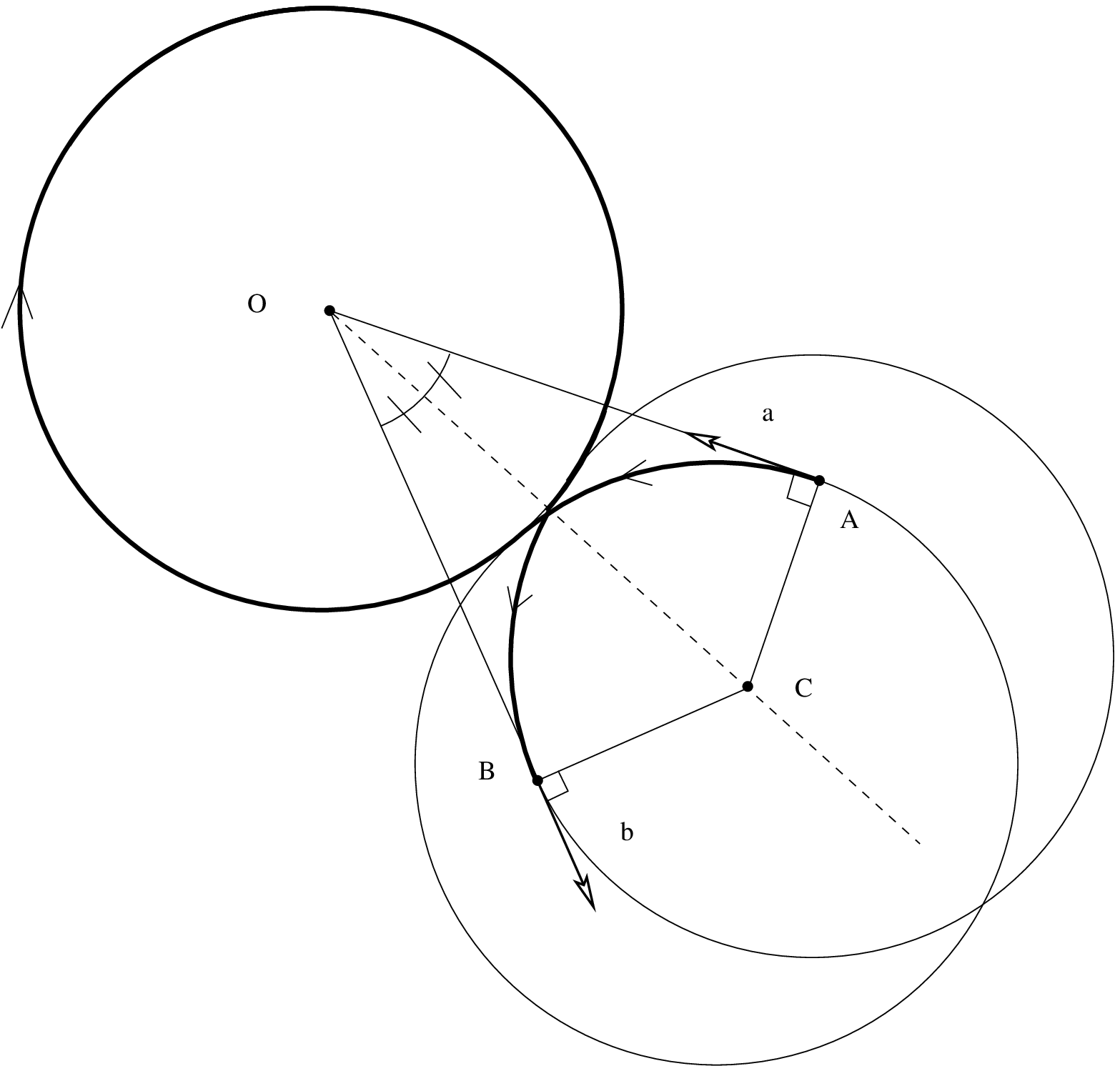, width=\lengthdiftypcodusc  cm}}
\caption{\label{diftypcodusc}Les différents types de courbes de Dubins définies par $R$  dans le cas symétrique.
$R_a(O,A,B)$
et 
$R_b(O,A,B)$ ne dépendent que de $O$, $A$ et $B$
\ifcase \cras \or /\textit{{The different kinds of Dubins's Curve defined by $R$ in the symmetric case.
$R_a(O,A,B)$
and
$R_c(O,A,B)$  depend only on $O$, $A$ and  $B$.}}\fi}
\end{figure}

Les cas $OA<OB$ ou $OA=OB$ se traitent de la même façon (voir figures
\refpreuve.\ref{codu0sc} et 
\refpreuve.\ref{coducl1sc}).

En revanche, on laisse au lecteur le soin de vérifier que 
  dès que $R$ dépasse strictement $R_a(O,A,B)$ 
on obtient des courbes qui présentent un changement strict du signe de la courbure. 
Voir figures \ref{diftypcodubis} et \ref{diftypcodusc}.
Le cas limite des  figures
\refpreuve.\ref{coducl1bis} 
et 
\refpreuve.\ref{coducl1sc} 
correspond exactement à la courbe décrite dans le lemme \ref{uniquecercledroitelem01}.
\ifcase \cras
\end{proof} 
\or
\end{preuve} 
\fi

\clearpage

\ifcase \cras
\section{Preuves complètes des résultats essentiels de cet article}
\or
\section{Preuves complètes des résultats essentiels de cette Note}
\fi
\label{preuvecomplete}

\ifcase \cras
\begin{proof}[Démonstration  de la proposition \ref{lemme100}]\
\or
\begin{preuve}[ de la proposition \ref{lemme100}]\
\fi

\label{tugudu00}

\begin{enumerate}
\item
Puisque $Z$ appartient à $\mathcal{E}$, on a par définition de $e$ :  
\begin{equation*}
\Omega=\theta(M)=\int_0^M \frac{d\theta}{ds} ds
=\int_0^M \left|\frac{d\theta}{ds}\right| ds
\leq eM,
\end{equation*}
et on a donc, d'après \eqref{lemme100eq10},  $M\geq \Omega/e\geq  \Omega R_a$, soit
\ajoutAE
\begin{equation}
\label{lemme100eq40a}
M \geq l.
\end{equation}
Les deux fonctions $X$ et $Z$ sont donc ainsi définies toutes les deux aux moins sur $[0,l]$ et il est légitime de poser :
\begin{equation}
\label{lemme100eq50}
\forall s\in [0,l],\quad
\zeta(s)=
-\left(\widehat x(s) -x(s)\right)\sin\left(\phi(s)\right)
+
\left(\widehat y(s)-  y(s) \right)\cos\left(\phi(s)\right).
\end{equation}
Ainsi, 
géométriquement, $\zeta$ correspond à 
\begin{equation*}
\zeta(s)=\prodsca{\begin{pmatrix}-\sin(\phi(s))\\\cos(\phi(s)\end{pmatrix}}{\overrightarrow{X(s)Z(s)}},
\end{equation*}
où, sur la partie circulaire de la courbe $X$ : 
\begin{equation*}
\begin{pmatrix}-\sin(\phi(s))\\\cos(\phi(s)\end{pmatrix}=\sigma(X'(s))=\vec N(s),
\end{equation*}
où $\sigma$ est la rotation vectorielle d'angle $\pi/2$ et $\vec N(s)$ désigne la normale extérieure à la courbe $X$ au point d'abscisse curviligne $S$.
Ainsi, $\zeta(s)$ désigne la composante du vecteur $\overrightarrow{X(s)Z(s)}$ sur $\vec N(s)$. 
\begin{figure}[h]    
\psfrag{a}{$Z(s)$}
\psfrag{c}{$X(s)$}
\psfrag{b}{$\vec X'(s)$}
\psfrag{d}{$\vec  N(s)$}
\begin{center} 
\epsfig{file=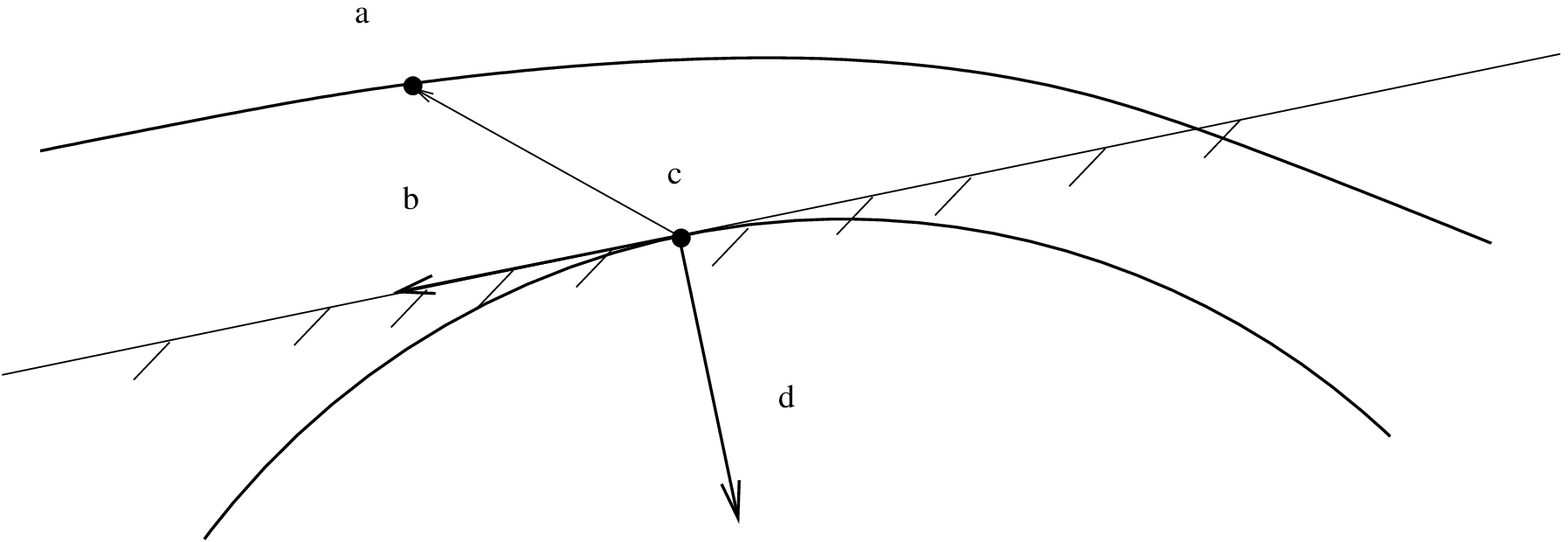, width=9 cm}  
\end{center} 
\caption{\label{composante_normale_courbe}$Z$ est toujours dans le demi-plan
défini par la tangente à courbe $X$, du côté opposé à la normale. 
\ifcase \cras \or \textit{
$Z(s)$ is always in the half-plane defined by the tangent of curve $X$ on the opposite side to the  outer-pointing normal.}\fi}
\end{figure}
Voir figure \ref{composante_normale_courbe}.
Nous allons étudier le signe de $\zeta(s)$ et montrer sous l'hypothèse \eqref{lemme100eq10} que $\zeta(s)$ est négatif,
ce qui signifie que la courbe $Z$ est toujours dans le demi-plan
défini par la tangente à courbe $X$, du côté opposé à la normale.

\item
Par définition de $\phi$, on a 
\begin{equation}
\label{lemme100eq100}
\forall s\in[0,l[,\quad
\phi'(s)=\frac{1}{R_a}.
\end{equation}
\'Etudions tous d'abord la fonction $\theta-\phi$ sur $[0,l]$. 
On a donc, presque partout sur $[0,l]$  : 
\begin{equation*}
\theta'(s)-\phi'(s)=
             \theta'(s)-\frac{1}{R_a}=
{\vnorm{Z''(s)}}-\frac{1}{R_a}
\leq 
\vnorm[L^{\infty}(0,L)]{\vnorm{Z''}}-\frac{1}{R_a}
=
e-\frac{1}{R_a},
\end{equation*}
et donc
\begin{equation*}
\forall s\in [0,l],\quad
\theta(s)-\phi(s)
\leq \theta(0)-\phi(0)+\left(e-\frac{1}{R_a}\right)s,
\end{equation*}
et donc
\begin{equation}
\label{lemme100eq110}
\forall s\in]0,l],\quad
\theta(s)-\phi(s)
\leq \left(e-\frac{1}{R_a}\right)s.
\end{equation}

\item
On a les relations habituelles  
\begin{subequations}
\label{eqderphinbbbbcrasnew}
\begin{align}
\label{eqderphinabcrasnew}
&\frac{d  \widehat x}{ds}=\cos   \theta,\\
\label{eqderphinbbcrasnew}
&\frac{d \widehat y}{ds}=\sin   \theta,\\
\label{eqderphinabnewcrasnew}
&\frac{d   x}{ds}=\cos   \phi,\\
\label{eqderphinbbnewcrasnew}
&\frac{d  y}{ds}=\sin   \phi,
\end{align}
\end{subequations}
On a donc, 
\begin{equation*}
\forall s\in [0,l],\quad
\frac{d}{ds}\left(\widehat x(s) -x(s)\right)
=\cos\left(\theta(s)\right)-\cos\left(\phi(s)\right),
\end{equation*}
et donc
\begin{subequations}
\label{lemme100eq120tot}
\begin{align}
\label{lemme100eq120a}
&\forall s\in[0,l],\quad
\frac{d}{ds}\left(\widehat x(s) -x(s)\right)
=
-2\sin\left(\frac{\theta(s)+\phi(s)}{2}\right)\sin\left(\frac{\theta(s)-\phi(s)}{2}\right).
\intertext{%
On a de même} 
\label{lemme100eq120b}
&\forall s\in[0,l],\quad
\frac{d}{ds}\left(\widehat y(s) -y(s)\right)
=
2\cos\left(\frac{\theta(s)+\phi(s)}{2}\right)\sin\left(\frac{\theta(s)-\phi(s)}{2}\right).
\end{align}
\end{subequations}
Par définition de $ \phi$, on a 
\begin{equation}
\label{lemme100eq130}
\forall s\in]0,l],\quad
0< \phi (s)<\pi.
\end{equation}
Par ailleurs, d'après les hypothèses 
\eqref{eq20},
\eqref{eq150} et 
\ajoutJJ
\label{ajoutJJ30}
\begin{equation*}
\forall s\in[0,l],\quad
0\leq  \frac{\theta(s)+\phi(s)}{2}\leq \Omega<\pi,
\end{equation*}
et puisque, $\phi$ est strictement croissant, on a 
\begin{subequations}
\label{lemme100eq140tot}
\begin{align}
\label{lemme100eq140a}
& \forall s\in]0,l],\quad
0< \frac{\theta(s)+\phi(s)}{2} <\pi.
\intertext{%
On a, de même, en utilisant \eqref{lemme100eq110}}
\label{lemme100eq150}
&\forall s\in]0,l],\quad
-\pi< \frac{\theta(s)-\phi(s)}{2} \leq \left(e-\frac{1}{R_a}\right)s.
\end{align}
\end{subequations}
De \eqref{lemme100eq120a} et \eqref{lemme100eq140tot}, on déduit
\begin{subequations}
\label{lemme100eq150totante}
\begin{align}
\label{lemme100eq150aante}
 \forall s\in]0,l],\quad
&\frac{d}{ds}\left(\widehat x(s) -x(s)\right) \geq 0,
\intertext{et}
\label{lemme100eq150bante}
e<\frac{1}{R_a}
\Longrightarrow 
 \forall s\in]0,l],\quad
&\frac{d}{ds}\left(\widehat x(s) -x(s)\right)>0,
\end{align}
\end{subequations}
et donc 
\begin{subequations}
\label{lemme100eq150tot}
\begin{align}
\label{lemme100eq150a}
 \forall s\in]0,l],\quad
&\widehat x(s) -x(s) \geq 0,
\intertext{et}
\label{lemme100eq150b}
e<\frac{1}{R_a}
\Longrightarrow 
 \forall s\in]0,l],\quad
&\widehat x(s) -x(s)>0.
\end{align}
\end{subequations}
On raisonne de la même façon avec les ordonnées, en remplaçant \eqref{lemme100eq140a}
par la majoration plus fine : 
\begin{equation}
\label{lemme100eq180}
\forall s\in]0,l],\quad
0< \frac{\theta(s)+\phi(s)}{2} <\frac{\pi}{2},
\end{equation}
qui provient de $ \frac{\theta(s)+\phi(s)}{2}\leq \Omega$
et de \eqref{eq20bis}.
On a alors
\begin{subequations}
\label{lemme100eq180totante}
\begin{align}
\label{lemme100eq180aante}
 \forall s\in]0,l],\quad
&\frac{d}{ds}\left(\widehat y(s) -y(s)\right) \leq 0,
\intertext{et}
\label{lemme100eq180bante}
e<\frac{1}{R_a}
\Longrightarrow 
 \forall s\in]0,l],\quad
&\frac{d}{ds}\left(\widehat y(s) -y(s)\right)<0,
\end{align}
\end{subequations}
et donc 
\begin{subequations}
\label{lemme100eq180tot}
\begin{align}
\label{lemme100eq180a}
 \forall s\in]0,l],\quad
&\widehat y(s) -y(s) \leq 0,
\intertext{et}
\label{lemme100eq180b}
e<\frac{1}{R_a}
\Longrightarrow 
 \forall s\in]0,l],\quad
&\widehat y(s) -y(s)<0.
\end{align}
\end{subequations}
Sous l'hypothèse \eqref{eq20bis}, on remplace alors \eqref{lemme100eq130} par la majoration plus fine : 
\begin{equation}
\label{lemme100eq190}
\forall s\in]0,l],\quad
0< \phi (s)<\frac{\pi}{2}.
\end{equation}

\item
Bref, grâce à 
\eqref{lemme100eq150tot},
\eqref{lemme100eq180tot} et 
\eqref{lemme100eq190}
et la définition \eqref{lemme100eq50} de $\zeta$, on a 
\begin{subequations}
\label{lemme100eq190tot}
\begin{align}
\label{lemme100eq190a}
 \forall s\in]0,l],\quad
&\zeta(s)\leq 0,
\intertext{et}
\label{lemme100eq190b}
e<\frac{1}{R_a}
\Longrightarrow 
 \forall s\in]0,l],\quad
&\zeta(s)<0
\end{align}
\end{subequations}
Cela est vrai en particulier pour $s=l$, d'où l'on tire 
\eqref{lemme100eq30c}
et 
\eqref{lemme100eq30d}.

\item
Concluons en montrant \eqref{lemme100eq30e}.
Remarquons que, d'après 
\eqref{lemme100eq150a},
\eqref{lemme100eq180a} et 
\eqref{lemme100eq190}
 et la définition \eqref{lemme100eq50} de $\zeta$, 
on a 
\begin{equation*}
0=\zeta_0=-\left(\widehat x\left(l\right) -x\left(l\right)\right) a
+
\left(\widehat y\left(l\right)-  y\left(l\right) \right) b
\end{equation*}
avec 
\begin{align*}
& a >0, \quad b >0,\\
&-\left(\widehat x\left(l\right) -x\left(l\right)\right) \leq 0,\quad 
  \widehat y\left(l\right)-  y\left(l\right)  \leq 0,
\end{align*}
et donc 
\begin{equation}
\label{lemme100eq200}
\widehat x\left(l\right) -x\left(l\right)=
\widehat y\left(l\right) -y\left(l\right)=0.
\end{equation}
Or on a vu, d'après 
\eqref{lemme100eq150aante} et 
\eqref{lemme100eq180aante},
 que les fonctions $\widehat x-x$ et $\widehat y-y$ étaient monotones sur $[0,l]$.
D'après \eqref{lemme100eq200},
on a donc :
\begin{equation*}
 \forall s\in[0,l],\quad
\widehat x(s)=x(s),\quad
\widehat y(s)=y(s).
\end{equation*}
On a donc
\begin{equation*}
 \forall s\in[0,l],\quad
\cos(\phi(s))=\cos(\theta(s)),\quad
\sin(\phi(s))=\sin(\theta(s)),
\end{equation*}
et, d'après \eqref{eq20} et \eqref{eq160}, on a donc
\begin{equation}
\label{lemme100eq201}
 \forall s\in[0,l],\quad
\phi(s)=\theta(s).
\end{equation}
En particulier $\Omega=\phi(l)=\theta(l )$. D'après 
\eqref{eq150} et \eqref{lemme100eq40a}, on a donc 
\begin{equation*}
\forall s\in [l,M],\quad
\Omega=\theta(M)\geq \theta(s) \geq \theta(l) =\phi(l)=\Omega,
\end{equation*}
Ainsi, 
\begin{equation}
\label{lemme100eq220}
\forall s\in [l,M],\quad
\phi(s)=\theta(s).
\end{equation}
\eqref{lemme100eq201} et \eqref{lemme100eq220}
impliquent que $\phi$ et $\theta$ coïncident sur $[0,M]$ et il en est de même pour $x$ et $\widehat  x$ et $y$ et $\widehat  y$.
En particulier en $s=M$ où $x$ et $y$ valent $x_B$ et $y_B$. Les deux courbes $X $ et $Z$ finissent donc au même point $B$ et donc $L=M$
et $X=Z$. On en déduit alors $e=1/R_a$.
\end{enumerate}
\ifcase \cras
\end{proof}
\or
\end{preuve}
\fi

\ifcase \cras
\begin{remark}
\or
\begin{remarque}
\fi
L'hypothèse \eqref{eq20bis} est fondamentale, via les inégalités \eqref{lemme100eq180}
et \eqref{lemme100eq190},
pour montrer le signe constant de $\zeta$.
\ifcase \cras
\end{remark}
\or
\end{remarque}
\fi

\ifcase \cras
\begin{proof}[Démonstration du théorème \ref{existenceunicitetheoprop01}]\
\or
\begin{preuve}[ du théorème \ref{existenceunicitetheoprop01}]\
\fi

\label{tugudu10}

Notons tout d'abord que $\inf_{Z\in \mathcal{E}}\vnorm[L^{\infty}(0,L(Z))]{\vnorm{Z''}}$ existe, puisque $\mathcal{E}$ est non vide et que 
$\vnorm[L^{\infty}(0,L(Z))]{\vnorm{Z''}}$ est toujours positif.

\begin{itemize}
\item[$\bullet$]

Premier cas :  l'hypothèse \eqref{eq20bis} est vraie. On peut alors utiliser directement la proposition 
\ref{lemme100}.

\begin{enumerate}
\item

Rappelons la définition \eqref{pdueq01} de $g$.
On a exhibé dans le lemme \ref{uniquecercledroitelem01} une fonction $X=\mathcal{J}(O,A,B)$ de $\mathcal{E}$ vérifiant $g(X)=1/R_a$ où $R_a=R_a(O,A,X)$. On a donc
\begin{equation}
\label{existenceunicitetheoprop01eq01}
\inf_{Z\in \mathcal{E}} g(Z)\leq \frac{1}{R_a}.
\end{equation}
\item
Montrons tout d'abord que 
\begin{equation}
\label{existenceunicitetheoprop01eq10}
\inf_{Z\in \mathcal{E}} g(Z)=\frac{1}{R_a}.
\end{equation}
Si c'était faux, d'après \eqref{existenceunicitetheoprop01eq01}, on 
aurait 
\begin{equation*}
\inf_{Z\in \mathcal{E}} g(Z)< \frac{1}{R_a},
\end{equation*}
et il existe donc une courbe de $\mathcal{E}$ notée $Z$ telle que 
\begin{equation}
\label{existenceunicitetheoprop01eq20}
 g(Z)< \frac{1}{R_a}.
\end{equation}
D'après la proposition  
\ref{lemme100}, on a donc $\zeta_0 < 0$.
Il existe donc un point de la courbe $Z$ dans le demi-plan ouvert défini par la tangente à la courbe $X$, notée $\mathcal{D}$
au point d'abscisse curviligne $l$, du côté opposé à la normale. 
\begin{figure}[h]    
\psfrag{a}{$Z(s)$}
\psfrag{A}{$A$}
\psfrag{B}{$B$}
\psfrag{a}{$\vec \alpha$}
\psfrag{b}{$\vec \beta$}
\psfrag{D}{$\mathcal{D}$}
\psfrag{Z}{$Z$}
\psfrag{M}{$Z\left(l\right)$}
\psfrag{X}{Z}
\begin{center} 
\epsfig{file=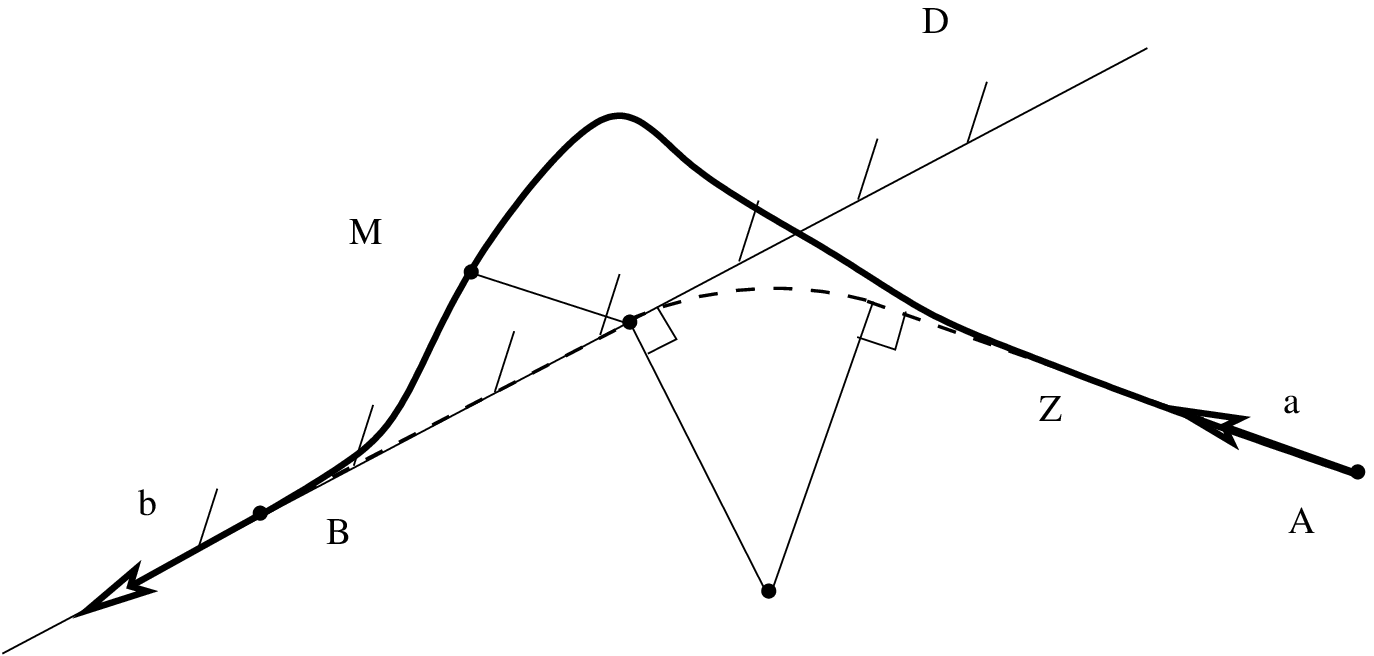, width=11 cm}  
\end{center} 
\caption{\label{point_autre_cote_normal}
Il existe donc un point de la courbe $Z\left(l\right)$ de $Z$ dans le demi-plan ouvert $\Pi$ défini par la tangente à la courbe $X$, notée $\mathcal{D}$
au point d'abscisse curviligne $l$, du côté opposé à la normale. 
\ifcase \cras \or \textit{%
There exists a point of curve in the open half-plane defined by the line
 $\mathcal{D}$, opposite to the 
 outer-pointing normal.}\fi}
\end{figure}
Voir figure \ref{point_autre_cote_normal}.
Cette droite $\mathcal{D}$ est aussi la droite passant par $B$ et portée par $\vec \beta$. 
Puisque la courbe $Z$ est dans $\mathcal{E}$, d'après le lemme \ref{lemmeconvexite},
la courbe $Z$ est incluse dans le demi-plan délimitée par la droite tangente à la courbe au point $B$, 
du côté la normale extérieure à la courbe en $B$, donc de l'autre côté du demi-plan $\Pi$  délimité par $\mathcal{D}$.
Il y a donc contradiction, ce qui achève la preuve de ce point et on a donc 
\begin{equation*}
\inf_{Z\in \mathcal{E}}g(Z)=g(X)=\min_{Z\in \mathcal{E}}g(Z)=\frac{1}{R_a}.
\end{equation*}

\item
Montrons enfin l'unicité de la courbe $X$ vérifiant 
\eqref{eq200} en montrant que cette courbe est la courbe $\mathcal{J}(O,A,B)$  donnée dans le lemme \ref{uniquecercledroitedef01}.
On utilise là encore la proposition  \ref{lemme100}.
Supposons qu'il existe une autre courbe $Z$ de $\mathcal{E}$ telle que 
\begin{equation}
\label{existenceunicitetheoprop01eq21}
g(Z)=\frac{1}{R_a},
\end{equation}
où $R_a=R_a(O,A,B)$. 
\begin{figure}[h]    
\psfrag{a}{$Z(s)$}
\psfrag{A}{$A$}
\psfrag{B}{$B$}
\psfrag{a}{$\vec \alpha$}
\psfrag{b}{$\vec \beta$}
\psfrag{D}{$\mathcal{D}$}
\psfrag{Z}{$Z$}
\psfrag{M}[][l]{$Z(\left(l\right)$}
\psfrag{X}{Z}
\begin{center} 
\epsfig{file=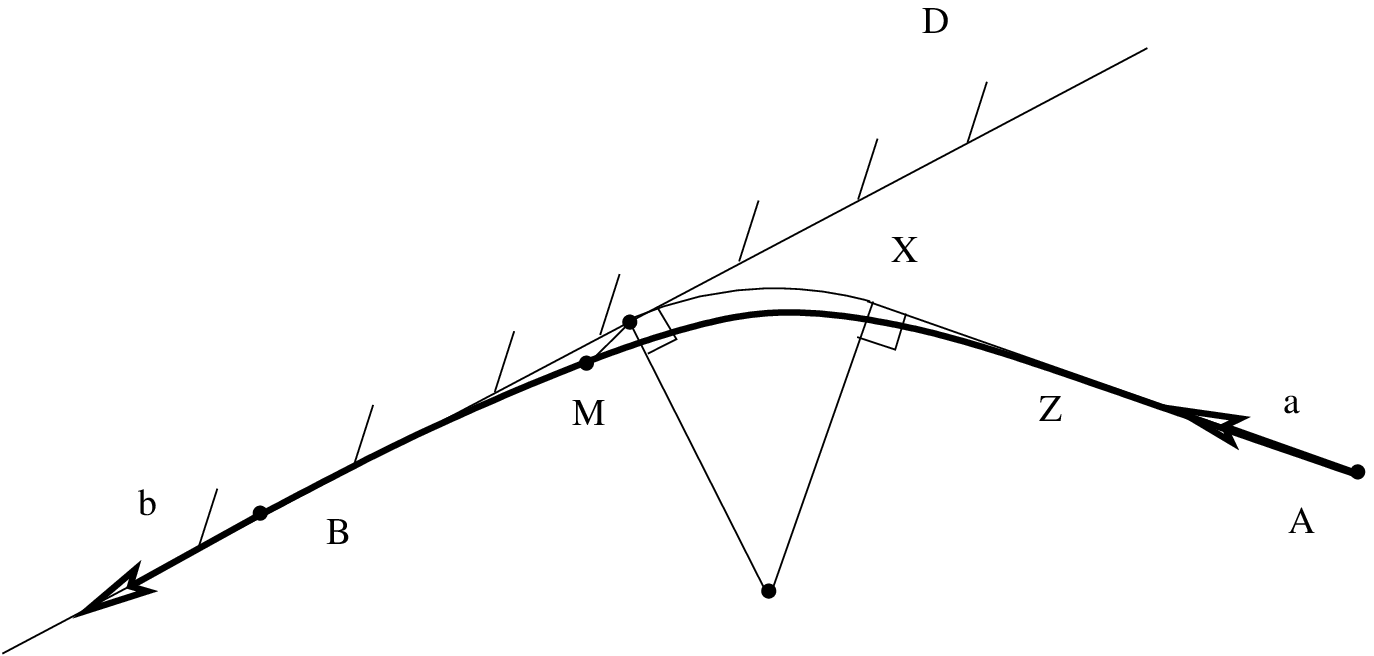, width=11 cm}  
\end{center} 
\caption{\label{point_autre_cote_normal_bis}La courbe $Z$ est toujours du côté de la normale extérieure.\ifcase \cras \or /\textit{The curve $Z$ is ever in the side of the outer-pointing normal.}\fi}
\end{figure}
Voir figure \ref{point_autre_cote_normal_bis}.
D'après le lemme \ref{lemmeconvexite} appliqué à cette courbe $Z$ et à la droite $\mathcal{D}$, la tangente à la courbe $Z$ au point $B$, on a donc
 avec les notations de la proposition 
\ref{lemme100}, $\zeta_0\geq 0$.
D'après cette  même proposition, on a $\zeta_0\leq 0$ et donc $\zeta_0=0$ et de nouveau d'après  la proposition
\ref{lemme100}, on a $Z=X$.

Le théorème \ref{existenceunicitetheoprop01} est donc vrai sous l'hypothèse \eqref{eq20bis}.
\end{enumerate}

\item[$\bullet$]

Second cas : 
l'hypothèse \eqref{eq20bis} n'est plus  vraie.
Nous allons décomposer le problème en deux sous-problèmes et à chacun d'eux, nous appliquerons 
le théorème \ref{existenceunicitetheoprop01}, sous l'hypothèse \eqref{eq20bis}.

\begin{enumerate}

\item
Commençons par montrer un résultat similaire à la proposition 
\ref{lemme100}. 
Supposons donc qu'il existe une courbe  $\mathcal{E}$ notée $Z$, de longueur $M$,  telle que, en notant $e=g(Z)$, 
\begin{equation}
\label{existenceunicitetheoprop01eq22}
e \leq \frac{1}{R_a}.
\end{equation}
Comme dans la proposition \ref{lemme100},
notons $\theta$ l'angle associé à cette courbe.
D'après le lemme \ref{lemmeconvexite}, la courbe est à la fois dans le demi-plan ouvert limité par $(OA)$ et contenant $B$
et le demi-plan ouvert limité par $(OB)$ et contenant $A$.
Considérons $s_1\in ]0,M[$, dont l'existence est assurée d'après les hypothèses sur $Z$,  tel que  
\begin{equation}
\label{existenceunicitetheoprop01eq30}
\theta(s_1)=\frac{\Omega}{2}.
\end{equation}
Considérons $Z_1$, la restriction de $Z$ à $[0,s_1]$ et $Z_2$ la restriction de $Z$ à $[s_1,M]$ et posons 
\begin{equation}
\label{existenceunicitetheoprop01eq40}
C=Z(s_1).
\end{equation}
\begin{figure}[h]    
\psfrag{A}{$A$}
\psfrag{B}{$B$}
\psfrag{O}{$O$}
\psfrag{C}{$C$}
\psfrag{D}{$D$}
\psfrag{E}{$E$}
\psfrag{F}{$F$}
\psfrag{Z1}{$Z_1$}
\psfrag{Z2}{$Z_2$}
\psfrag{d}{$Z'(s_1)$}
\psfrag{a}{$\vec \alpha$}
\psfrag{b}{$\vec \beta$}
\begin{center} 
\epsfig{file=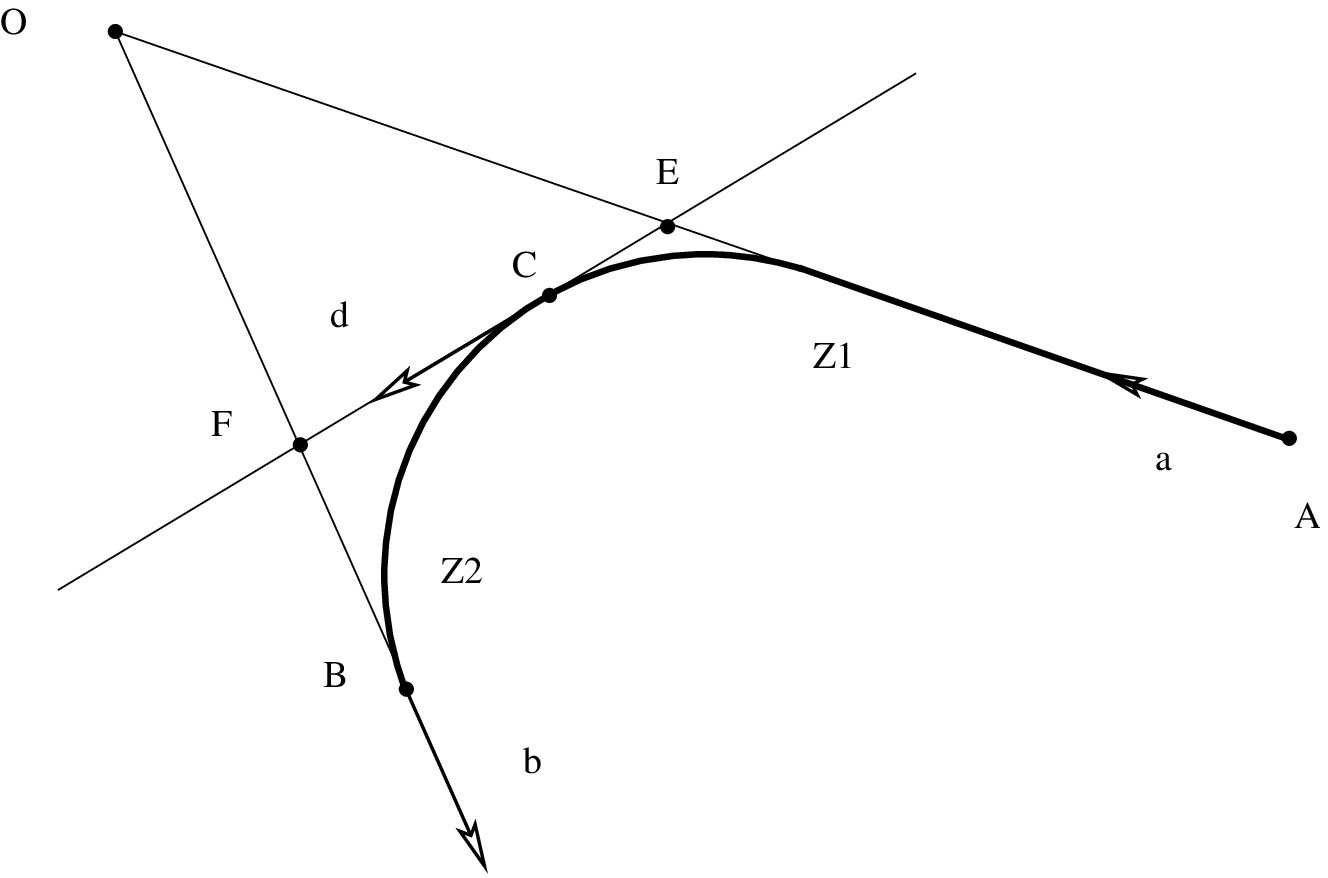, width=9 cm} 
\end{center}
\caption{\label{pointD_fig01}Les deux courbes $Z_1$ et $Z_2$. \ifcase \cras \or /\textit{{The both curves $Z_1$ and $Z_2$.}}\fi}
\end{figure}
La situation est représentée en figure 
\ref{pointD_fig01}. La courbe $Z_1$ relie donc les deux points $A$ et $C$ tandis que la courbe
$Z_2$ relie $C$ et $B$ et, d'après le lemme \ref{lemmeA}, les deux droites respectives passant par $A$, portée par $\vec \alpha$
et passant par $C$ portée par $Z'(s_1)$ se coupent en point $E$, distinct de $C$ et de $A$ avec $\alpha=\overrightarrow{AE}/{AE}$ et 
 $Z'(s_1)=\overrightarrow{EC}/{EC}$. 
De même, les deux droites respectives passant par $B$, portée par $\vec \beta$
et passant par $C$ portée par $Z'(s_1)$ se coupent en point $F$, distinct de $C$ et de $B$ avec $Z'(s_1)=\overrightarrow{CF}/{CF}$ et $\beta=
\overrightarrow{FB}/{FB}$.
Remarquons que, grâce à \eqref{existenceunicitetheoprop01eq22}, 
\begin{equation}
\label{sauvetot1}
\vnorm[L^{\infty}(0,s_1)]{\vnorm{Z_1''}}
\leq 
\vnorm[L^{\infty}(0,M)]{\vnorm{Z''}}=e\leq \frac{1}{R_a},
\end{equation}
et donc 
\begin{equation}
\label{pointD_eq01}
\vnorm[L^{\infty}(0,s_1)]{\vnorm{Z_1''}}\leq \frac{1}{R_a}.
\end{equation}
De même, 
\begin{equation}
\label{pointD_eq01b}
e<\frac{1}{R_a}  \Longrightarrow \vnorm[L^{\infty}(0,s_1)]{\vnorm{Z_1''}}< \frac{1}{R_a}.
\end{equation}
On peut donc appliquer le théorème \ref{existenceunicitetheoprop01}
 sous l'hypothèse \eqref{eq20bis}
aux trois points $A$, $E$ et $C$ puisque, d'après \eqref{existenceunicitetheoprop01eq30}
\begin{equation}
\label{pointD_eq07}
\left(\widehat{\vec \alpha,\frac{\overrightarrow{EC}}{{EC}}}\right)=
\frac{\Omega}{2}\in]0,\pi/2[.
\end{equation}
Il existe donc une unique courbe à courbure positive
$\widetilde Z_1=\mathcal{J}(E,A,C)$ passant par $A$ et $C$, dont les tangentes sont portées en ces points par $\vec\alpha$ et $Z'(s_1)$,
qui minimise $ \vnorm[L^{\infty}(0,L(z_1))]{\vnorm{z_1''}}$ pour $z_1$ décrivant 
l'ensemble ${\mathcal{E}}_1$ des courbes à courbure positive
passant par $A$ et $C$, dont les tangentes sont portées en ces points par $\vec\alpha$ et $Z'(s_1)$. 
D'après le théorème \ref{existenceunicitetheoprop01}, cette courbe est la réunion d'un arc de cercle de rayon $R_1$, 
d'angle $\Omega/2$ et d'un segment de droite de longueur $d_1\geq 0$. 
Puisque $Z_1$ est dans ${\mathcal{E}}_1$, 
on a donc 
\begin{equation}
\label{sauvetot2}
 \vnorm[L^{\infty}(0,s_1)]{\vnorm{Z_1''}}\geq   \vnorm[L^{\infty}(0,L(\widetilde Z_1))]{\vnorm{\widetilde Z_1}}  =\frac{1}{R_1}.
\end{equation}
ce qui implique, d'après \eqref{pointD_eq01} et \eqref{pointD_eq01b}
\begin{subequations}
\label{pointD_eq10}
\begin{align}
\label{pointD_eq10a}
& R_a\leq R_1,\\
\label{pointD_eq10b}
& e<\frac{1}{R_a} \Longrightarrow R_a<R_1.
\end{align}
\end{subequations}
Si on raisonne sur $Z_2$, on montre donc de même 
qu'il existe une unique courbe à courbure positive
$\widetilde Z_2=\mathcal{J}(F,C,B)$ passant par $C$ et $B$, dont les tangentes sont portées en ces points par  $Z'(s_1)$ et $\vec\beta$ ; cette courbe est la réunion d'un arc de cercle de rayon $R_2$, 
d'angle $\Omega/2$ et d'un segment de droite de longueur $d_2\geq 0$, avec 
\begin{subequations}
\label{pointD_eq20}
\begin{align}
\label{pointD_eq20a}
& R_a\leq R_2,\\
\label{pointD_eq20b}
&  e<\frac{1}{R_a} \Longrightarrow R_a<R_2.
\end{align}
\end{subequations}
Les deux courbes $\widetilde Z_1$ et $\widetilde Z_2$ peuvent être constituées respectivement et dans cet ordre : 
soit 
d'un arc de cercle et d'un segment de droite soit
d'un segment de droite et d'un arc de cercle.
Afin de couvrir tous les cas possibles, on considère donc l'unique courbe de $\mathcal{E}$, notée $\widetilde Z$, formée 
d'un segment de droite de longueur $d_1\geq 0$,
puis 
d'un arc de cercle de rayon $R_1$, 
d'angle $\Omega/2$,
puis 
d'un segment de droite de longueur $d_2\geq 0$,
puis 
d'un arc de cercle de rayon $R_2$, 
d'angle $\Omega/2$
enfin,
d'un 
segment de droite de longueur $d_3\geq 0$.
Le point $C$ est donc sur le segment de droite de longueur $d_2\geq 0$. 
On pose de façon analogue à 
\eqref{lemme100eq50}
\begin{equation}
\label{lemme100eq51}
\forall s\in [0,L(\widetilde Z)],\quad
\zeta(s)=
-\left(\widetilde  x\left(s\right) -x_B\right)\sin \Omega
+
\left(\widetilde   y\left(s\right)-  y_B \right)\cos\Omega,
\end{equation}
où $\widetilde Z(s)=(\widetilde  x\left(s\right),\widetilde   y\left(s\right))$.
On pose enfin
\begin{equation}
\label{pointD_eq30}
\zeta_0= \zeta(\widehat L),
\end{equation}
où $\widehat L$ est la somme des longueurs des deux premiers segments de droite et des deux arcs de cercle.
Comme dans la preuve de la proposition 
\ref{lemme100} on peut montrer que
\begin{subequations}
\label{pointD_eq32}
\begin{align}
\label{pointD_eq32a}
&\zeta_0 \leq 0,\\
\label{pointD_eq32b}
& \left(R_1 >R_a \text{ et } R_2>R_a\right) \Longrightarrow \zeta_0 < 0,\\
\label{pointD_eq32c}
&  \zeta_0=0 \Longrightarrow  \left(R_1=R_a \text{ et } R_2=R_a,\quad Z=\mathcal{J}(O,A,B)\right).
\end{align}
\end{subequations}
Le fait que $\zeta_0\leq 0$ signifie que la courbe $\widetilde Z$ a des points de l'autre côté (par rapport à la normale) de la tangente en $B$  la courbe.

Pour démontrer \eqref{pointD_eq32}, on procède comme suit :

On peut supposer, sans perte de généralité, comme dans la preuve de la proposition 
\ref{existenceunicitetheoprop01}, 
 que $OA\leq OB$, de sorte 
que la courbe $X=\mathcal{J}(O,A,B)\in \mathcal{E}$, de longueur $L$ décrite dans le lemme \ref{uniquecercledroitelem01}, commence d'abord par un arc de cercle de rayon $R_a=R_a(O,A,B)$. 
On pose 
\begin{equation}
\label{preuveD_eq16}
l=R_a\Omega.
\end{equation}
En utilisant les relations \eqref{eqderphinbbbb}, on obtient donc les expressions suivantes de $(x(s),y(s))=X(s)$ dans 
le repère orthonormé $\left(A, \vec {\alpha}, \vec {k}\right)$ :
\begin{subequations}
\label{preuveD_eq01tot}
\begin{align}
\label{preuveD_eq01a}
\forall s\in [0,l],\quad
x(s)&=R_a\sin\left(\frac{s}{R_a}\right),\\
\label{preuveD_eq01b}
y(s)&=R_a\left(1-\cos\left(\frac{s}{R_a}\right)\right),\\
\label{preuveD_eq01c}
\forall s\in [l,L],\quad
x(s)&=R_a\sin\left(\Omega\right)+(s-l)\cos\Omega,\\
\label{preuveD_eq01d}
y(s)&=R_a\left(1-\cos\left(\Omega\right)\right)+(s-l)\sin\Omega.
\end{align}
\end{subequations}
On a donc en particulier les coordonnées du point $B$ : 
\begin{subequations}
\label{preuveD_eq10tot}
\begin{align}
\label{preuveD_eq10a}
&x_B=R_a\sin\left(\Omega\right)+(L-l)\cos\Omega,\\
\label{preuveD_eq10b}
&y_B=R_a\left(1-\cos\left(\Omega\right)\right)+(L-l)\sin\Omega.
\end{align}
\end{subequations}
On pose 
\begin{subequations}
\label{preuveD_eq15tot}
\begin{align}
\label{preuveD_eq15a}
&l_1=\frac{R_1 \Omega}{2},\\
\label{preuveD_eq15b}
&l_2=\frac{R_2 \Omega}{2}.
\end{align}
\end{subequations}
On démontre ensuite de la même façon que l'on a les 
expressions suivantes pour les diverses parties, rectilignes ou circulaires, de la courbe 
$\widetilde Z$  :
\begin{subequations}
\label{preuveD_eq20tot}
\begin{align}
\label{preuveD_eq20a}
\forall s\in[0,d_1],\quad
&\widetilde  x\left(s\right)=s,\\
\label{preuveD_eq20b}
&\widetilde  y\left(s\right)=0,\\
\label{preuveD_eq20c}
\forall s\in[d_1,d_1+l_1],\quad
&\widetilde  x\left(s\right)
=
d_1+R_1\sin\left(\frac{s-d_1}{R_1}\right),\\
\label{preuveD_eq20d}
&\widetilde  y\left(s\right)
=
R_1\left(1-\cos\left(\frac{s-d_1}{R_1}\right)\right),\\
\label{preuveD_eq20e}
\forall s\in[d_1+l_1,d_1+l_1+d_2],\quad
&\widetilde  x\left(s\right)=
d_1+R_1\sin\left(\frac{\Omega}{2}\right)+(s-(d_1+l_1))\cos\left(\frac{\Omega}{2}\right),\\
\label{preuveD_eq20f}
&\widetilde  y\left(s\right)=
R_1\left(1-\cos\left(\frac{\Omega}{2}\right)\right)+(s-(d_1+l_1))\sin\left(\frac{\Omega}{2}\right),
\end{align}
{et enfin, pour tout $s\in[d_1+l_1+d_2,d_1+l_1+d_2+l_2]$,}
\begin{align}
\label{preuveD_eq20g}
&
\widetilde  x\left(s\right)=
d_1+R_1\sin\left(\frac{\Omega}{2}\right)+d_2 \cos\left(\frac{\Omega}{2}\right)+R_2\left(\sin\left(\frac{\Omega}{2}+\frac{s-(d_1+l_1+d_2)}{R_2}\right)-\sin\left(\frac{\Omega}{2}\right)\right),\\
\label{preuveD_eq20h}
&\widetilde  y\left(s\right)=
R_1\left(1-\cos\left(\frac{\Omega}{2}\right)\right)+d_2 \sin\left(\frac{\Omega}{2}\right)-R_2\left(\cos\left(\frac{\Omega}{2}+\frac{s-(d_1+l_1+d_2)}{R_2}\right)-\cos\left(\frac{\Omega}{2}\right)\right).
\end{align}
\end{subequations}
On déduit de ces deux dernières lignes :
\begin{subequations}
\label{preuveD_eq30tot}
\begin{align}
\label{preuveD_eq30a}
&\widehat L=d_1+l_1+d_2+l_2,\\
\label{preuveD_eq30b}
&
\widetilde  x\left(\widehat L\right)=
d_1+R_1\sin\left(\frac{\Omega}{2}\right)+d_2 \cos\left(\frac{\Omega}{2}\right)+R_2\left(\sin\left({\Omega}\right)-\sin\left(\frac{\Omega}{2}\right)\right),\\
\label{preuveD_eq30c}
&
\widetilde  y\left(\widehat L\right)=
R_1\left(1-\cos\left(\frac{\Omega}{2}\right)\right)+d_2 \sin\left(\frac{\Omega}{2}\right)-R_2\left(\cos\left({\Omega}\right)-\cos\left(\frac{\Omega}{2}\right)\right).
\end{align}
\end{subequations}
Ainsi, on 
a, grâce à \eqref{preuveD_eq10a} et \eqref{preuveD_eq30b},
\begin{align*}
\widetilde  x\left(\widetilde L\right) -x_B
&=
d_1+R_1\sin\left(\frac{\Omega}{2}\right)+d_2 \cos\left(\frac{\Omega}{2}\right)+R_2\left(\sin\left({\Omega}\right)-\sin\left(\frac{\Omega}{2}\right)\right)
-R_a\sin\left(\Omega\right)-(L-l)\cos\Omega,\\
&=
d_1+d_2\cos\left(\frac{\Omega}{2}\right)+(R_2-R_a)\sin\left(\Omega\right)+(R_1-R_2)\sin\left(\frac{\Omega}{2}\right) -(L-l)\cos\Omega.
\end{align*}
De même, 
\begin{equation*}
\widetilde  y\left(\widetilde L\right) -y_B
=d_2\sin\left(\frac{\Omega}{2}\right)+(R_1-R_a)+(R_a-R_2)\cos\left({\Omega}\right)+(R_2-R_1)\cos\left(\frac{\Omega}{2}\right)-(L-l)\sin\Omega.
\end{equation*}
On en déduit donc, grâce à \eqref{lemme100eq51} et \eqref{pointD_eq30} : 
\begin{multline*}
\zeta_0=
-
\left(d_1+d_2\cos\left(\frac{\Omega}{2}\right)+(R_2-R_a)\sin\left(\Omega\right)+(R_1-R_2)\sin\left(\frac{\Omega}{2}\right)\right)\sin \left(\Omega\right)+\\
\left(d_2\sin\left(\frac{\Omega}{2}\right)+(R_1-R_a)+(R_a-R_2)\cos\left({\Omega}\right)+(R_2-R_1)\cos\left(\frac{\Omega}{2}\right)\right)\cos \left(\Omega\right).
\end{multline*}
En posant 
\begin{subequations}
\label{preuveD_eq40tot}
\begin{align}
\label{preuveD_eq40a}
&a=\cos \left(\Omega\right)-\cos \left(\frac{\Omega}{2}\right)=-2\sin\left(\frac{3\Omega}{4}\right)\sin\left(\frac{\Omega}{4}\right),\\
\label{preuveD_eq40b}
&b=-1+\cos\left(\frac{\Omega}{2}\right),\\
\label{preuveD_eq40c}
&c=-\sin\left(\Omega\right),\\
\label{preuveD_eq40d}
&f=-\sin\left(\frac{\Omega}{2}\right),
\intertext{on a, après calculs,}
\label{preuveD_eq40e}
&\zeta_0=a(R_1-R_a)+b(R_2-R_a)+cd_1+fd_2.
\end{align}
\end{subequations}
Puisque l'hypothèse \eqref{eq20bis} n'est plus valable,
on a d'après \eqref{eq20} 
\begin{equation}
\label{eq20ter}
\Omega \in [\pi/2,\pi[,
\end{equation}
et  \eqref{preuveD_eq40tot} implique 
\begin{equation}
\label{preuveD_eq50}
a<0,\quad
b<0,\quad
c<0,\quad
f<0.
\end{equation}
De \eqref{preuveD_eq40tot}, 
\eqref{eq20ter}
et \eqref{preuveD_eq50}, on déduit donc
\eqref{pointD_eq32a} et \eqref{pointD_eq32b}.
De plus, si $\zeta_0=0$, on a donc $R_1=R_2=R_a$ et $d_1=d_2=0$. 
Donc, la courbe $\widetilde Z$ est constituée 
d'un arc de cercle de rayon $R_a$ et d'angle $\Omega$
et 
d'un 
segment de droite de longueur $d_3\geq 0$, qui est donc exactement $\mathcal{J}(O,A,B)$.
Puisque le point $C$ appartient à l'arc de cercle de rayon $R_a$ et d'angle $\Omega$, 
la courbe $\mathcal{J}(E,A,C)$ relie $A$ à $C$. De plus,   
d'après  
\eqref{pointD_eq01} et \eqref{sauvetot2}
\begin{equation*}
\frac{1}{R_a} \geq \vnorm[L^{\infty}(0,s_1)]{\vnorm{Z_1''}} \geq   \vnorm[L^{\infty}(0,L(\widetilde Z_1))]{\vnorm{\widetilde Z_1}}  =\frac{1}{R_1}=\frac{1}{R_a}.
\end{equation*}
et donc 
\begin{equation*}
\vnorm[L^{\infty}(0,s_1)]{\vnorm{Z_1''}}=\frac{1}{R_a}=\frac{1}{R_1}.
\end{equation*}
Ainsi, d'après le théorème \eqref{existenceunicitetheoprop01} 
grâce à \eqref{pointD_eq07},
par unicité, on a $Z_1=\mathcal{J}(E,A,C)$. De même, 
$Z_2=\mathcal{J}(F,C,B)$ et donc $Z=\mathcal{J}(O,A,B)$, ce qui achève  la preuve de \eqref{pointD_eq32c}.

\item

L'inégalité \eqref{existenceunicitetheoprop01eq01} est toujours vraie.
Commençons tout d'abord par montrer 
\eqref{existenceunicitetheoprop01eq10}
et supposons donc qu'il existe une courbe  $\mathcal{E}$ notée $Z$, de longueur $M$,  telle que 
\eqref{existenceunicitetheoprop01eq20} soit vraie.
D'après  
\eqref{pointD_eq10b},
\eqref{pointD_eq20b}
et 
\eqref{pointD_eq32b},
 on a $\zeta_0<0$, ce qui contredit, grâce au lemme
 \ref{lemmeconvexite}, le fait que $\widetilde Z$, construite précédemment à partir de $Z$,  est dans $\mathcal{E}$.

\item
 Montrons maintenant l'unicité et supposons  qu'il existe une autre courbe $Z$ de $\mathcal{E}$ 
 vérifiant \eqref{existenceunicitetheoprop01eq21}.
D'après  \eqref{pointD_eq32a}, on a $\zeta_0\leq 0$. Le  lemme
 \ref{lemmeconvexite}, implique  que puisque $\widetilde Z$, construite précédemment à partir de $Z$,  est dans $\mathcal{E}$  et on a  $\zeta_0\geq 0$. On a donc $\zeta_0=0$ et 
d'après le premier point et \eqref{pointD_eq32c}, on a $Z=\mathcal{J}(O,A,B)$.
\end{enumerate}

\end{itemize}
\ifcase \cras
\end{proof}
\or
\end{preuve}
\fi

\ifcase \cras
\begin{remark}
\or
\begin{remarque}
\fi
\label{remarquexcool}
Dans la démonstration  de la proposition \ref{lemme100}, sous l'hypothèse \eqref{eq20ter}, 
on a aussi, comme sous-produit, le fait que la réunion de deux arcs de cercles de même angle $\Omega/2$, 
et de trois segments de droite est dans $\mathcal{E}$ ssi cette réunion est réduite à 
un arc de cercle d'angle $\Omega$ et d'un segment de droite et donc retrouvé un cas particulier du cas discret optimal 
avec $p=4$  présenté dans 
\cite[section \numerosectioncasdiscretgeneral]{piece6_optimale_JB_2019_X4}.
Si on ne fait plus l'hypothèse \eqref{eq20ter},  le raisonnement n'est plus valable car 
\eqref{preuveD_eq50} n'est plus vrai.
En revanche, il est intéressant de constater que si l'on choisit $R_1=R_2$ et $d_2=0$, on remplace
\eqref{preuveD_eq40tot} par 
\begin{equation*}
a+b=\cos \left(\Omega\right)-1, \quad
c=-\sin\left(\Omega\right),
\end{equation*}
et on a 
\begin{equation*}
\zeta_0=(a+b)(R-R_a)+cd_1.
\end{equation*}
Dans ce cas, \eqref{eq20} implique
\begin{equation*}
a+b<0,\quad
c<0
\end{equation*}
Ainsi, si $\zeta_0=0$, on a donc $R_1=R_2=R_a$ et $d_1=0$. 
On conclut donc de nouveau que  $Z=\mathcal{J}(O,A,B)$.
Autrement dit,  la réunion d'un arc de cercle et de deux segments 
 est dans $\mathcal{E}$ ssi cette réunion est réduite à 
un arc de cercle d'angle $\Omega$ et d'un segment de droite et a donc retrouvé le cas discret optimal 
avec $p=2$  présenté dans 
\cite[section \numerosectioncasdiscretgeneral]{piece6_optimale_JB_2019_X4}.
\ifcase \cras
\end{remark}
\or
\end{remarque}
\fi

\ifcase \cras

\printbibliography

\fi

\end{document}

%% file: ajouts.tex

\newcommand{\ajoutA}{%
\uwave{D'après \eqref{eq100d} et \eqref{eq130}, on a}
\begin{equation*}
\uwave{\text{$\displaystyle{\forall s\in [0,L],\quad \phi(s)=\int_0^s c(t) dt,}$}}
\end{equation*}
\uwave{ et donc l'angle $\phi$ est bien continu.
Désormais,  dans}
\ifcase \cras
\uwave{tout cet article,}
\or
\uwave{toute cette Note,}
\fi
\uwave{on adoptera systématiquement la convention suivante :
pour tout couple de vecteurs} $\left(\vec x,\vec y\right)$ de $\left(\Er^2\right)^2$,
\ifcase \rar
\begin{subequations}
\label{xi}
\begin{align}
\label{xia}
&\text{\uwave{$\left( \widehat{\vec x,\vec y} \right)$ désigne la détermination principale de l'angle des deux vecteurs, appartenant à $]-\pi,\pi]$}}
\intertext{et}
\label{xib}
&\text{\uwave{chaque égalité d'angle est vraie strictement et non modulo-$2\pi$.}}
\end{align}
\end{subequations}
\or
\begin{align*}
&
\null
\hfill
\text{\uwave{$\left( \widehat{\vec x,\vec y} \right)$ désigne la détermination principale de l'angle des deux vecteurs, appartenant à $]-\pi,\pi]$}}
\hfill
\eqref{xia}
\intertext{et}
&
\null
\hfill
\text{\uwave{chaque égalité d'angle est vraie strictement et non modulo-$2\pi$.}}
\hfill
\eqref{xib}
\end{align*}
\fi
}

\newcommand{\ajoutB}{%
\sout{positive}
\uwave{positivity of the}
}

\newcommand{\ajoutC}{%
\uwave{points}
}

\newcommand{\ajoutD}{%
\uwave{Enfin, notons que ce problème est différent du problème de Dubins et ne peut pas se poser classiquement comme un problème de
contrôle optimal, le coût n'étant ni un coût intégral ni un coût final. Ce coût en norme $L^{\infty}$  du contrôle
(le contrôle est la courbure dans ce point de vue) est connu pour ne pas donner de problème de contrôle
optimal qui relève d'une théorie établie. 
Voir par exemple }\cite[Transparents 19 et 20]{piece6_optimale_JB_CRAS_2019_ICJ} \uwave{qui montre ce problème sous la forme d'un problème d'optimisation
en norme $L^{\infty}$.
Les travaux classiques de contrôle optimal, comme par exemple \cite{MR2062547,MR0220537,MR0166037,MR1354838},
ne traitent
pas les coûts $L^{\infty}$.
Par exemple, l'article \cite{MR2873265} étudie  la recherche d'une courbe confinée dont on cherche le minimum non pas du maximum de la courbure mais de sa norme $L^2$.
Dans}
\ifcase \cras
\uwave{cet article,}
\or
\uwave{cette Note,}
\fi
\uwave{nous proposerons donc une résolution de ce problème par des méthodes élémentaires.}
}

\newcommand{\ajoutF}{%
\sout{une détermination continue de}
}

\newcommand{\ajoutG}{%
\uwave{Les équations  \eqref{eq20} et \eqref{eq160} et la convention \eqref{xi}
impliquent aussi}
\ifcase \rar
\begin{equation}
\label{xibis}
\text{\uwave{$\forall s\in [0,L]$, $\quad $
$ \phi(s)\in [0,\pi[$.}}
\end{equation}
\or
\begin{equation*}
\null
\hfill
\text{\uwave{$\forall s\in [0,L]$, $\quad $
$ \phi(s)\in [0,\pi[$.}}
\hfill
\eqref{xibis}
\end{equation*}
\fi
\uwave{et, \textit{a posteriori}, on peut donc vérifier que 
\eqref{eq10totnew}
 et \eqref{eq120} respectent la convention \eqref{xi}.}
}

\newcommand{\ajoutHint}{\uwave{On se donne $O$, $A$ et $B$, trois points du plan deux à deux distincts puis deux vecteurs unitaires $\vec \alpha$ et $\vec \beta$
définis par \eqref{eq10tot}
et vérifiant 
\eqref{eq10totnew} et \eqref{eq20}.}}

\newcommand{\ajoutH}{%
\sout{Soient deux vecteurs unitaires $\vec \alpha$ et $\vec \beta$, vérifiant \eqref{eq20} et  $O$ donné par 
\eqref{defiEremeq01}.}
\ajoutHint
}

\newcommand{\ajoutHb}{%
\sout{Pour $O$, $A$ et $B$ trois points du plan deux à deux distincts et $\vec \alpha$ et $\vec \beta$ 
deux vecteurs vérifiant 
 \eqref{eq10tottot} et \eqref{eq20}, 
on}
\ajoutHint\  On 
}

\newcommand{\ajoutHc}{%
\sout{%
\eqref{eq10totnew},
\eqref{eq20},}}

\newcommand{\ajoutI}{%
\ifcase \rar
\ifcase \cras
\begin{remark}
\or
\begin{remarque}
\fi
\or
{\textit{Remarque \ref{sanshypodeconne}}}
\fi
\ifcase \rar
\label{sanshypodeconne}
\fi
\ifcase \rar
\uwave{Si on s'affranchit de l'hypothèse \eqref{eq20}, le problème est mal posé ; en effet, 
il est possible de construire une courbe vérifiant 
 \eqref{eq10tottot},  
\eqref{eq100tot},
\eqref{eq110},
\eqref{eq140} (ou \eqref{eq150}), avec un minimum de rayon de  courbure arbitrairement grand, comme le montre la figure 
\ref{codu0Rgrand} 
où le minimum de rayon de  courbure vaut $R$, choisi aussi grand que l'on veut.}
\or
\uwave{Si on s'affranchit de l'hypothèse \eqref{eq20}, le problème est mal posé ; en effet, 
il est possible de construire une courbe vérifiant 
 \eqref{eq10tottot},  
\eqref{eq100tot},
\eqref{eq110},
\eqref{eq140} (ou \eqref{eq150}), avec un minimum de rayon de  courbure arbitrairement grand, comme le montre la figure 
\ref{codu0Rgrandrar} 
où le minimum de rayon de  courbure vaut $R$, choisi aussi grand que l'on veut.}
\fi
\begin{figure}[h] 
\psfrag{R}{$R$}
\psfrag{A}{$A$}
\psfrag{B}{$B$}
\psfrag{a}{$\vec \alpha$}
\psfrag{b}{$\vec \beta$}
\begin{center} 
\epsfig{file=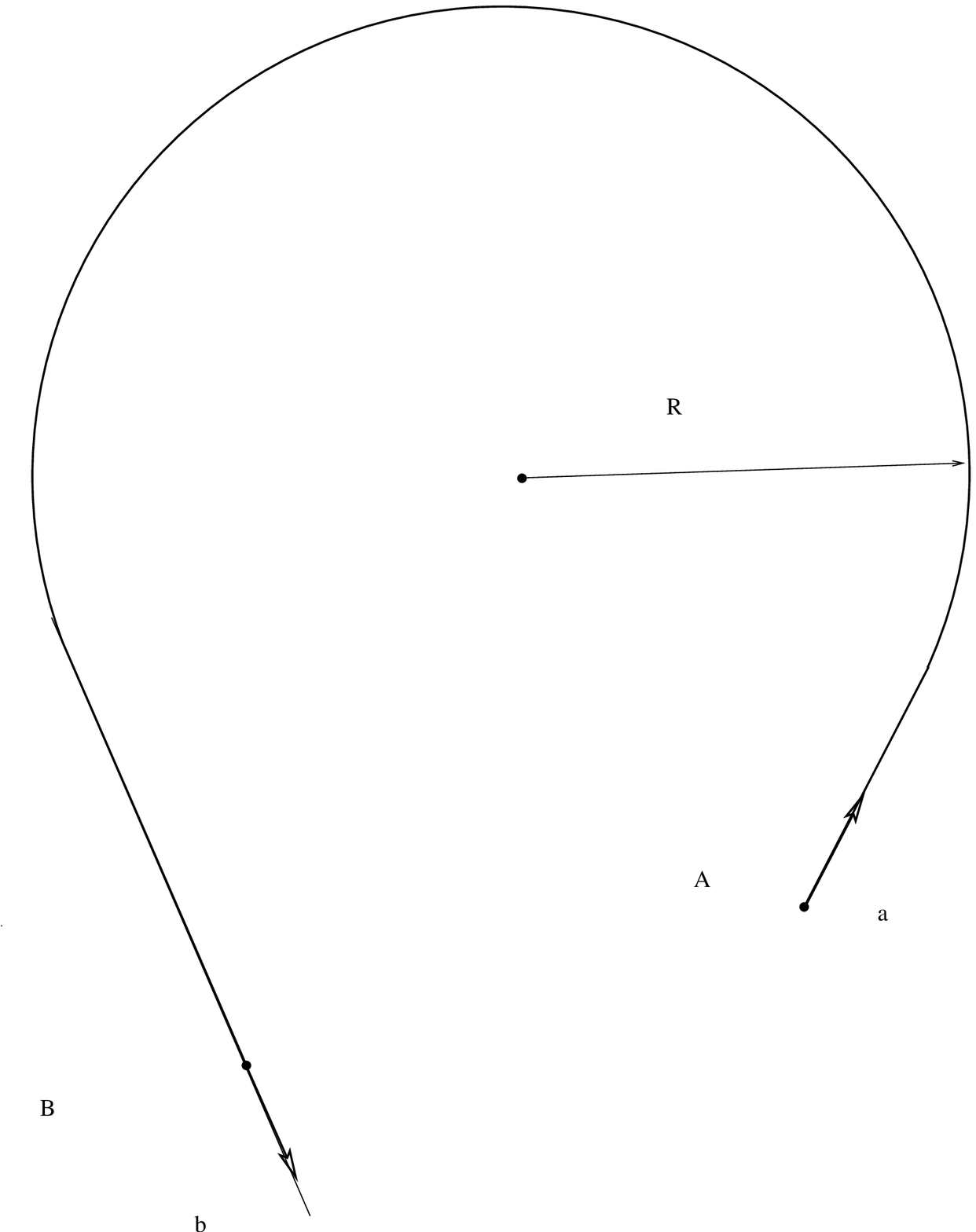, width=4  cm} 
\end{center} 
\ifcase \rar
\caption{\uwave{\label{codu0Rgrand}Une courbe de $\mathcal{E}$ formée de deux segments
et d'un arc de cercle de rayon $R$ arbitrairement grand.}\ifcase \cras \or /\textit{\uwave{A curve of $\mathcal{E}$  formed by two line segments
and an arc of a circle of arbitrarily large radius $R$.}}\fi}
\or
\caption{\uwave{\label{codu0Rgrandrar}Une courbe de $\mathcal{E}$ formée de deux segments
et d'un arc de cercle de rayon $R$ arbitrairement grand.}\ifcase \cras \or /\textit{\uwave{A curve of $\mathcal{E}$  formed by two line segments
and an arc of a circle of arbitrarily large radius $R$.}}\fi}
\fi
\end{figure}

\uwave{De même, si on s'affranchit de l'hypothèse \eqref{eq140} (ou \eqref{eq150} ou même \eqref{eq160}), le problème 
consistant à recherche une courbe vérifiant 
\eqref{eq10tottot}, 
\eqref{eq20}
\eqref{eq100tot} et 
\eqref{eq110},
est mal posé ; en effet, 
il est possible de construire une courbe de Dubins associée à un rayon $R$ aussi grand que l'on veut
mais cette courbe ne sera pas à courbure positive et \eqref{eq160} ne sera  pas vérifié.
Voir par exemple}
\ifcase \cras
les figures A.\ref{codu4} et  A.\ref{codu2sc}.
\or
\ifcase \versiondefinitive
\uwave{les figures A.\ref{codu4} et  A.\ref{codu2sc}.}  
\or
\uwave{\cite[\citerefhuit]{piece6_optimale_JB_arXiv_X6}.}
\fi
\fi
\ifcase \rar
\ifcase \cras
\end{remark}
\or
\end{remarque}
\fi
\or
\fi
}

\newcommand{\ajoutJJ}{%
\uwave{\eqref{xi},}
}

\newcommand{\ajoutK}{%
\sout{et l'angle $\Omega$ définis par
\eqref{eq10tottot},}
\uwave{définis par \eqref{eq10tot} et l'angle $\Omega$ défini par
\eqref{eq10totnew}
}%
}

\newcommand{\ajoutL}{%
\uwave{Le théorème \ref{existenceunicitetheoprop01} assurera que cette courbe est la courbe optimale. On pourra consulter }
\ifcase \cras
la section \ref{construceffct}.
\or
\uwave{\cite[\citerefneuf]{piece6_optimale_JB_arXiv_X6}.}
\fi
}

\newcommand{\ajoutM}{%
\sout{Trois}\uwave{Deux}
}

\newcommand{\ajoutN}{%
\sout{Une autre généralisation consisterait à s'affranchir de l'hypothèse \eqref{eq20}.}
}

\newcommand{\ajoutO}{%
\uwave{(voir la figure }\cite[figure 1f)]{enumeration_circuit_JB_2016}\uwave{)}
}

\newcommand{\ajoutOb}{%
\uwave{(données qui correspondent aux caractéristiques de la figure }\cite[figure 1f)]{enumeration_circuit_JB_2016}\uwave{)}
}

\newcommand{\ajoutP}{%
\uwave{Voir remarque \ref{sanshypodeconne}.}
}

\newcommand{\ajoutQ}{%
\sout{\eqref{eq150} implique}
\uwave{%
\eqref{eq10totnew},
\eqref{eq100d},
\eqref{eq100e},
\eqref{eq150}
et la convention \eqref{xi}  impliquent}
}

\newcommand{\ajoutR}{%
\sout{on n'a plus l'hypothèse \eqref{eq20bis}.}
\uwave{l'hypothèse \eqref{eq20bis} n'est plus valable.}%
}

\newcommand{\ajoutRb}{%
\sout{on a l'hypothèse \eqref{eq20bis}.}
\uwave{l'hypothèse \eqref{eq20bis} est valable.}%
}

\newcommand{\ajoutS}{%
\uwave{et de longueur minimale,}}

\newcommand{\ajoutSb}{%
\sout{tous les résultats d'unicité tombent}
les 
\ifcase \cras
\uwave{démonstrations}
\or
\uwave{preuves}
\fi
\uwave{des résultats d'unicité ne sont plus valables}}

\newcommand{\ajoutT}{%
\sout{le repère orthonormé $\left(A, \vec {\alpha}, \vec {k}\right)$}
\uwave{repère orthonormé direct $\left(A, \vec {\alpha}, \vec {k}\right)$}%
}

\newcommand{\ajoutU}{%
\uwave{direct}}

\newcommand{\ajoutV}{%
\sout{\eqref{eq10totnew},}
\uwave{\eqref{eq10tottot},}%
}

\newcommand{\ajoutW}{%
\sout{tels deux}%
\uwave{deux tels}%
}

\newcommand{\ajoutX}{%
\uwave{équivalent à \eqref{eq150}, }%
}

\newcommand{\ajoutY}{%
\uwave{Désormais, on }\sout{On}%
}

\newcommand{\ajoutZ}{%
\sout{$\phi $}\uwave{$\phi =\left( \widehat{\vec \alpha,X'} \right)$}%
}

\newcommand{\ajoutZb}{%
\sout{$\theta$}\uwave{$\theta =\left( \widehat{\vec \alpha,Z'} \right)$}%
}

\newcommand{\ajoutAA}{%
\sout{$\inf_{Z\in \mathcal{E}}\vnorm[L^{\infty}(0,L(X))]{\vnorm{X''}}$}
\uwave{$\inf_{Z\in \mathcal{E}}\vnorm[L^{\infty}(0,L(Z))]{\vnorm{Z''}}$}
}

\newcommand{\ajoutAAb}{%
\sout{$\vnorm[L^{\infty}(0,L(X))]{\vnorm{Z''}}$}
\uwave{$\vnorm[L^{\infty}(0,L(Z))]{\vnorm{Z''}}$}
}

\newcommand{\ajoutAB}{%
\sout{$ \vnorm[L^{\infty}(0,L(Z))]{\vnorm{z_1''}}$.}
\uwave{$ \vnorm[L^{\infty}(0,L(z1))]{\vnorm{z_1''}}$.}
}

\newcommand{\ajoutAC}{%
sous l'hypothèse \eqref{eq20bis}
}

\newcommand{\ajoutAD}{%
.}

\newcommand{\ajoutADb}{%
\sout{alors,  pour}
\uwave{Pour}
}

\newcommand{\ajoutAE}{%
d'après \eqref{lemme100eq20}
}

%% file: textelemmeconvexite.tex
\ifcase \cras
\begin{lemma}
\or
\ifcase \versiondefinitive
\begin{lemme}
\or
\begin{lemme}[\cite{piece6_optimale_JB_arXiv_X6}, Lemme A.1, Annexe A, page 14]
\fi
\fi
\label{lemmeconvexite}
Soit une courbe $X$ vérifiant 
\eqref{eq10tottot},
\eqref{eq20},
\eqref{eq100tot},
\eqref{eq110} et  
\eqref{eq150}\ajoutAD\
\ajoutADb\
tout  $s\in [0,L]$, la courbe est incluse dans le demi-plan délimité par la droite tangente à la courbe au point $X(s)$, 
du côté de $\vec N(s)$,  la normale extérieure à la courbe en $X(s)$. 
\ifcase \cras
\end{lemma}
\or
\end{lemme}
\fi